\documentclass{gtart_a}
\pdfoutput=1
\usepackage{amscd}
\usepackage{xypic}
\xyoption{all}

\title[Rounding corners of polygons and the embedded contact homology of $T^3$]{Rounding corners of polygons and the\\embedded contact homology of $T^3$}

\author{Michael Hutchings}
\givenname{Michael}
\surname{Hutchings}
\address{Department of Mathematics\\
University of California\\\newline
Berkeley CA 94720\\USA}
\email{hutching@math.berkeley.edu}
\urladdr{http://math.berkeley.edu/~hutching/}

\author{Michael Sullivan}
\givenname{Michael C}
\surname{Sullivan}
\address{Department of Mathematics and Statistics\\
University of Massachusetts\\\newline
Amhurst MA 01003-9305\\
USA}
\email{sullivan@math.umass.edu}
\urladdr{http://www.math.umass.edu/~sullivan/}

\volumenumber{10}
\issuenumber{}
\publicationyear{2006}
\papernumber{6}
\lognumber{0503}
\startpage{169}
\endpage{266}

\doi{10.2140/gt.2006.10.169}

\keyword{embedded contact homology}
\keyword{Floer homology}
\subject{primary}{msc2000}{57R58}
\subject{secondary}{msc2000}{57M27}

\received{5 October 2004}
\revised{27 December 2005}
\accepted{25 January 2006}
\published{26 March 2006}
\publishedonline{26 March 2006}
\proposed{Robion Kirby}
\seconded{Peter Ozsv\'ath, Tomasz Mrowka}

\def\bigmid{\,\,\vrule width .8pt height 12pt depth 5pt\,\,}
\let\frak\mathfrak

\def\cnewtheorem#1[#2]#3{\newtheorem{#1}{#3}[section]}
\let\widecheck\check

\newcommand{\mc}[1]{{\mathcal #1}}

\newtheorem{theorem}{Theorem}[section]
\cnewtheorem{proposition}[theorem]{Proposition}
\cnewtheorem{corollary}[theorem]{Corollary}
\cnewtheorem{lemma}[theorem]{Lemma}
\cnewtheorem{lemmadef}[theorem]{Lemma-Definition}
\cnewtheorem{sublemma}[theorem]{Sublemma}
\cnewtheorem{assumption}[theorem]{Assumption}
\cnewtheorem{conjecture}[theorem]{Conjecture}
\cnewtheorem{claim}[theorem]{Claim}

\theoremstyle{definition}
\cnewtheorem{definition}[theorem]{Definition}
\cnewtheorem{remark}[theorem]{Remark}
\cnewtheorem{remarks}[theorem]{Remarks}
\cnewtheorem{example}[theorem]{Example}
\cnewtheorem{convention}[theorem]{Convention}
\cnewtheorem{notation}[theorem]{Notation}
\cnewtheorem{openproblem}[theorem]{Open Problem}

\makeatletter
\let\c@proposition\c@theorem
\let\c@corollary\c@theorem
\let\c@lemma\c@theorem
\let\c@lemmadef\c@theorem
\let\c@sublemma\c@theorem
\let\c@assumption\c@theorem
\let\c@conjecture\c@theorem
\let\c@claim\c@theorem
\let\c@definition\c@theorem
\let\c@remark\c@theorem
\let\c@remarks\c@theorem
\let\c@convention\c@theorem
\let\c@notation\c@theorem
\let\c@openproblem\c@theorem
\makeatother
\makeautorefname{lemmadef}{Lemma-Definition}

\newcommand{\eqdef}{\;{:=}\;}
\newcommand{\fedqe}{\;{=:}\;}

\newcommand{\op}{\operatorname}

\newcommand{\Hom}{\op{Hom}}
\newcommand{\Ker}{\op{Ker}}
\newcommand{\Coker}{\op{Coker}}

\newcommand{\tensor}{\otimes}

\newcommand{\feed}{\,{-}\hspace{-1.6pt}\mbox{\raisebox{2.65pt}
{\mathsurround 0pt$\shortmid$}}\,}

\newcommand{\union}{\bigcup}

\newcommand{\Zmodule}{{\mskip0mu\underline{\mskip-0mu\mathbb{Z}\mskip-1mu}\mskip1mu}}

\let\xysavmatrix\xymatrix
\def\xymatrix{\disablesubscriptcorrection\xysavmatrix}
\AtBeginDocument{\def\notin{\not\in}}

\begin{document}

\begin{asciiabstract}
The embedded contact homology (ECH) of a 3-manifold with a contact
form is a variant of Eliashberg-Givental-Hofer's symplectic field
theory, which counts certain embedded J-holomorphic curves in the
symplectization.  We show that the ECH of T^3 is computed by a
combinatorial chain complex which is generated by labeled convex
polygons in the plane with vertices at lattice points, and whose
differential involves `rounding corners'.  We compute the homology
of this combinatorial chain complex.  The answer agrees with the
Ozsvath--Szabo Floer homology HF^+(T^3).
\end{asciiabstract}

\begin{htmlabstract}
The embedded contact homology (ECH) of a 3&ndash;manifold with a contact
form is a variant of Eliashberg&ndash;Givental&ndash;Hofer's symplectic field
theory, which counts certain embedded J&ndash;holomorphic curves in the
symplectization.  We show that the ECH of T<sup>3</sup> is computed by a
combinatorial chain complex which is generated by labeled convex
polygons in the plane with vertices at lattice points, and whose
differential involves "rounding corners".  We compute the homology
of this combinatorial chain complex.  The answer agrees with the
Ozsv&aacute;th&ndash;Szab&oacute; Floer homology HF<sup>+</sup>(T<sup>3</sup>).
\end{htmlabstract}

\begin{abstract}
The embedded contact homology (ECH) of a 3--manifold with a contact
form is a variant of Eliashberg--Givental--Hofer's symplectic field
theory, which counts certain embedded $J$--holomorphic curves in the
symplectization.  We show that the ECH of $T^3$ is computed by a
combinatorial chain complex which is generated by labeled convex
polygons in the plane with vertices at lattice points, and whose
differential involves ``rounding corners''.  We compute the homology
of this combinatorial chain complex.  The answer agrees with the
Ozsv\'{a}th--Szab\'{o} Floer homology $HF^+(T^3)$.
\end{abstract}

\maketitle

\section{Introduction}
\label{sec:introduction}

\subsection{Motivation}
\label{sec:motivation}

Let $Y$ be a closed oriented 3--manifold with a contact form, ie a
$1$--form $\lambda$ such that $\lambda\wedge d\lambda > 0$.  The
corresponding contact structure is the 2--plane field
$\xi\eqdef\Ker(\lambda)$; this is oriented by $d\lambda$.  Also
associated to $\lambda$ is the Reeb vector field $R$ characterized by
$\lambda(R)=1$ and $R\feed d\lambda=0$.  A periodic orbit of the Reeb
flow $R$ is called a Reeb orbit.  We assume that the Reeb orbits are
nondegenerate or Morse--Bott.  For $\Gamma\in H_1(Y)$, one can then
define the ``embedded contact homology'' $ECH_*(Y,\lambda;\Gamma)$, as
we explain in \fullref{sec:defineECH}.  This is the homology of a chain
complex which is generated by certain unions of Reeb orbits with total
homology class $\Gamma$.  The differential counts ``maximal index''
$J$--holomorphic curves in $\R\times Y$, for a suitable almost complex
structure $J$.  These maximal index curves turn out to be embedded,
except that they may contain multiple covers of ``trivial cylinders''
$\R\times\gamma$ where $\gamma$ is a Reeb orbit.  The embedded contact
homology is a relatively $\Z/N$--graded $\Z$--module, where $N$ is the
divisibility of the image of $c_1(\xi)+2\op{PD}(\Gamma)$ in
$\Hom(H_2(Y),\Z)$.  When $\Gamma=0$, there is a canonical absolute
$\Z/N$--grading.

Embedded contact homology (ECH) is analogous to the periodic Floer
theory (PFH) for mapping tori considered by Hutchings, Sullivan and
Thaddeus \cite{pfh1,pfh2,pfh3}.  ECH is similar to the symplectic
field theory (SFT) of Eliashberg--Givental--Hofer
\cite{eliashberg-givental-hofer00}, but has different generators and
grading, and counts more restricted $J$--holomorphic curves.  Unlike
SFT, which is highly sensitive to the contact structure, ECH is
conjectured to be a topological invariant, except that it detects the
Euler class of the contact structure via the identification
\eqref{eqn:turaev} below.  More precisely, we conjecture that
$ECH_*(Y,\lambda;\Gamma)$ agrees with the Seiberg--Witten Floer
homology $\widecheck{HM}_*(-Y)$ defined by Kronheimer and Mrowka
\cite{kronheimer-mrowka06} and summarized with
Ozsv\'ath and Szab\'o in \cite{kmos}, or the conjecturally isomorphic
Ozsv\'{a}th--Szab\'{o} Floer homology $HF^+_*(-Y)$ defined in
\cite{ozsvath-szabo0101}, as follows.  The latter two Floer homologies
depend on the choice of a spin-c structure $\frak{s}$ on $-Y$, which
is equivalent to a spin-c structure on $Y$.  Let $\op{Spin}^c(Y)$
denote the set of spin-c structures on $Y$.  Recall that an oriented
2--plane field on $Y$ determines a spin-c structure, cf Turaev
\cite{turaev97}.  Hence the contact structure $\xi$ gives rise to an
$H^2(Y;\Z)$--equivariant bijection

\vspace{-25pt}
\begin{equation}
\label{eqn:turaev}
\frak{s}_\xi\co H_1(Y) \stackrel{\simeq}{\longrightarrow} \op{Spin}^c(Y)
\end{equation}

\vspace{-10pt}
sending $0$ to the spin-c structure determined by the oriented
$2$--plane field $\xi$.

\begin{conjecture}
\label{conj:big}
Let $Y$ be a closed oriented $3$--manifold with a contact form
$\lambda$.  Then for $\Gamma\in H_1(Y)$, the embedded contact homology
is related to the Seiberg--Witten and Ozsv\'{a}th--Szab\'{o} Floer
homologies by
\begin{equation}
\label{eqn:BCI}
\widecheck{HM}_*(-Y,\frak{s}_\xi(\Gamma)) \simeq ECH_*(Y,\lambda;\Gamma) 
\simeq HF^+_*(-Y,\frak{s}_\xi(\Gamma)),
\end{equation}
up to a grading shift.
\end{conjecture}

Recall that Taubes's ``SW=Gr'' theorem \cite{taubes00} states that the
Seiberg--Witten invariant of a closed symplectic 4--manifold $X$ is
equivalent to the ``Gromov invariant'', which is a certain count of
embedded (except for multiply covered tori) $J$--holomorphic curves in
$X$.  The conjectural relation between ECH and Seiberg--Witten Floer
homology can be regarded as an analogue of ``SW=Gr'' for the
noncompact symplectic manifold $\R\times Y$.  This was the original
motivation for the definition of PFH and ECH.  Note also that the
Ozsv\'{a}th--Szab\'{o} Floer homology has been given a four-dimensional
reformulation by Lipshitz \cite{lipshitz}. It is possible that the
latter could be directly related to ECH by defining a more general
theory including both as special cases.

A proof of \fullref{conj:big}, while perhaps a long way off,
would have implications for contact dynamics.  For example, one
version of the {\em Weinstein conjecture\/} asserts that any contact
1--form $\lambda$ on a closed oriented 3--manifold $Y$ has a Reeb orbit,
see eg
Abbas--Cieliebak--Hofer \cite{ach}. If $\lambda$ has no Reeb orbit, then by definition, the
ECH chain complex has only one generator given by the empty set of
Reeb orbits, so
\[
ECH_*(Y,\lambda;\Gamma) = \left\{ \begin{array}{cl} \Z, &
\mbox{$\Gamma=0$ and $*=0$,}\\ 0, & \mbox{otherwise}.
\end{array}\right.
\]
However, Tom Mrowka has pointed out to us that by results
in Kronheimer--Mrowka \cite{kronheimer-mrowka06}, for any closed
oriented 3--manifold $Y$, if $\frak{s}$ is a spin-c structure with
$c_1(\frak{s})$ torsion then $\widecheck{HM}_*(Y,\frak{s})$ is infinitely
generated.  Since $TY$ is a trivial bundle, one can always find a spin-c
structure $\frak{s}$ with $c_1(\frak{s})=0$.  Therefore, the first part
of \fullref{conj:big} implies the Weinstein conjecture for every closed
oriented 3--manifold.

ECH has some additional structure, analogous to structures in the
Seiberg--Witten and Ozsv\'{a}th--Szab\'{o} Floer homologies.  For
example, there is a canonical element
\begin{equation}
\label{eqn:contactInvariant}
c(\lambda) \in ECH_0(Y,\lambda;0).
\end{equation}
In the ECH chain complex, the homology class $c(\lambda)$ is
represented by the empty set of Reeb orbits.  Under the conjectured
isomorphisms \eqref{eqn:BCI}, $c(\lambda)$ may agree with the
Ozsv\'{a}th--Szab\'{o} contact invariant \cite{ozsvath-szabo0210}, and
the Seiberg--Witten analogue implicit in
the paper by Kronheimer and Mrowka \cite{kronheimer-mrowka97}.

Further motivation for studying ECH (and PFH) is that it is expected to
be the recipient of (yet to be defined) relative Gromov invariants of
symplectic 4--manifolds with boundary.  For example, Taubes has proposed
\cite{taubes98,taubes01} that the Gromov invariant may be extended to
near-symplectic 4--manifolds by counting $J$--holomorphic curves in
the complement of the circles where the near-symplectic form vanishes.
We expect such a counting invariant to take values in the embedded
contact homology of a disjoint union of
$S^1\times S^2$'s, one for each vanishing circle, with the contact
form studied by Taubes in
\cite{taubes02}.  Also, the relative Gromov invariants should enter
into gluing formulas for Gromov invariants of closed symplectic
4--manifolds cut along 3--manifolds.

Much of the embedded contact homology story is still conjectural.  In
particular, a proof that $ECH_*(Y,\lambda;\Gamma)$ is well-defined is
currently in preparation; the precise statement is given here as
\fullref{conj:ECHDefined}.  In any case, the results in this
paper from \fullref{sec:preliminaries} to \fullref{sec:interlude}, while
motivated by this conjecture, are logically independent of it.

\subsection{The embedded contact homology of $T^3$}

This paper is concerned with computations of embedded contact
homology.  We will restrict attention to the example of $Y=T^3$,
although the methods developed here are applicable to some other
simple contact manifolds such as $S^1\times S^2$, or $T^2$--bundles
over $S^1$, see \fullref{sec:otherManifolds}.  For each positive integer
$n$, there is a standard contact form $\lambda_n$ on $T^3$ defined as
follows.  We choose the following coordinates on $T^3$ that depend on $n$:
\begin{equation}
\label{eqn:T^3}
T^3 = S^1\times T^2 = (\R/2\pi n\Z)_\theta \times (\R^2/\Z^2)_{x,y}
\end{equation}
Then
\begin{equation}
\label{eqn:lambda_n}
\lambda_n\eqdef \cos \theta\, dx + \sin\theta\, dy.
\end{equation}
The associated Reeb vector field is given by
\[
R = \cos\theta\,\partial_x + \sin\theta\,\partial_y.
\]
In particular, $\lambda_n$ is a Morse--Bott contact form; for each
$\theta\in\R/2\pi n\Z$ with $\tan\theta\in\Q\cup\{\infty\}$, there is
an $S^1$--family of Reeb orbits in $\{\theta\}\times T^2$.

It turns out that to compute ECH in this example, for suitable almost
complex structures $J$, the relevant $J$--holomorphic curves can be
counted quite explicitly.  For this purpose we modify some arguments
from our previous paper \cite{pfh3} on the PFH of a Dehn twist, and use
some results of Taubes \cite{taubes02}.  Consequently, for
\[
\Gamma \in \Z^2 = H_1(T^2) \subset H_1(S^1\times T^2),
\]
we can define
a combinatorial chain complex $\wwbar{C}_*(2\pi n;\Gamma)$, see
\fullref{sec:combinatoricsIntro} and \fullref{sec:polygonComplexes}, whose
homology $\wwbar{H}_*(2\pi n;\Gamma)$ agrees with the embedded
contact homology of $T^3$.  (Throughout this paper we adopt
the convention that changing the letter `$C$' to `$H$' indicates
passing from a chain complex to its homology.)  Namely:

\begin{theorem}
\label{thm:echt3}
If \fullref{conj:ECHDefined} holds (so that ECH is
well-defined), then
\[
ECH_*(T^3,\lambda_n;\Gamma) \simeq \wwbar{H}_*(2\pi n;\Gamma).
\]
\end{theorem}
Note that since all Reeb orbits have homology classes in the subgroup
$H_1(T^2)$, the ECH automatically vanishes for $\Gamma\in H_1(S^1\times
T^2)\setminus H_1(T^2)$, because the chain complex has no generators.

As will be explained in \fullref{sec:additionalStructure}, ECH has some
variants and additional structure.  In particular, there is a degree
$-2$ operation
\[
U\co  ECH_*(Y,\lambda;\Gamma) \longrightarrow ECH_{*-2}(Y,\lambda;\Gamma),
\]
which counts $J$--holomorphic curves with a marked point mapping to a
chosen point in $\R\times Y$.  In the case of $T^3$, the operation $U$
corresponds to a combinatorial chain map
\[
U\co \wwbar{C}_*(2\pi n;\Gamma) \longrightarrow
\wwbar{C}_{*-2}(2\pi n;\Gamma)
\]
defined in \fullref{sec:U}.  We can now state our main computational
result:

\begin{theorem}
\label{thm:main}
For every positive integer $n$:
\begin{enumerate}
\item[(a)] If $\Gamma\neq 0$, then
$
\wwbar{H}_*(2\pi n;\Gamma) = 0.
$
\item[(b)]
For $\Gamma=0$,
\[
\wwbar{H}_i(2\pi n;0) \simeq \left\{\begin{array}{cl} \Z^3,& i\ge
0,\\
0, & i<0.
\end{array}
\right.
\]
\item[(c)] For all $i\ge 2$, the map $U$ induces an isomorphism
\[
U\co  \wwbar{H}_i(2\pi n;0) \stackrel{\simeq}{\longrightarrow}
\wwbar{H}_{i-2}(2\pi n;0).
\]
\end{enumerate}
\end{theorem}

In particular, the ECH of $(T^3,\lambda_n)$ does not depend on $n$.
By contrast, the simplest version of SFT, namely cylindrical contact
homology, distinguishes the contact structures
$\xi_n=\Ker(\lambda_n)$; see Eliashberg--Givental--Hofer
\cite[Theorem~1.9.9]{eliashberg-givental-hofer00}, and for
generalizations see
Bourgeois--Colin \cite{bourgeois-colin}.  On the other hand, the contact invariant
\eqref{eqn:contactInvariant} for $\lambda_n$ does depend on $n$, see
\fullref{sec:contactInvariant}.

The above computation of the ECH of $T^3$, together with some
additional structure on it described in \fullref{sec:h1action}, agree
perfectly with $HF^+(T^3)$ as computed in
Ozsv\'ath--Szab\'o \cite{ozsvath-szabo0110}, and also $\widecheck{HM}_*(T^3)$
as computed in Kronheimer--Mrowka \cite{kronheimer-mrowka06}.
This provides a nontrivial check of \fullref{conj:big}.

\subsection{Rounding corners of polygons}
\label{sec:combinatoricsIntro}

We now introduce the combinatorial chain complex $\wwbar{C}_*(2\pi
n;\Gamma)$, along with two variants $\wwtilde{C}_*(2\pi n;\Gamma)$
and $C_*(2\pi n;\Gamma)$, in the simplest case where $n=1$ and
$\Gamma=0$.

\subsubsection{The generators}
The complex $C_*=C_*(2\pi;0)$ is a free $\Z$--module.  A generator of
$C_*$ is a convex polygon in $\R^2$, possibly a $2$--gon or a point,
such that the corners are lattice points, and every edge is labeled
either `$e$' or `$h$'.  To fix the signs in the differential, we
choose an ordering of the `$h$' edges, and we declare that a
reordering of the `$h$' edges multiplies the generator by the sign of
the reordering permutation.

\subsubsection{The grading}
The grading, or index, of a generator $\alpha$ is defined by
\begin{equation}
\label{eqn:index1}
I(\alpha)\eqdef 2(\#L(\alpha)-1)-\#h(\alpha).
\end{equation}
Here $\#L(\alpha)$ denotes the cardinality of the set $L(\alpha)$ of
lattice points on the polygon or enclosed by it, and $\#h(\alpha)$
denotes the number of `$h$' edges.  By Pick's formula for the area of
a lattice polygon, equation \eqref{eqn:index1} is equivalent to
\begin{equation}
\label{eqn:pick}
I(\alpha) = 2\op{Area}(\alpha)+\#\ell(\alpha)-\#h(\alpha).
\end{equation}
Here $\op{Area}(\alpha)$ denotes the area enclosed by $\alpha$, and
$\ell(\alpha)$ denotes the sum of the divisibilities of the edges of
the polygon.

For example, equation \eqref{eqn:pick} implies that $I(\alpha)\ge 0$,
and the only index zero generators are the following:
\begin{itemize}
\item
Points, which we denote by $p(u)$ where
$u\in\Z^2$.
\item $2$--gons with vertices $u,v\in\Z^2$ with $u-v$ indivisible and
  with both edges labeled `$h$'.  If the edge from $u$ to $v$ is first
  in the ordering, then we denote this generator by $h(u,v)$.  Thus
  $h(v,u)=-h(u,v)$.
\end{itemize}

\subsubsection{The differential}
\label{sec:introduceDifferential}

The differential $\delta\co C_*\to C_{*-1}$ is defined as follows.
Roughly, if $\alpha$ is a generator, then $\delta\alpha$ is the signed
sum of all generators $\beta$ obtained by ``rounding a corner'' and
``locally losing one `$h${'}''.  More precisely:

\begin{itemize}
\item ``Rounding a corner'' means that $L(\beta)=L(\alpha)\setminus c$
  where $c$ is a corner of $\alpha$.
  
\item ``Locally losing one `$h${'}'' means the following.  First, at
least one of the two edges in $\alpha$ adjacent to $c$ must be labeled
  `$h$'.  Second, of the edges in $\beta$ that are created or
  shortened by rounding the corner $c$, all are labeled `$e$', except
  for one when both edges adjacent to $c$ are labeled `$h$'.  Finally,
  all other edges in $\beta$ have the same labels as the corresponding
  edges in $\alpha$.
  
\item To determine the sign, let $\theta$ denote (one of) the `$h$'
  edge(s) of $\alpha$ adjacent to $c$.  Without loss of generality,
  $\theta$ is last in the ordering of the `$h$' edges of $\alpha$,
  while the remaining `$h$' edges of $\alpha$ are ordered the same way
  as the `$h$' edges of $\beta$ under the obvious bijection between
  them.  Then the differential coefficient
  $\langle\delta\alpha,\beta\rangle = +1$ if $\theta$ comes
  immediately after $c$ as we traverse $\alpha$ counterclockwise.  If
  $\theta$ comes immediately before $c$, then
  $\langle\delta\alpha,\beta\rangle = -1$.
\end{itemize}

Here is a random example of $\delta$.  In the pictures below, the
unmarked edges are labeled `$e$', and on the left side the bottom
`$h$' edge is first in the ordering.
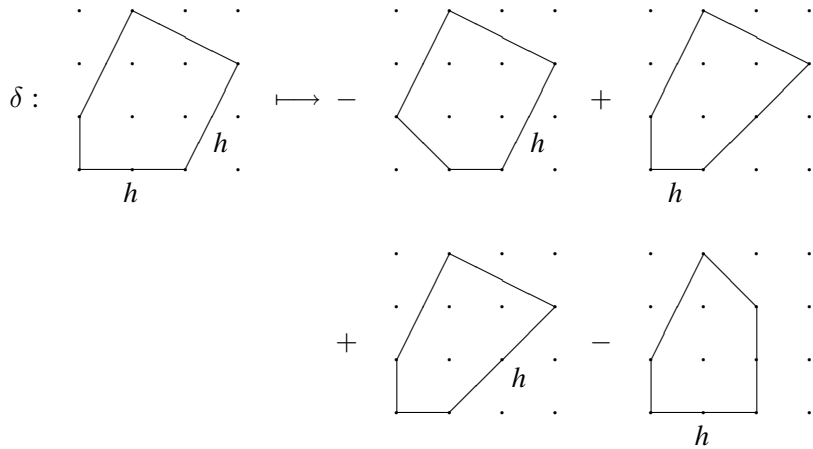
\[
\begin{split}
\delta:
\raisebox{-44pt}{
\begin{picture}(77,88)(-8,-20)
\put(0,0){.}
\put(20,0){.}
\put(40,0){.}
\put(60,0){.}
\put(0,20){.}
\put(20,20){.}
\put(40,20){.}
\put(60,20){.}
\put(0,40){.}
\put(20,40){.}
\put(40,40){.}
\put(60,40){.}
\put(0,60){.}
\put(20,60){.}
\put(40,60){.}
\put(60,60){.}
\put(1.5,0.5){\line(1,0){40}}
\put(41.5,0.5){\line(1,2){20}}
\put(1.5,0.5){\line(0,1){20}}
\put(1.5,20.5){\line(1,2){20}}
\put(21.5,60.5){\line(2,-1){40}}
\put(18,-12){$h$}
\put(52,7){$h$}
\end{picture}}
\;
\longmapsto
\;
&-
\raisebox{-44pt}{
\begin{picture}(77,88)(-8,-20)
\put(0,0){.}
\put(20,0){.}
\put(40,0){.}
\put(60,0){.}
\put(0,20){.}
\put(20,20){.}
\put(40,20){.}
\put(60,20){.}
\put(0,40){.}
\put(20,40){.}
\put(40,40){.}
\put(60,40){.}
\put(0,60){.}
\put(20,60){.}
\put(40,60){.}
\put(60,60){.}
\put(21.5,0.5){\line(1,0){20}}
\put(41.5,0.5){\line(1,2){20}}
\put(1.5,20.5){\line(1,-1){20}}
\put(1.5,20.5){\line(1,2){20}}
\put(21.5,60.5){\line(2,-1){40}}
\put(52,7){$h$}
\end{picture}}
\;
+
\raisebox{-44pt}{
\begin{picture}(77,88)(-8,-20)
\put(0,0){.}
\put(20,0){.}
\put(40,0){.}
\put(60,0){.}
\put(0,20){.}
\put(20,20){.}
\put(40,20){.}
\put(60,20){.}
\put(0,40){.}
\put(20,40){.}
\put(40,40){.}
\put(60,40){.}
\put(0,60){.}
\put(20,60){.}
\put(40,60){.}
\put(60,60){.}
\put(1.5,0.5){\line(1,0){20}}
\put(21.5,0.5){\line(1,1){40}}
\put(1.5,0.5){\line(0,1){20}}
\put(1.5,20.5){\line(1,2){20}}
\put(21.5,60.5){\line(2,-1){40}}
\put(8,-12){$h$}
\end{picture}}
\\
&+
\raisebox{-44pt}{
\begin{picture}(77,88)(-8,-20)
\put(0,0){.}
\put(20,0){.}
\put(40,0){.}
\put(60,0){.}
\put(0,20){.}
\put(20,20){.}
\put(40,20){.}
\put(60,20){.}
\put(0,40){.}
\put(20,40){.}
\put(40,40){.}
\put(60,40){.}
\put(0,60){.}
\put(20,60){.}
\put(40,60){.}
\put(60,60){.}
\put(1.5,0.5){\line(1,0){20}}
\put(21.5,0.5){\line(1,1){40}}
\put(1.5,0.5){\line(0,1){20}}
\put(1.5,20.5){\line(1,2){20}}
\put(21.5,60.5){\line(2,-1){40}}
\put(45,10){$h$}
\end{picture}}
\;
-
\raisebox{-44pt}{
\begin{picture}(77,88)(-8,-20)
\put(0,0){.}
\put(20,0){.}
\put(40,0){.}
\put(60,0){.}
\put(0,20){.}
\put(20,20){.}
\put(40,20){.}
\put(60,20){.}
\put(0,40){.}
\put(20,40){.}
\put(40,40){.}
\put(60,40){.}
\put(0,60){.}
\put(20,60){.}
\put(40,60){.}
\put(60,60){.}
\put(1.5,0.5){\line(1,0){40}}
\put(41.5,0.5){\line(0,1){40}}
\put(1.5,0.5){\line(0,1){20}}
\put(1.5,20.5){\line(1,2){20}}
\put(21.5,60.5){\line(1,-1){20}}
\put(18,-12){$h$}
\end{picture}}
\end{split}
\]

It is shown in \fullref{sec:combinatorialDifferential} that $\delta$ has
degree $-1$ and $\delta^2=0$.

\subsubsection{The homology}

To further illustrate the definitions, let us compute the degree--0
homology $H_0$.  The index one generators are the following:
\begin{itemize}
\item
$2$--gons with no lattice points in the interiors of the edges and one
edge labeled `$e$' and one edge labeled `$h$'.  We denote such a
generator by $e(u,v)$, where $u$ and $v$ are the two corners and the
edge from $u$ to $v$ is labeled `$e$'.
\item
Triangles enclosing no lattice points except the corners $u,v,w$, with
all three edges labeled `$h$'.  We denote such a generator by
$h(u,v,w)$, where $u,v,w$ are listed in counterclockwise order and the
edges are ordered counterclockwise.
\end{itemize}
It follows from the definitions that
\begin{gather*}
\delta p(u) = \delta h(u,v) = 0,\\
\delta e(u,v) = p(u)-p(v),\\
\delta h(u,v,w)=h(u,v)+h(v,w)+h(w,u).
\end{gather*}
Therefore $H_0$ is generated by the
homology classes $[p(u)]$ and $[h(u,v)]$
modulo the relations
\[
\begin{split}
[p(u)] &= [p(v)],\\
[h(u,v)]+[h(v,w)] &= [h(u,w)].
\end{split}
\]
The computation of the higher homology is more complicated, but the
result is simpler.
Let $\Lambda$ be a convex polygon that encloses and contains $k$
lattice points with $k\ge 2$.  Let $E_\Lambda$ denote the generator
consisting of the polygon $\Lambda$ with all edges labeled `$e$'.  Let
$H_\Lambda$ denote the sum of all generators consisting of the polygon
$\Lambda$ with one edge labeled `$h$' and all other edges labeled
`$e$'.  It follows from the definitions that $E_\Lambda$ and
$H_\Lambda$ are cycles.  We will see in \fullref{sec:beginComputation}
that the homology classes of these cycles depend only on $k$.
Moreover,
\[
H_{2k-2}=\Z\left\{E_\Lambda\right\}, \quad \quad H_{2k-3} =
\Z\left\{H_\Lambda\right\}.
\]

\subsubsection{Variants}
\label{sec:variantsIntro}
Let $\wwtilde{C}_*$ denote the chain complex $C_*$ regarded as a
module over the group ring $\Z[\Z^2]$, where $\Z^2$ acts on the
generators by translation in the plane.  The homology
$\wwtilde{H}_*$ can be read off from the preceding calculations.
Let $\mc{I}(\Z^2)$ denote the augmentation ideal in $\Z[\Z^2]$, ie
the kernel of the augmentation map $\Z[\Z^2]\to\Z$ sending a group
ring element to the sum of its coefficents.  Also, let $\Zmodule$
denote the $\Z[\Z^2]$--module with one generator on which $\Z^2$ acts
by the identity, ie
\[
\Zmodule \eqdef \{a\mid (x-1)a=(y-1)a=0\}
\]
where $x$ and $y$ denote generators of $\Z^2$.  Then we have
\[
\wwtilde{H}_* \simeq \left\{\begin{array}{cl}
0, & *<0,\\
\mc{I}(\Z^2)\oplus \Zmodule, & *=0,\\
\Zmodule, & *>0.
\end{array}\right.
\]
Here the isomorphism $\wwtilde{H}_0\simeq
\mc{I}(\Z^2)\oplus\Zmodule$ sends $[p(u)]$ to the generator of
$\Zmodule$ and $[h(u,v)]$ to $u-v\in \mc{I}(\Z^2)$.

We can obtain another complex $\wwbar{C}_*$ over $\Z$ by declaring
generators to be equivalent when they differ by translation in the
plane.  The homology $\wwbar{H}_*$ of $\wwbar{C}_*$ is partially
but not entirely determined by $\wwtilde{H}_*$, via a ``universal
coefficient spectral sequence''
\begin{equation}
\label{eqn:TorSS}
E^2_{p,q}=\op{Tor}_p\left(\wwtilde{H}_q,\Zmodule\right) \Longrightarrow
\wwbar{H}_{p+q}.
\end{equation}
It turns out that some of the differentials in the spectral sequence
\eqref{eqn:TorSS} are nonzero, and we will find in \fullref{sec:HBar}
that
\[
\wwbar{H}_*\simeq\left\{\begin{array}{cl}
0, & *<0,\\
\Z^3, & *\ge 0.
\end{array}\right.
\]

\subsubsection{Geometric interpretation}

The rough idea of the relation between $\wwbar{H}_*$ and
$ECH_*(T^3,\lambda_1;0)$ is as follows.  For a generator of
$\wwbar{C}_*$, each edge of the polygon corresponds to a circle of
Reeb orbits in $T^3$.  The labels `$e$' and `$h$' reflect the fact
that each circle of Reeb orbits for the Morse--Bott contact form
$\lambda_1$ can be perturbed into two nondegenerate orbits, one
elliptic and one hyperbolic.  It turns out that in $\R\times T^3$,
every $J$--holomorphic curve counted by the ECH differential consists
of one embedded genus zero component with two positive ends
(corresponding to the edges adjacent to the corner being rounded) and
an arbitrary number of negative ends (corresponding to the edges
created by rounding), together with some $\R$--invariant cylinder
components (corresponding to the edges or parts thereof not involved
in the corner rounding).  The chain complex $\wwtilde{C}_*$
corresponds to a ``partially twisted'' version of ECH in which one
keeps track of some information about the relative homology classes of
the $J$--holomorphic curves.  

\subsubsection{The rest of the paper}
The combinatorial chain complexes for general $n$ and $\Gamma$,
roughly speaking, involve left-turning polygonal paths of rotation
number $n$ and period $\Gamma$.  The precise definitions require some
care and are given in \fullref{sec:polygonComplexes}, after some
combinatorial preliminaries in \fullref{sec:preliminaries}.  The
combinatorial chain map $U$ is defined in \fullref{sec:U}; aside from
its significance for ECH, it will help compute the homology of the
combinatorial chain complexes.  In
\fullref{sec:beginComputation}--\fullref{sec:HBar} we compute
the homology of all of the above combinatorial chain complexes and in
particular prove \fullref{thm:main}.

In \fullref{sec:uniqueness} we establish an axiomatic characterization of
the combinatorial chain complexes.  In \fullref{sec:interlude} we recall
and prove some relevant facts about $J$--holomorphic curves in
$\R\times T^3$.  In \fullref{sec:ECH} we outline the definition of
ECH and prove \fullref{thm:echt3}.  The proof uses the results in
\fullref{sec:interlude} to show that the chain complex computing the ECH
of $T^3$ satisfies the axioms in \fullref{sec:uniqueness}.  In
\fullref{sec:concludingRemarks} we make some concluding remarks.

\subsubsection{Acknowledgments}

We thank Y Eliashberg, P Kronheimer, T Mrowka, P Ozsv\'{a}th, Z Szab\'{o},
and C Taubes for enlightening discussions.  We thank the anonymous
referee for many helpful comments.  The first author was partially
supported by NSF grant DMS-0204681 and the Alfred P Sloan Foundation.
The second author was partially supported by NSF grant DMS-0305825.

\subsubsection{Index of frequently used notation}\qua

\begin{tabular}{ll}
$\Zmodule$ & the $\Z[\Z^2]$--module with one generator on which \\

& $\Z^2$ acts trivially, see \fullref{sec:variantsIntro}\\

$\mc{I}(\Z^2)$ & augmentation ideal, see \fullref{sec:variantsIntro} and
also \fullref{sec:zero} \\

$\Lambda\setminus c$ & rounding of $\Lambda$ at $c$, see
\fullref{sec:roundingCorners} \\

$\Lambda' \le \Lambda$ & $\Lambda'$ is to the left of $\Lambda$, see
\fullref{sec:partialOrder}\\

$\mc{A}$ & length of a polygonal path defined in equation
\eqref{eqn:defineLength}, and\\

& symplectic action of an orbit set defined in \fullref{def:action}
\\

 $C_*(2\pi n;\Gamma)$ & combinatorial chain complex generated by
polygons\\

& of rotation number $n$ and period $\Gamma$, see \fullref{sec:mgpc} \\

$\wwtilde{C}_*(2\pi n;\Gamma)$ & above complex regarded as a
$\Z[\Z^2]$--module\\

$\wwbar{C}_*(2\pi n;\Gamma)$ & above complex modulo translation of
polygons\\

$C_*(\Lambda)$ & complex of polygons to the left of $\Lambda$, see
\fullref{sec:auxiliary}\\

$I$ & combinatorial relative index defined in \fullref{sec:mgpc}, and\\

& analytical relative index defined in \fullref{sec:II} \\

$\delta$ & combinatorial differential defined in \fullref{sec:mgpc} \\

$E_\Lambda, H_\Lambda$ & two distinguished cycles defined in
\fullref{sec:distinguishedCycles}\\

$c_\theta$ & corner at $\theta$, see \fullref{sec:U}\\

$U$ & combinatorial degree $-2$ chain map defined in \fullref{sec:U},
and\\

& analytical degree $-2$ map defined in \fullref{sec:ECHU}\\

$C^{(j)}$ &
subcomplex where $I-\#h=j$, see \fullref{definition:subcomplex}\\

$\partial_c$ & connecting homomorphism in long exact sequence in
\fullref{sec:exactSequence}\\

$F_\theta$ & flattening chain map, see \fullref{def:flattening}
\end{tabular}

\begin{tabular}{ll}
$Z_n(a,b)$ & a degree zero cycle defined in \fullref{sec:zero} \\

$CX_*$ & subcomplex consisting of $x$--axis polygons, see
\fullref{sec:splicing}
\\

$S$ & splicing chain map, see \fullref{sec:splicing}\\

$H_2(Y,\alpha,\beta)$ & relative homology classes as in
\fullref{def:affine} \\

$\mc{M}^J$ & moduli space of $J$--holomorphic
curves, see
\fullref{sec:interlude} \\

$\op{ind}$ & SFT index, see \fullref{def:ind} \\

$ECH_*(Y,\lambda;\Gamma)$ & embedded contact homology, see
\fullref{sec:defineECH} \\

$\widetilde{ECH}_*(Y,\lambda;\Gamma,G)$ & twisted embedded contact
homology, see \fullref{sec:defineECH}\\

$\partial$ & embedded contact homology differential, see
\fullref{sec:defineECH}
\end{tabular}

\section{Rounding corners of polygonal paths}
\label{sec:preliminaries}

In this section we lay the foundations for our combinatorial
investigations.  We are preparing to define (in
\fullref{sec:polygonComplexes}) a general combinatorial chain complex
$C_*(2\pi n;\Gamma)$ in which the generating polygons have rotation
number $n\ge 1$, and when $\Gamma\neq 0$ are periodic with period
$\Gamma$ rather than closed.  In \fullref{sec:admissiblePaths} we define
the relevant classes of polygonal paths, which we call ``admissible
paths''.  In \fullref{sec:roundingCorners} we define the corner rounding
operation for these polygonal paths.  In \fullref{sec:partialOrder} we
introduce a closely related partial order on admissible paths, and we
prove various facts about corner rounding and the partial order which
will be needed later.  A dictionary between some of this combinatorics
and the geometry of pseudoholomorphic curves in $\R\times T^3$ will be
given later in \fullref{sec:interlude}.

\subsection{Admissible (left-turning polygonal) paths}
\label{sec:admissiblePaths}

We now define three types of ``admissible paths'': ``open'',
``closed'', and ``periodic''.  Closed and periodic admissible paths
will be used to define the chain complex
$$C_*(2\pi n;\Gamma)$$
for
$\Gamma=0$ and $\Gamma\neq0$ respectively, while we will use open
admissible paths to help compute its homology.

Let $\Theta$ denote the set of $\theta\in\R$ such that
$\tan\theta\in\Q\cup\{\infty\}$.  For each $\theta\in \Theta$, there is a
unique integer vector
$$\begin{pmatrix}x_\theta\\y_\theta\end{pmatrix}\in\Z^2$$
such that
$x_\theta,y_\theta$ are relatively prime and such that 
$(\cos\theta, \sin\theta)$
is a positive real multiple of
$(x_\theta,y_\theta)$.

\begin{definition}
  An {\em open admissible path\/} defined on an interval $I\subset\R$
  is a locally constant map $\Lambda\co I\setminus T\to\Z^2$, where
  $T\subset \Theta\cap \op{int}(I)$ is a finite set of {\em edges\/}, and
  there is a {\em multiplicity\/} function $m\co T\to\Z_{>0}$ such that
  $\Lambda$ satisfies the ``jumping condition''
\begin{equation}
\label{eqn:jump1}
\frac{d}{dt}\Lambda(t)=\sum_{\theta\in T}
m(\theta)\begin{pmatrix}x_\theta\\y_\theta\end{pmatrix}\delta_\theta(t).
\end{equation}
Here $\delta_\theta$ denotes the delta function supported at $\theta$.
A {\em corner\/} of $\Lambda$ is a component of $I\setminus T$, other
than the first and last components.  The {\em endpoints\/} of
$\Lambda$ are its values on the first and last components of
$I\setminus T$.  A {\em kink\/} of $\Lambda$ is a corner between
consecutive edges $\theta_1,\theta_2$ with $\theta_2-\theta_1>\pi$.
\end{definition}

That is, $\Lambda$ is a polygonal path in the plane with
corners in $\Z^2$, parametrized so that it is usually stopped at a
corner, and jumps discontinuously to the next corner at time $t$ when
the vector $(\cos t,\sin t)$ points in the direction of the
corresponding edge.  That is, a smooth locally convex curve is
naturally parametrized by its tangent direction, and we are extending
this notion to polygonal paths.  We will sometimes abuse notation and
pretend that $\Lambda$ is a continuous, piecewise linear path, moving
along straight line segments from one vertex to the next.  The path
$\Lambda$ turns to the left at its corners, except at kinks where it
rotates by more than $\pi$ and hence may turn in any direction. Note
that by equation \eqref{eqn:jump1}, $\Lambda$ determines $m$, while
$m$ determines $\Lambda$ up to a $\Z^2$ translation ambiguity.

Now fix a positive integer $n$ and an integer vector $\Gamma\in\Z^2$.
Let $p\co \R\to\R/2\pi n\Z$ denote the projection.

\begin{definition}
\label{def:PAP}
A {\em periodic admissible path\/} of rotation number $n$ and period
$\Gamma$ is a locally constant map $\Lambda\co \R\setminus
p^{-1}(T)\to\Z^2$, where $T\subset p(\Theta)$ is a finite set of {\em
  edges\/}, and there is a {\em multiplicity\/} function
$m\co T\to\Z_{>0}$ such that $\Lambda$ satisfies the jumping condition
\begin{equation}
\label{eqn:jump2}
\frac{d}{dt}\Lambda(t)=\sum_{\theta\in p^{-1}(T)}
m(p(\theta))
\begin{pmatrix}x_\theta\\y_\theta\end{pmatrix}\delta_\theta(t).
\end{equation}
We also assume that $\Lambda$ satisfies 
the periodicity condition
\[
\Lambda(t+2\pi n) = \Lambda(t) + \Gamma,
\]
which by \eqref{eqn:jump2} is equivalent to
\[
\sum_{\theta\in T}m(\theta)\begin{pmatrix}x_\theta \\
y_\theta\end{pmatrix} = \Gamma.
\]

A {\em corner\/} of $\Lambda$ is a component of $(\R/2\pi n\Z)\setminus
T$.  A {\em kink\/} is a corner
of length greater than $\pi$.
If $\Gamma=0$, we say that $\Lambda$ is a {\em closed admissible
  path\/} of rotation number $n$.
\end{definition}

\begin{example}
\label{ex:simplestPath}
In the simplest case of this definition, where $n=1$ and $\Gamma=0$, a
closed admissible path $\Lambda$ of rotation number $1$ is equivalent
to a convex polygonal region $P$ in $\R^2$ (possibly a $2$--gon or a
$0$--gon) with corners in $\Z^2$.  The path $\Lambda$ traverses the
boundary of $P$ counterclockwise.
\end{example}

Here is an example of a closed admissible path with rotation number
$n=2$.  The bottom corner is a kink at which the path turns by angle
$5\pi/4$.
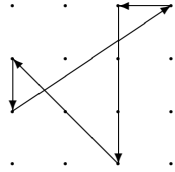
\begin{center}
\begin{picture}(77,70)(-8,-2)
\put(0,0){.}
\put(20,0){.}
\put(40,0){.}
\put(60,0){.}
\put(0,20){.}
\put(20,20){.}
\put(40,20){.}
\put(60,20){.}
\put(0,40){.}
\put(20,40){.}
\put(40,40){.}
\put(60,40){.}
\put(0,60){.}
\put(20,60){.}
\put(40,60){.}
\put(60,60){.}
\put(41.5,0.5){\vector(-1,1){40}}
\put(1.5,40.5){\vector(0,-1){20}}
\put(1.5,20.5){\vector(3,2){60}}
\put(61.5,60.5){\vector(-1,0){20}}
\put(41.5,60.5){\vector(0,-1){60}}
\end{picture}
\end{center}
Here is an example of a periodic admissible path with
$$n=1 \quad\text{and}\quad
\Gamma=\scriptsize\begin{pmatrix}2\\1\end{pmatrix}.$$
(The path continues infinitely in both directions.)
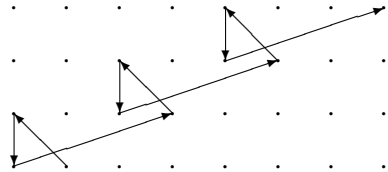
\begin{center}
\begin{picture}(157,70)(-8,-2)
\put(0,0){.}
\put(20,0){.}
\put(40,0){.}
\put(60,0){.}
\put(80,0){.}
\put(100,0){.}
\put(120,0){.}
\put(140,0){.}
\put(0,20){.}
\put(20,20){.}
\put(40,20){.}
\put(60,20){.}
\put(80,20){.}
\put(100,20){.}
\put(120,20){.}
\put(140,20){.}
\put(0,40){.}
\put(20,40){.}
\put(40,40){.}
\put(60,40){.}
\put(80,40){.}
\put(100,40){.}
\put(120,40){.}
\put(140,40){.}
\put(0,60){.}
\put(20,60){.}
\put(40,60){.}
\put(60,60){.}
\put(80,60){.}
\put(100,60){.}
\put(120,60){.}
\put(140,60){.}
\put(21.5,0.5){\vector(-1,1){20}}
\put(1.5,20.5){\vector(0,-1){20}}
\put(1.5,0.5){\vector(3,1){60}}
\put(61.5,20.5){\vector(-1,1){20}}
\put(41.5,40.5){\vector(0,-1){20}}
\put(41.5,20.5){\vector(3,1){60}}
\put(101.5,40.5){\vector(-1,1){20}}
\put(81.5,60.5){\vector(0,-1){20}}
\put(81.5,40.5){\vector(3,1){60}}
\end{picture}
\end{center}
For closed admissible paths with $n>1$, our definition keeps track of
the parametrization of the path; that is, $\Z/n$ acts nontrivially on
the set of closed admissible paths of rotation number $n$ by
precomposing $\Lambda$ with translations on $\R$ by multiples of
$2\pi$.  For periodic admissible paths with $\Gamma\neq 0$,
precomposing $\Lambda$ with translations by multiples of $2\pi$ gives
a $\Z$ action on the set of periodic admissible paths.  (The
parametrization is not indicated in the above two pictures.)

We say that two admissible paths $\Lambda,\Lambda'$ are
{\em of the same type\/} if they are both open and defined on the same
interval $I$ and have the same endpoints, or if they are both closed or
periodic for the same $n$ and $\Gamma$.

\subsection{Rounding corners of admissible paths}
\label{sec:roundingCorners}

We now define the corner rounding operation for admissible paths of
any type.

\begin{definition}
  Let $\Lambda$ be an admissible path, let $c$ be a corner of
  $\Lambda$ adjacent to consecutive edges $\theta_1,\theta_2$, and
  assume that $c$ is not a kink, ie $\theta_2-\theta_1\in(0,\pi]$.
  Define a new admissible path $\Lambda\setminus c$ of the
  same type as follows.
\begin{itemize}
\item If $\Lambda$ is an open admissible path defined on $I$, then
  $\Lambda\setminus c=\Lambda$ on $I\setminus [\theta_1,\theta_2]$.  If
  $\Lambda$ is periodic or closed, then $\Lambda\setminus c=\Lambda$ on
  $\R\setminus p^{-1}[\theta_1,\theta_2]$.
\item If $m$ and $m_c$ denote the multiplicity functions for $\Lambda$
  and $\Lambda\setminus c$ respectively, then
  $m_c(\theta_i)=m(\theta_i)-1$ for $i=1,2$.
\item If $\Lambda$ is an open admissible path, let $W\subset\Z^2$ be
  the set of lattice points enclosed by the triangle (or $2$--gon when
  $\theta_2-\theta_1=\pi$) in the plane whose corners are
  $$\Lambda(c)-\begin{pmatrix}x_{\theta_1}\\ y_{\theta_1}\end{pmatrix},\quad
    \Lambda(c),\quad
  \Lambda(c)+\begin{pmatrix}x_{\theta_2}\\ y_{\theta_2}\end{pmatrix}.$$
  Then $(\Lambda\setminus
  c)|_{(\theta_1,\theta_2)}$ traverses counterclockwise the boundary
  of the convex hull of $W\setminus\{\Lambda(c)\}$, except for
  the edge from
  $$\Lambda(c)+\begin{pmatrix}x_{\theta_2}\\
    y_{\theta_2}\end{pmatrix}\quad \text{to}\quad \Lambda(c)-\begin{pmatrix}
    x_{\theta_1}\\ y_{\theta_1}\end{pmatrix}.$$
  If $\Lambda$ is
  periodic or closed, then $\Lambda\setminus c$ is defined the same
  way on each component of $p^{-1}(\theta_1,\theta_2)$.
\end{itemize}
We say that $\Lambda\setminus c$ is obtained from $\Lambda$ by {\em
rounding the corner at $c$\/}.
\end{definition}


\begin{example}
\label{ex:simplestRounding}
Suppose $\Lambda$ is a closed admissible path of rotation number $1$,
corresponding to the boundary of a convex polygonal region $P$.  If
$c$ is a corner of $\Lambda$ mapping to a corner $\Lambda(c)$ of $P$,
then $\Lambda\setminus c$ corresponds to the boundary of the convex
hull of the set of lattice points in $P\setminus \{\Lambda(c)\}$.
\end{example}

\subsection{A partial order on admissible paths}
\label{sec:partialOrder}

We now introduce a partial order on the set of admissible paths and
collect some useful facts about it.  This partial order plays a
fundamental role both in our computation of the homology of the chain
complexes $C_*(2\pi n;\Gamma)$ and in the
connection with pseudoholomorphic curves in $\R\times T^3$.

Let $\Lambda$ and $\Lambda'$ be admissible paths of
the same type, as defined in \fullref{sec:admissiblePaths}.

\begin{definition}
\label{def:left}
We say that $\Lambda'$ is {\em to the left\/} of $\Lambda$, and we
write $\Lambda'\le\Lambda$, if
\begin{equation}
\label{eqn:left}
\det\left(\begin{matrix}\begin{matrix}
\cos t \\
\sin t
\end{matrix}
&
\Lambda'(t) - \Lambda(t)
\end{matrix}
\right) \ge 0
\end{equation}
for all $t$ in $I$ (if $\Lambda,\Lambda'$
are open) or $\R$ (if $\Lambda,\Lambda'$ are periodic or closed).
\end{definition}

Note that the left side of \eqref{eqn:left} is defined even if $t$ is
an edge of $\Lambda$ or $\Lambda'$, because by equation
\eqref{eqn:jump1} or \eqref{eqn:jump2}, the determinant in
\eqref{eqn:left} extends to a continuous function defined for all $t$
in $I$ or $\R$.

\begin{remark}
  Let $L_t$ and $L_t'$ denote the oriented lines through $\Lambda(t)$
  and $\Lambda'(t)$ respectively in the direction $(\cos t,\sin t)$.
  Then $\Lambda'\le\Lambda$ iff for all
  $t$, the oriented line $L_t'$ is (not necessarily strictly) to the
  left of $L_t$ in the usual sense.
\end{remark}

\begin{proposition}
\label{prop:partialOrder}
$\le$ is a partial order.
\end{proposition}

\begin{proof}
  By \eqref{eqn:left}, we have $\Lambda\le\Lambda$, and if
  $\Lambda_1\le\Lambda_2$ and $\Lambda_2\le\Lambda_3$, then
  $\Lambda_1\le \Lambda_3$.  

Now suppose that $\Lambda'\le \Lambda$ and $\Lambda\le\Lambda'$;
we must show that $\Lambda'=\Lambda$.  We have
\[
\det\left(\begin{matrix}\begin{matrix}
\cos t \\
\sin t
\end{matrix}
&
\Lambda'(t) - \Lambda(t)
\end{matrix}
\right) = 0
\]
for all $t$ for which $\Lambda,\Lambda'$ are defined.  Now
$\Lambda'(t)$ and $\Lambda(t)$ are locally constant and defined on the
complement of the edges of $\Lambda$ and $\Lambda'$.  Since the vector
$(\cos t,\sin t)$
rotates as $t$ varies, the vanishing of the determinant implies that
$\Lambda'(t)-\Lambda(t)=0$ on each interval between edges of $\Lambda$
and $\Lambda'$.
\end{proof}

\begin{example}
\label{example:leftInclusion}
For closed admissible paths $\Lambda$ and $\Lambda'$ of rotation
number $1$, corresponding to the boundaries of convex polygonal
regions $P$ and $P'$, we have
\[
\Lambda' \le \Lambda \Longleftrightarrow P'\subset P.
\]
\end{example}

\begin{proof}
$(\Leftarrow)$\qua Suppose $P'\subset P$.  We need to show that
$\Lambda' \le \Lambda$; equivalently, for each $t$ the point
$\Lambda'(t)$ is to the left of the line $L_t$.  But this
holds because convexity of $P$ implies that all of $P$ is to the
left of $L_t$, and $P'\subset P$.

$(\Rightarrow)$\qua Suppose $P'\not\subset P$.  Let $x\in P'$ maximize
distance to $P$.  Let $y\in P$ minimize distance to $x$.  If $x$ is
a corner of $P'$ and $y$ is a corner of $P$, then there exists
$t\in\R/2\pi\Z$ such that $x=\Lambda'(t)$ and $y=\Lambda(t)$, and
$x$ is strictly to the right of the line $L_t$, so
\begin{equation}
\label{eqn:right}
\det\left(\begin{matrix}\begin{matrix}
\cos t \\
\sin t
\end{matrix}
&
\Lambda'(t) - \Lambda(t)
\end{matrix}
\right)
<0.
\end{equation}
If $x$ is in the interior of an edge of $P'$ and $y$ is in the
interior of an edge of $P$, then these edges both point in the same
direction $(\cos t,\sin t)$
and the inequality \eqref{eqn:right} still holds.  If $x$ is a corner
of $P'$ and $y$ is in the interior of an edge of $P$ with tangent direction
$(\cos t,\sin t)$
then $x=\Lambda'(t)$ and \eqref{eqn:right} still holds; likewise if just
$y$ is a corner.
\end{proof}

There is a useful characterization of the corner rounding operation in
terms of the partial order $\le$.  Namely, $\Lambda\setminus c$ is the
maximal admissible path which is to the left of $\Lambda$ and which
does not ``go through'' the corner $c$ of $\Lambda$.  More precisely:

\begin{proposition}
\label{prop:roundingMaximal}
  Let $\Lambda$ be an admissible path and let $c$ be a corner
  of $\Lambda$ which is not a kink.  Then:
\begin{enumerate}
\item[(a)] $\Lambda\setminus c \le \Lambda$.
\item[(b)] If $\Lambda'\le \Lambda$, and if $\Lambda'$ disagrees
  with $\Lambda$ somewhere on the interval corresponding to $c$,
  then $\Lambda'\le\Lambda\setminus c$.
\end{enumerate}
\end{proposition}

\begin{lemma}
\label{lem:wedge}
Suppose $\Lambda'\le\Lambda$.  Let $t_1,t_2$ be real numbers in the
domain of $\Lambda,\Lambda'$ with $t_2-t_1\in(0,\pi)$.  Then the
restriction of $\Lambda'$ to $[t_1,t_2]$ maps to the wedge consisting
of those $x\in\R^2$ that are to the left of $L_{t_1}$ and $L_{t_2}$.
\end{lemma}

\begin{proof}
For $i=1,2$ and $t\in[t_1,t_2]$ define a piecewise constant function
\[
f_i(t) \eqdef \det\left(\begin{matrix}\begin{matrix}
\cos t_i \\
\sin t_i
\end{matrix}
&
\Lambda'(t) - \Lambda(t_i)
\end{matrix}
\right).
\]
Note that this is well defined even if $t_i$ is an edge of $\Lambda$,
as one can then take $\Lambda(t_i)$ to be any point on the
corresponding line in $\R^2$.  Now by equation \eqref{eqn:jump1} or
\eqref{eqn:jump2} and our assumption that $t_2-t_1\in(0,\pi)$, we have
\[
\frac{d}{dt}f_1(t)\ge 0,\quad\quad \frac{d}{dt}f_2(t)\le 0
\]
for $t\in[t_1,t_2]$.  On the other hand, since we assumed that
$\Lambda'\le\Lambda$, putting $t=t_i$ into the definition of $\le$ gives
\[
f_1(t_1)\ge 0,\quad\quad f_2(t_2)\ge 0.
\]
Therefore for all $t\in[t_1,t_2]$, we have $f_1(t),f_2(t)\ge 0$, which
means that $\Lambda'(t)$ is to the left of $L_{t_1}$ and $L_{t_2}$.
\end{proof}

\begin{proof}[Proof of \fullref{prop:roundingMaximal}]
To simplify notation we assume that $\Lambda$ is an open admissible
path; the proof when $\Lambda$ is a closed or periodic admissible path
works the same way.  Let $\theta_1,\theta_2$ be the edges of $\Lambda$
adjacent to $c$.  We assume that $\theta_2-\theta_1<\pi$ and leave the
easier, extreme case when $\theta_2-\theta_1=\pi$ as an exercise.
  
(a)\qua Equation~\eqref{eqn:left} holds for $t\notin[\theta_1,\theta_2]$,
as then $(\Lambda\setminus c)(t)=\Lambda(t)$ by definition of
$\Lambda\setminus c$.  Now suppose that $t\in[\theta_1,\theta_2]$.  We
need to show that $(\Lambda\setminus c)(t)$ is to the left of the line
$L_t$.  But this is clear since the entire triangle $W$ in the
definition of $\Lambda\setminus c$ is to the left of the line $L_t$.
  
(b)\qua Suppose $\Lambda'\le\Lambda$ and $\Lambda'\not\le\Lambda\setminus
c$; we will show that $\Lambda'$ agrees with $\Lambda$ on all of the
interval $(\theta_1,\theta_2)$. For $i=1,2$ let $V$ be the wedge as in
\fullref{lem:wedge} for $t_i=\theta_i$.  Let $P$ be the path in $V$
traced out by $(\Lambda\setminus c)|_{[\theta_1,\theta_2]}$.  The path
$P$ separates $V$ into two components, one bounded and one unbounded.
By \fullref{lem:wedge}, $\Lambda'(t)\in V$ for all
$t\in[\theta_1,\theta_2]$ for which $\Lambda'(t)$ is defined.  Since
$\Lambda'\not\le\Lambda\setminus c$, there exists
$t\in(\theta_1,\theta_2)$ such that if $L_t''$ denotes the oriented
line through $(\Lambda\setminus c)(t)$ in the direction $(\cos t,\sin t)$,
then $\Lambda'(t)$ is
strictly to the right of the line $L_t''$.  Since $P$ is convex and to
the left of $L_t''$, it follows that $\Lambda'(t)$ is in the bounded
component of $V\setminus P$.  By definition of $\Lambda\setminus c$,
there are no lattice points in the interior of the bounded component
of $V\setminus P$, so we must have $\Lambda'(t)=\Lambda(c)$.  Then
$\Lambda'$ must agree with $\Lambda$ on the entire interval
$(\theta_1,\theta_2)$, since if $\Lambda'$ had any jumps on this
interval, then by equation \eqref{eqn:jump1} for $\Lambda'$ it would
escape the wedge $V$.
\end{proof}

The following proposition shows that the relation $\le$ imposes strong
restrictions in the presence of kinks.

\begin{proposition}
\label{prop:stuckAtKink}
  Suppose $\Lambda$ has a kink at $c$ and
  $\Lambda'\le \Lambda$.  Then $\Lambda'$ agrees with $\Lambda$
  on the interval corresponding to $c$.
\end{proposition}

\begin{proof}
  Let $\theta_1,\theta_2$ be the edges of $\Lambda$ adjacent to $c$, and let
  $t\in(\theta_1,\theta_2)$; we must show that
  $\Lambda'(t)=\Lambda(c)$.  Since $\theta_2-\theta_1>\pi$, we can
  find $t_0$ with
\begin{equation}
\label{eqn:t0}
\theta_1 < t_0 < t < t_0+\pi < \theta_2.
\end{equation}
By \fullref{lem:wedge}, if
\[
t_0<t_1<t_2<t_0+\pi,
\]
then $\Lambda'[t_1,t_2]$ is contained in the wedge of $x\in\R^2$
to the left of $L_{t_1}$ and $L_{t_2}$.  
Taking the limit as $t_1\searrow t_0$ and $t_2\nearrow
t_0+\pi$, we conclude that $\Lambda'(t)$ is in the intersection of all
such wedges, which is half of the line $L_{t_0}$.
Now we can perturb $t_0$ so as to still satisfy equation
\eqref{eqn:t0}, so $\Lambda'(t)$ must lie on another nearby such line, and
these two lines intersect only at the point $\Lambda(c)$, so
$\Lambda'(t)=\Lambda(c)$.
\end{proof}

Finally, there is a sort of converse to
\fullref{prop:roundingMaximal} which characterizes the partial
order $\le$ in terms of the corner rounding operation.

\begin{proposition}
\label{prop:roundingSequence}
Let $\Lambda,\Lambda'$ be admissible paths of the same type.  Then
$\Lambda'\le \Lambda$ if and only if one can obtain $\Lambda'$ from
$\Lambda$ by a finite sequence of corner roundings.
\end{proposition}

The proof of this proposition uses
induction on the ``length'' of an admissible path.  If $\Lambda$ is an
admissible path, define its {\em length\/} $\mc{A}(\Lambda)$ to be the
sum of the lengths of its edges:
\begin{equation}
\label{eqn:defineLength}
\mc{A}(\Lambda) \eqdef \sum_{\theta\in T}m(\theta)\sqrt{x_\theta^2 +
  y_\theta^2}.
\end{equation}
Since the set of possible lengths of admissible paths (namely finite
sums of square roots of nonnegative integers) is a discrete set of
nonnegative real numbers, it is valid to perform induction on length.

\begin{lemma}
\label{lem:decreaseLength}
Rounding corners decreases length, ie if $c$ is not a kink of
$\Lambda$ then
\[
\mc{A}(\Lambda\setminus c) < \mc{A}(\Lambda).
\]
\end{lemma}

\begin{proof}
  Let $\theta_1,\theta_2$ be the two edges adjacent to $c$; since $c$
  is not a kink, $\theta_2-\theta_1\in(0,\pi]$.  The lemma is
  immediate if $\theta_2-\theta_1=\pi$, as then
\[
\mc{A}(\Lambda\setminus c) = \mc{A}(\Lambda) - 2\sqrt{x_{\theta_1}^2 +
  y_{\theta_1}^2}.
\]
So assume that $\theta_2-\theta_1\in(0,\pi)$.  Let $P$ be the path in
the plane consisting of the line segments from
$$\Lambda(c)- \begin{pmatrix} x_{\theta_1} \\ y_{\theta_1}\end{pmatrix}
\quad\text{to}\quad
\Lambda(c)
\quad\text{to}\quad
\Lambda(c)+\begin{pmatrix} x_{\theta_2} \\ y_{\theta_2}\end{pmatrix}.$$
In passing from $\Lambda$ to
$\Lambda\setminus c$, $P$ is replaced by a path $Q$, which we regard
as a continuous embedded path in the plane.  We need to show that $P$
is longer than $Q$.  Define a map
\[
f\co Q\setminus\{\mbox{corners}\} \longrightarrow P
\]
as follows.  If $x\in Q$ is not a corner, let $L_x$ be the line
through $x$ perpendicular to $Q$, and define $f(x)=L_x\cap P$.  Since
$Q$ is convex, the map $f$ is injective (although not surjective if
$Q$ has corners) and increases length.
\end{proof}

\begin{proof}[Proof of \fullref{prop:roundingSequence}]
$(\Leftarrow)$\qua This
follows immediately from \fullref{prop:roundingMaximal}(a) and
the transitivity of $\le$.

$(\Rightarrow)$\qua Suppose $\Lambda'\le\Lambda$; we need to show that
there is a sequence of corner roundings from $\Lambda$ to $\Lambda'$.
If $\Lambda'=\Lambda$ there is nothing to prove.  If
$\Lambda'\neq\Lambda$, then there is a corner $c$ of $\Lambda$ such
that $\Lambda'$ disagrees with $\Lambda$ somewhere on the interval
corresponding to $c$. By Propositions~\ref{prop:stuckAtKink} and
\ref{prop:roundingMaximal}(b), $c$ is not a kink of $\Lambda$ and
$\Lambda'\le\Lambda\setminus c$.  By \fullref{lem:decreaseLength},
$\mc{A}(\Lambda\setminus c)<\mc{A}(\Lambda)$, so by induction on
length there is a sequence of corner roundings from $\Lambda\setminus
c$ to $\Lambda'$.
\end{proof}

\section{Polygon complexes}
\label{sec:polygonComplexes}

In this section we define the combinatorial chain complexes of
interest, and we prove that the combinatorial differential $\delta$ has
degree $-1$ and satisfies $\delta^2=0$.

\subsection{The chain complex $C_*(2\pi n;\Gamma)$ and variants}
\label{sec:mgpc}

Fix a positive integer $n$ and an integer vector $\Gamma\in\Z^2$.  We
now define the combinatorial chain complex $C_*(2\pi n;\Gamma)$, and
its variants $\wwtilde{C}_*(2\pi n;\Gamma)$ and $\wwbar{C}_*(2\pi
n; \Gamma)$ which are relevant to the embedded
contact homology of $T^3$.

\subsubsection{The generators}

\begin{definition}
\label{def:generators}
$C_*(2\pi n;\Gamma)$ is the $\Z$--module generated by triples
$(\Lambda,l,o)$, where:
\begin{itemize}
\item $\Lambda$ is a periodic admissible path with rotation number $n$
  and period $\Gamma$, as defined in \fullref{sec:admissiblePaths}.  (We
  say that $(\Lambda,l,o)$ has ``underlying admissible path
  $\Lambda$''.)
\item
$l$ is a labeling of each of the edges of $\Lambda$ by `$e$' or `$h$'.
\item
$o$ is an ordering of the set of edges that are labeled `$h$'.
\end{itemize}
We impose the relations that $(\Lambda,l,o)=(\Lambda,l,o')$ if the
orderings $o$ and $o'$ differ by an even permutation, and
$(\Lambda,l,o)=-(\Lambda,l,o')$ if $o$ and $o'$ differ by an odd
permutation.  Thus $C_*(2\pi n;\Gamma)$ is a free $\Z$--module with one
generator (with no canonical sign) for each pair $(\Lambda,l)$.
\end{definition}

\subsubsection{The grading}

We now define the grading on $C_*(2\pi n;\Gamma)$.  If
$\alpha=(\Lambda,l,o)$ is a generator, consider the sum of the
multiplicities of the edges,
\[
\ell(\alpha) \eqdef \sum_{\theta\in T}m(\theta).
\]
That is, $\ell(\alpha)$ is the number of lattice points traversed by
the path $\Lambda\fedqe\Lambda(\alpha)$, counted with repetitions, but
modulo translation by $2\pi n$.  Let $\#h(\alpha)$ denote the number
of edges of $\Lambda$ that are labeled `$h$'.

\begin{definition}
\label{def:closedIndex}
If $\Gamma=0$, define the index of a generator
$\alpha$ to be
\begin{equation}
\label{eqn:closedIndex}
I(\alpha)\eqdef 2\int_{\Lambda(\alpha)}x\,dy + \ell(\alpha) - \#h(\alpha)
\in \Z.
\end{equation}
Here $\int_{\Lambda(\alpha)}x\,dy$ is the area, counted with
multiplicity, enclosed by the polygonal path $\Lambda(\alpha)$; twice
this area is an integer by Pick's formula.
\end{definition}

If $\Gamma=0$ and $n>1$ then the index can be any integer, since there
are then admissible paths with negative area.  (An example is shown in
\fullref{sec:admissiblePaths}.)  If $\Gamma\neq 0$ then there is no
canonical absolute grading, only a relative grading which is defined as
follows.

\begin{definition}
  Let $\alpha$ and $\beta$ be generators of $C_*(2\pi n;\Gamma)$.
  Choose
  $$t\in\R\setminus p^{-1}\Bigl(T(\alpha)\union T(\beta)\Bigr)$$
  and let
  $\eta$ be any path in $\R^2$ from $\Lambda(\alpha)(t)$ to
  $\Lambda(\beta)(t)$.  Let $\eta'$ denote the the translation of
  $\eta$ by $\Gamma$.  Define the composite loop
\[
P\eqdef \Lambda(\alpha)|_{[t,t+2\pi n]} + \eta'-\Lambda(\beta)|_{[t,t+2\pi n]}
-\eta. 
\]
Define the relative index
\begin{equation}
\label{eqn:relativeIndex}
I(\alpha,\beta) \eqdef 2\int_{P}x\,dy + (\ell(\alpha)-
\#h(\alpha))-(\ell(\beta)- \#h(\beta)) \in \Z.
\end{equation}
\end{definition}

Note that if $\Gamma=0$, then $I(\alpha,\beta)=I(\alpha)-I(\beta)$.
This still holds when $\Gamma\neq 0$ if one regards the index of a
single generator as taking values in an affine space over $\Z$.
There is also a canonical
mod 2 index
\[
I_2(\alpha)\eqdef \#h(\alpha)\op{mod} 2.
\]
This satisfies
$I_2(\alpha)-I_2(\beta)\equiv I(\alpha,\beta)\op{mod}2$.

\subsubsection{The differential}

\begin{definition}
\label{def:delta}
Define the differential
\[
\delta\co  C_*(2\pi n;\Gamma) \to
C_{*-1}(2\pi n;\Gamma)
\]
as follows.
  Let $\alpha$ be a generator of $C_*(2\pi n;\Gamma)$.  We define
  $\delta\alpha$ to be the signed sum of all ways to round a corner of
  $\alpha$ and ``locally lose one $h$''.  More precisely, the
  differential coefficient $\langle\delta\alpha,\beta\rangle\neq 0$ if
  and only if:
\begin{itemize}
\item
If $\alpha$ has underlying admissible path $\Lambda$, then $\beta$ has
underlying admissible path $\Lambda\setminus c$, where $c$ is a corner
of $\Lambda$ which is not a kink, ie between consecutive edges
$\theta_1,\theta_2$ with $\theta_2-\theta_1\in(0,\pi]$.
\item The edge labels of $\alpha$ and $\beta$ are identical for edges
  outside of the closed interval $[\theta_1,\theta_2]$, and either:
\begin{itemize}
\item Exactly one of the edges $\theta_1,\theta_2$ for $\alpha$ is
  labeled `$h$', and all the edges in $[\theta_1,\theta_2]$ for
  $\beta$ are labeled `$e$', or:
\item Both of the edges $\theta_1,\theta_2$ for $\alpha$ are labeled
  `$h$', and exactly one of the edges in $[\theta_1,\theta_2]$ for
  $\beta$ is labeled `$h$'.
\end{itemize}
\end{itemize}
\item
In this case $\langle\delta\alpha,\beta\rangle\eqdef\pm 1$, with 
the sign determined as in \fullref{sec:introduceDifferential}.
\end{definition}

\subsubsection{Variants}
\label{sec:variants}

The group $\Z^2$ acts on the chain complex $C_*(2\pi n;\Z)$ by
translation in the plane.  Thus there is a variant
$\wwtilde{C}_*(2\pi n;\Gamma)$, which is just $C_*(2\pi n;\Gamma)$
regarded as a $\Z[\Z^2]$--module.  There is also a variant
$\wwbar{C}_*(2\pi n;\Gamma)$ in which we mod out by translation of
polygons.  Formally,
\begin{equation}
\label{eqn:CBar}
\wwbar{C}_*(2\pi n;0)
= \wwtilde{C}_*(2\pi n;0)\tensor_{\Z[\Z^2]}\Zmodule,
\end{equation}
where $\Zmodule$ is defined in \fullref{sec:variantsIntro}.  We regard
$\wwbar{C}_*(2\pi n;\Gamma)$ as a $\Z$--module.

When
$$\Gamma=\begin{pmatrix}\Gamma_1\\ \Gamma_2\end{pmatrix}\neq 0,$$
the relative grading on $\wwbar{C}_*(2\pi n;\Gamma)$ takes values
in $\Z/2\op{gcd}(\Gamma_1,\Gamma_2)$.  The reason is that if $\alpha$
is a generator of $C_*(2\pi n;\Gamma)$ and if $\Psi_w$ denotes the
translation in the plane by a vector
$$w=\begin{pmatrix}w_1 \\
w_2\end{pmatrix}\in\Z^2,$$
then by the definition of the relative index
\eqref{eqn:relativeIndex}, we have the ``index ambiguity formula''
\begin{equation}
\label{eqn:indexAmbiguity}
I(\alpha,\Psi_w\alpha)=2\det\begin{pmatrix}
\Gamma_1&w_1\\ \Gamma_2& w_2\end{pmatrix}.
\end{equation}

\subsubsection{The auxiliary chain complex $C_*(\Lambda)$}
\label{sec:auxiliary}

We now introduce an auxiliary complex which will be used in the
computation of the homology of $C_*(2\pi n;\Gamma)$.  Let $\Lambda$ be
an admissible path of any type (see the end of
\fullref{sec:admissiblePaths}).  Define a chain complex $C_*(\Lambda)$
as follows.

\begin{definition}
  $C_*(\Lambda)$ is the $\Z$--module generated by triples
  $(\lambda,l,o)$ where $\lambda$ is an admissible path of the same
  type as $\Lambda$ with $\lambda\le\Lambda$, $l$ is a labeling of the
  edges of $\lambda$ by `$e$' or `$h$', and $o$ is an ordering of the
  `$h$' edges, with $(\lambda,l,o)=\pm(\lambda,l,o')$ as in
  \fullref{def:generators}.  The differential $\delta$ on
  $C_*(\Lambda)$ is defined just as in \fullref{def:delta}.
\end{definition}

Note that $\delta$ sends $C_*(\Lambda)$ to itself, by \fullref{prop:roundingMaximal}(a) and transitivity of $\le$.  If
$\Lambda\le \Lambda'$, then $C_*(\Lambda)$ is a subcomplex of
$C_*(\Lambda')$.  If $\Lambda$ is closed or periodic of rotation
number $n$ and period $\Gamma$, then $C_*(\Lambda)$ is a subcomplex of
$C_*(2\pi n;\Gamma)$.

The chain complex $C_*(\Lambda)$ has a relative $\Z$--grading; if
$\Lambda$ is open, this is defined by equation
\eqref{eqn:relativeIndex}, with $P$ the difference between the two
admissible paths.  If $\Lambda$ is closed, or if $\Lambda$ is open and
its initial and final endpoints agree, then $C_*(\Lambda)$ has a
canonical absolute $\Z$--grading, defined by equation
\eqref{eqn:closedIndex}.

If $\Lambda$ is an open admissible path parametrized by an interval
$I$, then the ordering of $I$ gives rise to a canonical ordering of
the `$h$' edges of any generator of $C_*(\Lambda)$, so we can regard
$C_*(\Lambda)$ as generated by pairs $(\lambda,l)$.
 
\subsubsection{The cycles $E_\Lambda$ and $H_\Lambda$}
\label{sec:distinguishedCycles}

For any admissible path $\Lambda$, the chain complex $C_*(\Lambda)$
contains two special cycles which will play a fundamental role in our
calculations.

\begin{definition}
For any admissible path $\Lambda$, define 
$E_\Lambda,H_\Lambda\in C_*(\Lambda)$ as follows:
\begin{itemize}
\item
$E_\Lambda$ is the path $\Lambda$ with all edges labeled `$e$'.
\item $H_{\Lambda}$ is the sum of all ways of taking the path
  $\Lambda$ and labeling one edge `$h$' and all other edges `$e$'.  In
  particular, if $\Lambda$ has no edges, ie if $\Lambda$ is a
  constant path, then $H_\Lambda\eqdef0$.
\end{itemize}
\end{definition}

Note that since $E_\Lambda$ is a generator with no `$h$' edge, and
$H_\Lambda$ is a sum of generators each containing exactly one `$h$'
edge, there is no choice to make here in ordering the `$h$' edges.

\begin{lemma}
\label{lem:canonicalCycles}
If $\Lambda$ is any admissible path, then $\delta E_\Lambda=\delta
H_\Lambda=0$.
\end{lemma}

\begin{proof}
$E_\Lambda$ is automatically a cycle because it is a generator with no
`$h$' edges.

To see that $H_\Lambda$ is a cycle, let $\theta_0<\cdots<\theta_k$ be
the edges of $\Lambda$, and for $0\le i \le k$ let $\alpha_i$ be the
summand in $H_\Lambda$ in which the edge $\theta_i$ is labeled `$h$'.
Let $c_i$ be the corner between $\theta_{i-1}$ and $\theta_i$.  If
$0<i<k$ and if $c_i$ and $c_{i+1}$ are not kinks, then by the
definition of $\delta$,
\[
\delta\alpha_i=E_{\Lambda\setminus c_i} - E_{\Lambda\setminus
c_{i+1}}.
\]
A modified version of this formula holds for any $0\le i \le k$, where
the first term on the right is omitted if $c_i$ is a kink or $i=0$,
and the second term on the right is omitted if $c_{i+1}$ is a kink or
$i=k$.  Thus $\delta H_\Lambda = \sum_{i=0}^k \delta\alpha_i$ consists
of two copies of $E_{\Lambda\setminus c}$ for each corner $c$ of
$\Lambda$ which is not a kink, and these two copies appear with
opposite sign.
\end{proof}

\subsubsection{Concatenation}

It will be important in our calculations to consider concatenations of
the cycles $E_\Lambda$ and $H_\Lambda$.  Suppose that $\Lambda_1$ and
$\Lambda_2$ are open admissible paths parametrized by closed intervals
$I_1$ and $I_2$ such that the right endpoint of $I_1$ (resp.\
$\Lambda_1$) agrees with the left endpoint of $I_2$ (resp.
$\Lambda_2$).  Then we can concatenate $\Lambda_1$ and $\Lambda_2$ to
obtain an open admissible path which we denote by
$\Lambda_1\Lambda_2$, which is parametrized by $I_1\cup I_2$, and
which has a corner $c$ at the concatenation point.  There is a natural
inclusion
\begin{equation}
\label{eqn:concatenation}
C_*(\Lambda_1)\tensor C_*(\Lambda_2) \longrightarrow
C_*(\Lambda_1\Lambda_2)
\end{equation}
which concatenates generators at $c$, and which we denote by
juxtaposition of symbols.  For example, $E_{\Lambda_1}E_{\Lambda_2} =
E_{\Lambda_1\Lambda_2}$.

The map \eqref{eqn:concatenation} is in general not a chain map.  If
$\alpha$ and $\beta$ are generators of $C_*(\Lambda_1)$ and
$C_*(\Lambda_2)$ respectively, then it follows from the definition of
$\delta$ that
\begin{equation}
\label{eqn:concatenationDifferential}
\delta(\alpha\beta) = (-1)^{\#h(\beta)}(\delta\alpha)\beta +
\alpha(\delta\beta) + \delta_c(\alpha\beta).
\end{equation}
Here $\#h(\beta)$ denotes the number of edges of $\beta$ labeled `$h$',
and $\delta_c$ denotes the contribution to $\delta$, if any, from
rounding at the concatenation corner $c$.

\begin{lemma}
\label{lem:concatenation}
If the concatenation corner $c$ of $\Lambda_1\Lambda_2$ is not a kink, then
in $C_*(\Lambda_1\Lambda_2)$,
\[
\begin{split}
\delta(E_{\Lambda_1}H_{\Lambda_2}) &= E_{\Lambda_1\Lambda_2\setminus
c},\\
\delta(H_{\Lambda_1}E_{\Lambda_2}) &= - E_{\Lambda_1\Lambda_2\setminus
c},\\
\delta(H_{\Lambda_1}H_{\Lambda_2}) &= H_{\Lambda_1\Lambda_2\setminus
c}.
\end{split}
\]
Moreover, for each of these equations $\delta x = y$, for each
generator $\gamma$ in $y$, there is a unique generator $\beta$ in $x$
with $\langle\delta\beta,\gamma\rangle \neq 0$. 
\end{lemma}

\begin{proof}
  By \fullref{lem:canonicalCycles} and equation
  \eqref{eqn:concatenationDifferential}, in the left hand side of each
  equation, all terms involving rounding at corners other than $c$
  cancel out.  It follows directly from the definition of $\delta$
  that rounding at $c$ gives the right hand side of each equation, and
  that the last sentence of the lemma holds.
\end{proof}

\subsection{$\delta$ has degree $-1$ and $\delta^2=0$}
\label{sec:combinatorialDifferential}

Consider the chain complex $C_*(2\pi n;\Gamma)$ or $C_*(\Lambda)$ with
its differential $\delta$.  We first check that $\delta$ has degree
$-1$.  Let $\alpha$ and $\beta$ be generators of the chain complex.

\begin{lemma}
\label{lem:pick}
Suppose the admissible path underlying $\beta$ is obtained from that
of $\alpha$ by rounding a corner.  Then
\[
I(\alpha,\beta) = 2 - \#h(\alpha) + \#h(\beta).
\]
\end{lemma}

\begin{proof}
  Let $P$ be the polygon in the plane whose oriented boundary moves
  forward along the two edges (or parts of edges) of $\alpha$ that are
  rounded to obtain $\beta$, and then backwards along the edges of
  $\beta$ that are created in the rounding process.  Let $m$ denote
  the sum of the multiplicities of these new edges in $\beta$.  By the
  definition of the relative index in equation
  \eqref{eqn:relativeIndex},
\[
I(\alpha,\beta)= 2\int_P x\,dy + 2 - \#h(\alpha) - m + \#h(\beta).
\]
Note that $P$ is a simple closed polygon (except in the extremal case
when the angle of the corner rounded in $\alpha$ is $\pi$, in which
case $P$ is a $2$--gon).  By Pick's formula, since there are no lattice
points in the interior of $P$, we have
\[
2\int_Px\,dy = m.
\]
Combining the above two equations proves the lemma.
\end{proof}

\begin{corollary}
\label{cor:deg-1}
  $\delta$ decreases degree by $1$: if
  $\langle\delta\alpha,\beta\rangle\neq 0$, then $I(\alpha,\beta)=1$. \qed
\end{corollary}

We now turn to the proof that $\delta^2=0$.  An important part of this
is to check that the roundings of an admissible path at two different
corners commute.  Let $a$ and $b$ be distinct corners of an admissible
path $\Lambda$ (of any type).  Suppose that $a$ is not a kink of
$\Lambda$ so that $\Lambda\setminus a$ is defined.  Then $b$ induces a
corner $b'$ of $\Lambda\setminus a$.  The interval corresponding to
$b'$ contains the interval corresponding to $b$.  From now on we
denote $b'$ simply by $b$.

\begin{lemma}
\label{lem:roundingCommutes}
Let $a$ and $b$ be distinct corners of an admissible path $\Lambda$.
Suppose that $\Lambda\setminus a \setminus b$ is defined (ie $a$ is
not a kink of $\Lambda$ and $b$ is not a kink of $\Lambda\setminus
a$).  Then $\Lambda\setminus b\setminus a$ is defined and
\[
\Lambda\setminus a \setminus b = \Lambda\setminus b \setminus a.
\]
\end{lemma}

\begin{proof}
  We can assume that $a$ and $b$ are consecutive and that an edge
  between them has multiplicity $1$.  (Otherwise the lemma is obvious
  as then the roundings at $a$ and $b$ involve disjoint portions of
  the domain of $\Lambda$ and thus do not affect each other.)  Without
  loss of generality, $a$ precedes $b$.
  
  By \fullref{prop:roundingMaximal}(a), $\Lambda\setminus
  a\setminus b \le \Lambda$.  Moreover, $\Lambda\setminus a \setminus
  b$ disagrees with $\Lambda\setminus a$ on all of the interval
  corresponding to $b$ in $\Lambda\setminus a$, which contains the
  interval of $\Lambda$ corresponding to $b$. So by Propositions
  \ref{prop:roundingMaximal}(b) and \ref{prop:stuckAtKink}, $b$ is not
  a kink of $\Lambda$ and
\[
\Lambda\setminus a \setminus b \le \Lambda\setminus b.
\]
Now $\Lambda\setminus a \setminus b$ disagrees with $\Lambda\setminus
b$ on an initial segment of the interval corresponding to $a$ in
$\Lambda$, hence on some of the interval corresponding to $a$ in
$\Lambda\setminus b$, so by Propositions~\ref{prop:roundingMaximal}(b)
and \ref{prop:stuckAtKink}, $a$ is not a kink of $\Lambda\setminus b$
and
\[
\Lambda\setminus a \setminus b \le \Lambda\setminus b \setminus
  a.
\]
Since $b$ is not a kink of $\Lambda$, and $a$ is not a kink of
$\Lambda\setminus b$, the mirror image of the above argument shows
that $\Lambda\setminus b \setminus a \le \Lambda\setminus a \setminus
b$.  Since $\le$ is a partial order, it follows that
$\Lambda\setminus a\setminus b = \Lambda\setminus b \setminus a$.
\end{proof}

\begin{lemma}
\label{lem:deltaSquared}
Let $\alpha,\beta,\gamma$ be generators of $C_*(2\pi n;\Gamma)$ or
$C_*(\Lambda)$ with
$\langle\delta\alpha,\beta\rangle$ and
$\langle\delta\beta,\gamma\rangle$ nonzero.  Then:
\begin{enumerate}
\item[(a)] there is a unique (up to sign) generator $\beta'\neq
\pm\beta$ with $\langle\delta\alpha,\beta'\rangle,
\langle\delta\beta',\gamma\rangle\neq 0$.
\item[(b)]
For $\beta'$ as in (a),
\[
\langle\delta\alpha,\beta\rangle\langle\delta\beta,\gamma\rangle +
\langle\delta\alpha,\beta'\rangle \langle\delta\beta',\gamma\rangle =
0.
\]
\end{enumerate}
\end{lemma}

\begin{proof}
  Let $a$ denote the corner of $\alpha$ that is rounded to obtain
  $\beta$, and let $b$ denote the corner of $\beta$ that is rounded to
  obtain $\gamma$.  We consider three cases.

  \textbf{Case 1}\qua The corner $b$ in $\beta$ comes from a corner $b$
  of $\alpha$ which is not adjacent to $a$ and thus not affected by the
  rounding at $a$.  Then $\beta'$ is obtained by performing the rounding
  at $b$ first.

  \textbf{Case 2}\qua The corner $b$ is created by the
  rounding at $a$.  (We include here the extreme case
  where $a$ and $b$ correspond to the same interval of
  length $\pi$.)  Since $\langle\delta\alpha,\beta\rangle$ and
  $\langle\delta\beta,\gamma\rangle$ are nonzero, both edges of $\alpha$
  adjacent to $a$ must be labeled `$h$', and one edge of $\beta$ adjacent
  to $b$ must be labeled `$h$'.  Then $\beta'$ is obtained from $\beta$
  by switching the labels of the edges adjacent to $b$.

  \textbf{Case 3}\qua The corner $b$ in $\beta$ comes from a corner (also
  denoted by $b$) of $\alpha$ which is adjacent to $a$.  Without loss
  of generality, $b$ comes after $a$ in $\alpha$; the case where $b$
  comes before $a$ is proved by the mirror image of the argument below.

  Any $\beta'$ must be obtained from $\alpha$ by rounding at $a$ or $b$,
  since otherwise $\beta'$ would be strictly to the left of $\alpha$
  somewhere outside the union of the intervals corresponding to $a$
  and $b$, and hence so would $\gamma$, which is a contradiction.
  Therefore this is a local problem and we may assume without loss
  of generality that the path underlying $\alpha$ is an open path
  $\Lambda_1\Lambda_2\Lambda_3$ (so we are in $C_*(\Lambda)$ for some open
  path $\Lambda$) where each $\Lambda_i$ is a single edge, the corner $a$
  is between $\Lambda_1$ and $\Lambda_2$, and the corner $b$ is between
  $\Lambda_2$ and $\Lambda_3$.  Since $\langle\delta\alpha,\beta\rangle$
  and $\langle\delta\beta,\gamma\rangle$ are nonzero, at least two of
  the three edges of $\alpha$ must be labeled `$h$'.  There are four
  ways this can happen.

The first possibility is that
\[
\alpha =
H_{\Lambda_1}H_{\Lambda_2}E_{\Lambda_3}.
\]
By Lemmas~\ref{lem:concatenation} and \ref{lem:roundingCommutes},
\[
\begin{split}
  \delta(\alpha) &= H_{\Lambda_1\Lambda_2\setminus a}E_{\Lambda_3} -
  H_{\Lambda_1}E_{\Lambda_2\Lambda_3\setminus b},\\
  \delta(H_{\Lambda_1\Lambda_2\setminus a}E_{\Lambda_3}) &=
  - E_{\Lambda_1\Lambda_2\Lambda_3\setminus a\setminus b},\\
  \delta(H_{\Lambda_1}E_{\Lambda_2\Lambda_3\setminus b}) &= -
  E_{\Lambda_1\Lambda_2\Lambda_3\setminus b\setminus a}= -
  E_{\Lambda_1\Lambda_2\Lambda_3\setminus a\setminus b}.
\end{split}
\]
These equations imply that $\delta^2\alpha=0$.  Here there is a
unique $\beta'$, namely $\pm
H_{\Lambda_1}E_{\Lambda_2\Lambda_3\setminus b}$. Likewise
$\delta^2\alpha=0$ in the other three cases
$\alpha=H_{\Lambda_1}E_{\Lambda_2}H_{\Lambda_3}$ etc.  In each case
the last sentence of \fullref{lem:concatenation} implies that
$\beta'$ is unique.
\end{proof}

\begin{corollary}
$\delta^2=0$.
\end{corollary}

\begin{proof}
\fullref{lem:deltaSquared} shows that
$\langle\delta^2\alpha,\gamma\rangle=0$ for any generators
$\alpha,\gamma$.
\end{proof}

\section{The homology operation $U$}
\label{sec:U}

We now define a degree $-2$ chain map
\[
U\co \wwtilde{C}_*(2\pi n;\Gamma)\longrightarrow
\wwtilde{C}_{*-2}(2\pi n;\Gamma).
\]
The map $U$ will not appear very often in the rest of the paper.
However we will need it in \fullref{sec:UCSS} to help compute
$\wwbar{H}_*$, and it also has a geometric counterpart in the
embedded contact homology of $T^3$ discussed in \fullref{sec:ECHU}.

Fix $\theta\in\R/2\pi n\Z$ with $\tan\theta$ irrational.  The
admissible path underlying any generator $\alpha\in C_*(2\pi
n;\Gamma)$ has a {\it{distinguished corner}} $c_\theta$, which is the
component of $(\R/2\pi n\Z)\setminus T$ containing $\theta$.  The
definition of $U$ is similar to the definition of $\delta$, but here
we preserve the number of `$h$' edges and only round at the
distinguished corner.

\begin{definition}
  For a generator $\alpha\in C_*(2\pi n;\Gamma)$, define $U_\theta\alpha$ to
  be the sum of all generators $\beta$ such that:
\begin{itemize}
\item The admissible path underlying $\beta$ is obtained from that of
  $\alpha$ by rounding the distinguished corner $c_\theta$.
\item Of the edges created or shortened by the rounding process, let
  $\{\theta_i\}$ denote those edges coming before $\theta$, and let
  $\{\theta_j'\}$ denote those edges coming after $\theta$.  If the
  edge of $\alpha$ before $c_\theta$ is labeled `$h$', then exactly
  one of the edges $\theta_i$ of $\beta$ is labeled `$h$'; otherwise
  all the edges $\theta_i$ are labeled `$e$'.  Likewise, if the edge
  of $\alpha$ after $c_\theta$ is labeled `$h$', then exactly one of
  the edges $\theta_j'$ of $\beta$ is labeled `$h$'; otherwise all the
  edges $\theta_j'$ are labeled `$e$'.
\item For all edges of $\alpha$ not adjacent to $c_\theta$, the
  corresponding edges of $\beta$ have the same labels as the
  corresponding edges of $\alpha$.  The ordering of the `$h$' edges of
  $\beta$ is induced from the ordering of the `$h$' edges of $\alpha$
  under the obvious bijection between them.
\end{itemize}
It is implicit above that $U_\theta\alpha= 0$ when $c_\theta$ is
a kink in $\alpha$.  When a fixed $\theta$ is understood, we write
$U\eqdef U_\theta$.
\end{definition}

\begin{proposition}
\label{prop:UProperties}
\begin{enumerate}
\item[(a)]
  If $\langle U\alpha,\beta\rangle\neq 0$, then $I(\alpha,\beta)=2$.
\item[(b)] $Ux = xU$, and $Uy =yU$, where $x$ and $y$ denote
  translation in the $x$-- and $y$--directions.
\item[(c)]
If the distinguished corner $c_\theta$ of $\Lambda$ is not a kink, then
\begin{equation}
\notag
U(E_\Lambda) = E_{\Lambda \setminus c_\theta}, \quad\quad 
U(H_\Lambda) = H_{\Lambda \setminus c_\theta}.
\end{equation}
\end{enumerate}
\end{proposition}

\begin{proof}
  Property (a) follows from \fullref{lem:pick}.  Properties (b) and
  (c) are immediate from the definition of $U$.
\end{proof}

\begin{proposition}
\label{U1.prop}
$U$ is a chain map: $\delta U = U \delta$.  
\end{proposition}

\begin{proof}
This will become clear after we compute $U$ of a concatenation of
$E_\Lambda$'s and $H_\Lambda$'s.  Let $\Lambda$ be a closed or
periodic admissible path of rotation number $n$ and period $\Gamma$.
Suppose that $\Lambda$ is a cyclic concatenation
\begin{equation}
\label{eqn:CC1}
\Lambda = \Lambda_1\cdots\Lambda_k
\end{equation}
where $\Lambda_i$ is an open admissible path parametrized by the
interval $[\theta_{i-1},\theta_i]$;
$\Lambda_i(\theta_i)=\Lambda_{i+1}(\theta_i)$, for $i=1,\ldots,k-1$;
$\theta_k=\theta_0+2\pi n$; and
$\Lambda_k(\theta_k)=\Lambda_0(\theta_0)+\Gamma$.  Let
$c_0,c_1,\ldots,c_k=c_0$ denote the concatenation corners, regarded as
open intervals in $\R/2\pi n\Z$.  If $c_\theta$ is not a kink of
$\Lambda$, then we can write $\Lambda\setminus c_\theta$ as a cyclic
concatenation
\begin{equation}
\label{eqn:CC2}
\Lambda\setminus c_\theta = \Lambda_1'\cdots\Lambda_k'
\end{equation}
where
\begin{equation}
\label{eqn:roundConcatenation}
\Lambda_i' \eqdef \left\{\begin{array}{cl}
\left(\Lambda_{i-1}\Lambda_i\setminus c_\theta\right)|_{[\theta,\theta_i]}, &
\theta\in c_{i-1},\\
\Lambda_i\setminus c_\theta, & \max(c_{i-1}) < \theta < \min(c_i),\\
\left(\Lambda_{i}\Lambda_{i+1}\setminus
c_\theta\right)|_{[\theta_{i-1},\theta]}, & \theta\in c_i,\\
\Lambda_i, & \mbox{otherwise}.
\end{array}\right.
\end{equation}
(Here we interpret $\Lambda_0$ and $\Lambda_{k+1}$ as appropriate
translates of $\Lambda_k$ and $\Lambda_1$ respectively.) The
decomposition \eqref{eqn:CC2} is chosen so that if the distinguished
corner $c_\theta$ agrees with a concatenation corner $c_i$ in
\eqref{eqn:CC1}, then the same is true for \eqref{eqn:CC2}; and
otherwise the concatenation corners of \eqref{eqn:CC1} and
\eqref{eqn:CC2} agree. For each $i=1,\ldots,k$, let $X_{\Lambda_i}$
denote either $E_{\Lambda_i}$ or $H_{\Lambda_i}$, and let
$X_{\Lambda_i'}$ denote $E_{\Lambda_i'}$ or $H_{\Lambda_i'}$
respectively.  Then it follows from the definition of $U$ that
\begin{equation}
\label{eqn:UConcatenation}
U\left(X_{\Lambda_1}\cdots X_{\Lambda_k}\right) =
\left\{\begin{array}{cl} 
X_{\Lambda_1'}\cdots X_{\Lambda_k'}, & \mbox{$c_\theta$ not a kink of
$\Lambda$},\\
0, & \mbox{otherwise.}
\end{array}\right.
\end{equation}

Also if $k=1$ then $\delta X_{\Lambda_1}=0$ by
\fullref{lem:canonicalCycles}; and if $k>1$, then as in
\fullref{lem:concatenation},
\begin{equation}
\label{eqn:deltaConcatenation}
\delta\left(X_{\Lambda_1}\cdots X_{\Lambda_k}\right) = \sum_{i=1}^k
\pm X_{\Lambda_1}\cdots X_{\Lambda_{i-1}}
Y_{\Lambda_i\Lambda_{i+1}\setminus c_i} X_{\Lambda_{i+2}}\cdots
X_{\Lambda_k}.
\end{equation}
Here
\[
Y_{\Lambda_i\Lambda_{i+1}\setminus c} \eqdef \left\{\begin{array}{cl}
E_{\Lambda_i\Lambda_{i+1}\setminus c}, & \mbox{if $X_{\Lambda_i} =
E_{\Lambda_i}$ and $X_{\Lambda_{i+1}} = H_{\Lambda_{i+1}}$,}\\
E_{\Lambda_i\Lambda_{i+1}\setminus c}, & \mbox{if $X_{\Lambda_i} =
H_{\Lambda_i}$ and $X_{\Lambda_{i+1}} = E_{\Lambda_{i+1}}$,}\\
H_{\Lambda_i\Lambda_{i+1}\setminus c}, & \mbox{if $X_{\Lambda_i} =
H_{\Lambda_i}$ and $X_{\Lambda_{i+1}} = H_{\Lambda_{i+1}}$,}\\
\end{array}\right.
\]
if $c_i$ is not a kink of $\Lambda_i\Lambda_{i+1}$, and
$Y_{\Lambda_i\Lambda_{i+1}\setminus c} \eqdef 0$ otherwise.

Now any generator of $C_*(2\pi n;\Gamma)$ can be written as a cyclic
concatenation $X_{\Lambda_1}\cdots X_{\Lambda_k}$ where each
$\Lambda_i$ has one edge.  Hence to prove that $\delta U=U\delta$, by
equations \eqref{eqn:UConcatenation} and
\eqref{eqn:deltaConcatenation} it suffices to show that
\[
\left(\Lambda_i'\Lambda_{i+1}'\right) \setminus c_i =
\left(\Lambda_i\Lambda_{i+1}\setminus c_i\right)',
\]
and that one side of this equation is defined if and only if the other
side is.  This follows directly from equation
\eqref{eqn:roundConcatenation} and \fullref{lem:roundingCommutes}
applied to $c_i$ and $c_\theta$.
\end{proof}

We now consider the dependence of $U_\theta$ on $\theta$.

\begin{definition}
Let $\theta_1,\theta_2\in\R/2\pi n\Z$ with $\tan \theta_1, \tan
\theta_2$ irrational.  Define a $\Z[\Z^2]$--linear map
\[
K_{\theta_1,\theta_2} \co  \wwtilde{C}_*(2\pi n;\Gamma) \longrightarrow
\wwtilde{C}_{*-1}(2\pi n;\Gamma)
\]
as follows.  If $\alpha$ is a generator of $\wwtilde{C}_*(2\pi
n;\Gamma)$, then $K_{\theta_1,\theta_2}(\alpha)$ is the sum of all
ways of relabeling an `$e$' edge of $\alpha$ in between $\theta_1$ and
$\theta_2$ by `$h$' and making this edge last in the ordering of the
`$h$' edges.
\end{definition}

\begin{proposition}
\label{prop:UChainHomotopy}
$\delta K_{\theta_1,\theta_2} + K_{\theta_1,\theta_2} \delta=
U_{\theta_1} - U_{\theta_2}$.
\end{proposition}

\begin{corollary}
\label{cor:UDefined}
The induced homomorphisms on homology
\[
\begin{split}
(U_\theta)_*\co \wwtilde{H}_*(2\pi n, \Gamma) & \longrightarrow
\wwtilde{H}_{*-2}(2\pi n; \Gamma),\\
\wwbar{H}_*(2\pi n, \Gamma) & \longrightarrow
\wwbar{H}_{*-2}(2\pi n; \Gamma)
\end{split}
\]
do not depend on the choice of $\theta\in\R/2\pi n\Z$.
\qed
\end{corollary}

\begin{proof}[Proof of \fullref{prop:UChainHomotopy}]  First note that
the statement of the proposition also makes sense in $C_*(\Lambda)$
where $\Lambda$ is an open admissible path and $\theta_1,\theta_2\in
\R$.  Here for $\theta\in\R$ and $\alpha\in C_*(\Lambda)$, it is
understood that $U_\theta\alpha=0$ if $\theta$ is not contained in any
corner of $\alpha$.  Also, if $\theta_1>\theta_2$, then
$K_{\theta_1,\theta_2}$ is understood to sum over ways of relabeling
an `$e$' edge that is greater than $\theta_1$ or less than $\theta_2$.

Now let $U_i$ and $K$ denote $U_{\theta_i}$ and
$K_{\theta_1,\theta_2}$ respectively. We want to show that
\begin{equation}
\label{eqn:KGoal}
\delta K\alpha + K
\delta\alpha=
U_{1}\alpha -
U_{2}\alpha
\end{equation}
for every generator $\alpha$ of $\wwtilde{C}_*(2\pi n;\Gamma)$.
Every term in this equation, up to edge labels, is obtained from
$\alpha$ by rounding a single corner $c$ (depending on the term).
Hence we need only check that for each corner $c$ of $\alpha$, the
contributions to both sides of the equation involving rounding at $c$
agree.  By the definition of $U$, contributions to the right hand side
involving $c$ fix all edges not adjacent to $c$.  The same is true for
the left hand side, except for contributions in which $K$ relabels an
edge that is neither adjacent to $c$ nor created by rounding at $c$;
and by our sign conventions these terms cancel in pairs.  Therefore it
is enough to prove equation \eqref{eqn:KGoal} in $C_*(\Lambda)$, where
$\Lambda$ is an open admissible path with two edges and one corner
$c$, and $\alpha$ is one of the four possible labelings of $\Lambda$.
Also we can assume that $c$ is not a kink of $\Lambda$, since
otherwise both sides of \eqref{eqn:KGoal} immediately vanish.

Without loss of generality, $\theta_1\le \theta_2$.  The reason is
that if $\theta_1>\theta_2$, then equation \eqref{eqn:KGoal} follows
by subtracting the case where $\theta_1$ and $\theta_2$ are switched
from the case where $\theta_1<\theta_2$ and the interval
$(\theta_1,\theta_2)$ contains the domain of the open path $\Lambda$.

Denote the four possibilities for $\alpha$ in the obvious manner by
$ee$, $eh$, $he$, and $hh$.  Let $\lambda\eqdef\Lambda\setminus c$,
let $\lambda^{<\theta_i}$ denote
the part of $\lambda$ up to angle $\theta_i$, let
$\lambda^{>\theta_i}$ denote the part of $\lambda$ after angle
$\theta_i$, and let $\lambda^{(\theta_1,\theta_2)}$ denote the part of
$\lambda$ between $\theta_1$ and $\theta_2$.  For example, in this
notation $\lambda =
\lambda^{<\theta_1}\lambda^{(\theta_1,\theta_2)}\lambda^{>\theta_1}$.
Let $E\eqdef E_{\lambda}$, $H\eqdef H_{\lambda}$, $H^{<\theta_i}\eqdef
H_{\lambda^{<\theta_i}}$ etc.

Suppose first that $\theta_1,\theta_2\in c$.  Then $K\alpha=0$, and
\eqref{eqn:KGoal} follows from the computations
\[\begin{split}
& U_i(ee) = E, \quad\quad\quad\quad\quad\;\;\; U_i(eh) = E^{<\theta_i}
H^{>\theta_i},\\ & U_i(he) = H^{<\theta_i} E^{>\theta_i}, \quad\quad
U_i(hh) = H^{<\theta_i}H^{>\theta_i},\\
& K\delta(ee) = K(0) = 0,\\
& K\delta(eh) = K(E) = E^{<\theta_1} H^{(\theta_1,\theta_2)} E^{>\theta_2},\\
& K\delta(he) = K(-E) = - E^{<\theta_1} H^{(\theta_1,\theta_2)}
  E^{>\theta_2},\\
& K\delta(hh) = K(H) =
H^{<\theta_1}H^{(\theta_1,\theta_2)}E^{>\theta_2} -
E^{<\theta_1}H^{(\theta_1,\theta_2)}H^{>\theta_2}.
\end{split}
\]
Note that in the last line, there are no terms with two `$h$' edges in
$\lambda$ between $\theta_1$ and $\theta_2$, because each such term
arises twice in $K(H)$ with opposite sign.

The remaining cases where $\theta_1\notin c$ and/or $\theta_2\notin c$
follow by similar, straightforward calculations.
\end{proof}

\section{Some preliminary homology calculations}
\label{sec:beginComputation}

The closed admissible paths $\Lambda$ of rotation number $n$ form a
directed set under the partial order $\le$ defined in
\fullref{sec:partialOrder}.  Part of our strategy for computing
$H_*(2\pi n;0)$ is to realize $C_*(2\pi n;0)$ as the direct limit of
the subcomplexes $C_*(\Lambda)$ defined in \fullref{sec:auxiliary},
spanned by generators to the left of a given $\Lambda$.  In this
section, as a preliminary step, we calculate most of the homology
$H_*(\Lambda)$ when $\Lambda$ is a closed admissible path of rotation
number $1$.  This homology is computed inductively using a long exact
sequence introduced below.  To carry out the induction, we will also
need to calculate the homology $H_*(\Lambda)$ for certain open
admissible paths $\Lambda$.

Throughout the homological calculations, the following decomposition
will be useful.

\begin{definition}
\label{definition:subcomplex}
Let $C_*^{(j)}(2\pi n;0)$ denote the subcomplex of $C_*(2\pi n;0)$
spanned by generators in which the index minus the number of `$h$'
edges equals $j$.  Define $C_*^{(j)}(\Lambda)$ the same way if
$\Lambda$ is a closed admissible path or an open admissible path whose
endpoints agree, so that $C_*(\Lambda)$ has a canonical $\Z$--grading.
Define $\wwtilde{C}_*^{(j)}(2\pi n;0)$ to be the corresponding
subcomplex of $\wwtilde{C}_*(2\pi n;0)$.
\end{definition}

When it should not cause confusion, we
will use the same symbols to denote both cycles and the homology
classes that they represent.

\subsection{The rounding/breaking long exact sequence}
\label{sec:exactSequence}

If $\Lambda$ is an admissible path, the homology $H_*(\Lambda)$ fits
into a long exact sequence which provides a scheme for computing it by
induction on the length of $\Lambda$.  In the statement of the exact
sequence we adopt the following:

\begin{convention}
If $c$ is a kink of $\Lambda$, so $\Lambda\setminus c$ is
undefined, we interpret $H_*(\Lambda\setminus c)\eqdef 0$.
\end{convention}

There are two versions of this exact sequence
depending on whether $\Lambda$ is open or closed/periodic.  We first
consider the case where $\Lambda$ is open.

\begin{proposition}
\label{prop:openExactSequence}
Let $\Lambda$ be an open admissible path.  Suppose $\Lambda$ has a
corner $c$ which splits it into open paths
$\Lambda_1$ and $\Lambda_2$.  Then there is a long exact sequence
\begin{equation}
\label{eqn:openExactSequence}
\cdots\to H_*(\Lambda\setminus c) \to H_*(\Lambda) \to
H_*(C_*(\Lambda_1)\tensor C_*(\Lambda_2))
\stackrel{\partial_c}{\longrightarrow} H_{*-1}(\Lambda\setminus
c)\to \cdots.
\end{equation}
The first arrow is induced by inclusion, the second arrow is induced
by projection (see below), and the map $\partial_c$ is defined by
concatenating paths and computing the part of $\delta$ involving
rounding at the corner $c$.
\end{proposition}

\begin{proof}
Suppose first that $c$ is not a kink.  Split $C_*(\Lambda)$ as the
  direct sum of two submodules where the generators are those whose
  underlying paths lie to the left of $\Lambda \setminus c$ and those
  whose paths do not.  The former submodule is the subcomplex
  $C_*(\Lambda \setminus c)$.  The latter submodule, which is not a
  subcomplex, is naturally isomorphic to $C_*(\Lambda_1)\tensor
  C_*(\Lambda_2)$, via concatenation of paths.  For
  \fullref{prop:roundingMaximal}(b) shows that any generator
  of $C_*(\Lambda)$ not in $C_*(\Lambda\setminus c)$ is obtained by
  concatenating generators of $C_*(\Lambda_1)$ and $C_*(\Lambda_2)$,
  and the concatenation operation is clearly injective.
  
  Define a differential on $C_*(\Lambda_1)\tensor C_*(\Lambda_2)$ by
  starting with the differential $\delta$ on $C_*(\Lambda)$ and
  discarding terms involving rounding at the corner $c$.  We now have
  a short exact sequence of chain complexes
\[
0\longrightarrow C_*(\Lambda\setminus c) \longrightarrow C_*(\Lambda)
\longrightarrow C_*(\Lambda_1)\tensor C_*(\Lambda_2) \longrightarrow
0.
\]
The above differential on $C_*(\Lambda_1)\tensor C_*(\Lambda_2)$ is
given more explicitly as follows.
If $\alpha$ and $\beta$ are generators of
$C_*(\Lambda_1)$ and $C_*(\Lambda_2)$ respectively, and if
$\#h(\beta)$ denotes the number of edges of $\beta$ labeled `$h$',
then by equation \eqref{eqn:concatenationDifferential},
\begin{equation}
\label{eqn:concatenationTensor}
\delta(\alpha\tensor \beta) = (-1)^{\#h(\beta)} (\delta \alpha)\tensor
\beta + \alpha \tensor \delta \beta.
\end{equation}
This is the tensor product differential on
$C_*(\Lambda_2)\tensor C_*(\Lambda_1)$, as determined by the canonical
$\Z/2$--grading on $C_*(\Lambda_i)$ given by $\#h\!\!\mod 2$.  Thus we
obtain a long exact sequence on homology as claimed.  The description
of the connecting homomorphism is immediate from the definitions.

Suppose now that $c$ is a kink.  Then by
\fullref{prop:stuckAtKink}, every generator of $C_*(\Lambda)$
is obtained by concatenating generators of $C_*(\Lambda_1)$ and
$C_*(\Lambda_2)$.  As above, equation
\eqref{eqn:concatenationDifferential} then gives
\[
H_*(\Lambda) \simeq H_*(C_*(\Lambda_1)\tensor C_*(\Lambda_2)).
\proved
\]
\end{proof}

Next suppose that $\Lambda$ is a closed or periodic admissible path
with rotation number $n$.  Let $c$ be a corner of $\Lambda$.  We can
cut $\Lambda$ at $c$ to obtain an open admissible path $\Lambda^c$
parametrized by the interval $(\theta_0,\theta_0+2\pi n)$.  Here
$\theta_0\in\R$ is a lift of a point in $\R/2\pi n\Z$ in the interval
corresponding to $c$.  Note that the two endpoints of $\Lambda^c$ will
differ by the period $\Gamma$ of $\Lambda$.  If $\Gamma\neq 0$, then
$\Lambda^c$ depends on the choice of $\theta_0$, and the different
possibilities for $\Lambda^c$ differ by translation by multiples of
$2\pi n$ in the domain and $\Gamma$ in the range.

\begin{proposition}
\label{prop:closedExactSequence}
If $\Lambda$ is a closed or periodic admissible path, if $c$ is a
corner of $\Lambda$, and if $\Lambda^c$ is an open path obtained by
cutting at $c$ as above, then there is a long exact sequence
\begin{equation}
\label{eqn:closedExactSequence}
\cdots\longrightarrow H_*(\Lambda\setminus c) \longrightarrow
H_*(\Lambda) \longrightarrow H_*(\Lambda^c)
\stackrel{\partial_c}{\longrightarrow} H_{*-1}(\Lambda\setminus c)
\longrightarrow \cdots.
\end{equation}
\end{proposition}

\begin{proof}
This is a straightforward variant of
\fullref{prop:openExactSequence}.
\end{proof}

\subsection{Convex open paths with distinct endpoints}

\begin{definition}
We say that an open admissible path $\Lambda$ is {\em convex\/} if it
is parametrized by an interval of length $\le 2\pi$ and if it
traverses a subset of the boundary of a convex polygon, possibly a
2--gon.
\end{definition}

\begin{proposition}
\label{prop:step1}
Let $\Lambda$ be a convex open admissible path with distinct
endpoints.  Then $H_*(\Lambda)$ is the free $\Z$--module generated by
the homology classes of $E_\Lambda$ and $H_\Lambda$ (see
\fullref{sec:auxiliary}).
\end{proposition}

\begin{proof}
When $\Lambda$ is a straight line, the proposition is trivial since
$E_\Lambda$ and $H_\Lambda$ are the only two generators in
$C_*(\Lambda)$, and the differential vanishes.

Now suppose that $\Lambda$ has a corner $c$.  The corner $c$ splits
$\Lambda$ into two open paths $\Lambda_1$ and $\Lambda_2$.  Observe
that $\Lambda_1$, $\Lambda_2$, and $\Lambda\setminus c$ all satisfy
the hypotheses of the proposition.  By induction on the length
$\mc{A}(\Lambda)$ and using \fullref{lem:decreaseLength}, we may
assume that the proposition holds for $\Lambda_1$, $\Lambda_2$, and
$\Lambda\setminus c$.

By \fullref{prop:openExactSequence} there is a long exact
sequence
\[
\cdots\to H_*(\Lambda\setminus c) \to H_*(\Lambda) \to
H_*(\Lambda_1)\tensor H_*(\Lambda_2)
\stackrel{\partial_c}{\longrightarrow} H_{*-1}(\Lambda\setminus
c)\to \cdots.
\]
Here we have replaced $H_*(C_*(\Lambda_1)\tensor C_*(\Lambda_2))$ by
$H_*(\Lambda_1)\tensor H_*(\Lambda_2)$, since we know by inductive
hypothesis that $H_*(\Lambda_1)$ and $H_*(\Lambda_2)$ have no torsion.
We can replace this long exact sequence by the short exact sequence
\[
0 \longrightarrow \Coker(\partial_c) \longrightarrow H_*(\Lambda)
\longrightarrow \Ker(\partial_c) \longrightarrow 0.
\]
By \fullref{lem:concatenation}, the connecting homomorphism
$\partial_c$ is given by
\begin{equation}
\begin{aligned}
\partial_c(E_{\Lambda_1}\tensor E_{\Lambda_2}) &= 0,\\
\label{eqn:connecting}
\partial_c(E_{\Lambda_1}\tensor H_{\Lambda_2}) &= -
\partial_c(H_{\Lambda_1}\tensor E_{\Lambda_2}) = E_{\Lambda\setminus c}, \\
\partial_c(H_{\Lambda_1}\tensor H_{\Lambda_2}) &= H_{\Lambda\setminus c}.
\end{aligned}
\end{equation}
Thus $\Coker(\partial_c)=0$, and $H_*(\Lambda)\simeq\Ker(\partial_c)$ is
freely generated by the homology classes of the cycles
$E_{\Lambda_1}E_{\Lambda_2} = E_\Lambda$ and
$E_{\Lambda_1}H_{\Lambda_2} + H_{\Lambda_1}E_{\Lambda_2} = H_\Lambda$.
\end{proof}

\subsection{Convex open paths with identical endpoints}

Our next task is to compute $H_*(\Lambda)$ where $\Lambda$ is a
(non-constant) convex open admissible path whose two endpoints are the
same lattice point $a$, so that $\Lambda$ traverses all of the
boundary of a convex polygonal region $P_\Lambda$. Let $k_\Lambda$
denote the number of lattice points in $P_\Lambda$.

In this case a new homology generator appears.  To describe it, let
$b\neq a$ be a lattice point in $P_\Lambda$, and suppose that
$b-a\in\Z^2$ is indivisible.
Define
\[
Z_1(a,b) \in C_*(\Lambda)
\]
to be the 2--gon from $a$ to $b$ and back with both edges labeled
`$h$', and with the edge from $a$ to $b$ first in the ordering.  By
the definition of $\delta$, we have $\delta Z_1(a,b) = 0$.  It is not
hard to show that for the open path $\Lambda$, the homology class of
$Z_1(a,b)$ in $C_*(\Lambda)$ does not depend on $b$, although we will
not need this.

\begin{proposition}
\label{prop:step2}
  Let $\Lambda$ be a nonconstant convex open admissible path whose two
  endpoints are the same lattice point $a$.  Then $H_*(\Lambda)$ is
  freely generated by the homology classes of $E_\Lambda$,
  $H_\Lambda$, and $Z_1(a,b)$ for a single $b$ with $b-a$ indivisible.
\end{proposition}

\begin{proof}
We use induction on $k_\Lambda$, the number of lattice points in $P_\Lambda$.

If $k_\Lambda=2$, then the chain complex $C_*(\Lambda)$ has five
generators: four 2--gons which we denote in the obvious manner by $ee$,
$eh$, $he$, and $hh$; and the constant path at $a$ which we denote by
$p(a)$.  The differential is given by $\delta(ee)=0$,
$\delta(eh)=-\delta(he)=p(a)$, $\delta(hh)=0$, and $\delta(p(a))=0$.
The proposition in this case follows by inspection.

If $k_\Lambda>2$, let $c$ be a corner of $\Lambda$, which cuts
$\Lambda$ into convex open paths $\Lambda_1$ and $\Lambda_2$ with
distinct endpoints.  We may inductively assume that the proposition
holds for $\Lambda\setminus c$.  By \fullref{prop:step1},
$H_*(\Lambda_i)$ is freely generated by the homology classes of
$E_{\Lambda_i}$ and $H_{\Lambda_i}$.  By
\fullref{prop:openExactSequence}, there is a long exact
sequence
\[
\cdots\to H_*(\Lambda\setminus c) \to H_*(\Lambda) \to
H_*(\Lambda_1)\tensor H_*(\Lambda_2)
\stackrel{\partial_c}{\longrightarrow} H_{*-1}(\Lambda\setminus
c)\to \cdots.
\]
Here we have used the fact that $H_*(\Lambda_i)$ is free to commute
homology and tensor product in the third term.  The above long exact
sequence gives rise to a short exact sequence
\[
0 \longrightarrow \Coker(\partial_c) \longrightarrow H_*(\Lambda)
\longrightarrow \Ker(\partial_c) \longrightarrow 0.
\]
The connecting homomorphism $\partial_c$ is again computed by the
equations \eqref{eqn:connecting}.  Thus $\Ker(\partial_c)$ is freely
generated by the homology classes of $E_\Lambda$ and $H_\Lambda$ as
before, while now $\Coker(\partial_c)$ is freely generated by the
homology class of $Z_1(a,b)$.
\end{proof}

\subsection{Closed paths of rotation number 1}
\label{sec:step3}

We now compute most of $H_*(\Lambda)$ where $\Lambda$ is a closed
admissible path with rotation number 1.  Again, $\Lambda$ traverses
the boundary of a convex polygonal region $P_\Lambda$; let $k_\Lambda$
denote the number of lattice points in $P_\Lambda$.

If $\lambda\le \Lambda$ is another closed admissible path with
rotation number 1, then $E_\lambda$ and $H_\lambda$ are cycles in
$C_*(\Lambda)$.  (Recall that when $\lambda$ is a constant path,
$H_\lambda$ is defined to be zero.)  We have the following relations
between these cycles.

\begin{lemma}
\label{lem:convenient}
  Let $\Lambda$ be a closed admissible path with rotation number 1 and
  let $\lambda_1,\lambda_2\le \Lambda$ with
  $k_{\lambda_1}=k_{\lambda_2}=k\ge 1$.  Then in the homology
  $H_*(\Lambda)$ we have
\begin{align}
\label{eqn:ERelation}
E_{\lambda_1} = E_{\lambda_2},\\
\label{eqn:HRelation}
H_{\lambda_1} = H_{\lambda_2}.
\end{align}
\end{lemma}

\begin{proof}
  If $\lambda_1=\Lambda$ then necessarily $\lambda_2=\Lambda$ and the
  lemma is trivial.  If $\lambda_1\neq\Lambda$ then by
  \fullref{prop:roundingSequence} there exist corners
  $a_1,a_2$ of $\Lambda$ with $\lambda_i\le\Lambda\setminus a_i$.  By
  induction on $k_\Lambda$, we may assume that the lemma holds for
  $\Lambda\setminus a_i$.  Consequently, it is enough to choose
  convenient $\lambda_i\le\Lambda\setminus a_i$ with $k_{\lambda_i}=k$
  and verify the relations \eqref{eqn:ERelation} and
  \eqref{eqn:HRelation}.  Without loss of generality, $a_1$ and $a_2$
  are consecutive corners of $\Lambda$ with $a_1$ coming first.
  
  Choose any $\mu\le\Lambda$ with $k_{\mu}=k+1$.  The path $\mu$ has
  at most one edge $\theta$ whose initial corner maps to the same
  lattice point as $a_1$ and/or whose final corner maps to the same
  lattice point as $a_2$.  If $\mu$ has no such edge, let $\theta$ be
  an arbitrary edge of $\mu$.  Let $b_1$ and $b_2$ denote the corners
  of $\mu$ before and after the edge $\theta$.  Since $\mu$ is a
  closed admissible path of rotation number 1, it contains no kinks.
  Choose $\lambda_i \eqdef \mu\setminus b_i$; then $\lambda_i \le
  \Lambda\setminus a_i$.
  
  Let $\alpha$ be the generator of $C_*(\Lambda)$ with underlying path
  $\mu$ with the edge $\theta$ labeled `$h$' and all other edges
  labeled `$e$'.  Let $\beta\in C_*(\Lambda)$ be the chain obtained
  from $\alpha$ by summing over all ways of relabeling one of the
  `$e$' edges by `$h$' and ordering it first.  Then similarly to
  \fullref{lem:concatenation},
\[
\begin{split}
\delta\alpha &= E_{\lambda_1} - E_{\lambda_2},\\
\delta\beta &= H_{\lambda_1} - H_{\lambda_2}.
\end{split}
\]
Hence \eqref{eqn:ERelation} and \eqref{eqn:HRelation} hold in
homology.
\end{proof}

Note that by the index formula \eqref{eqn:index1}, the indices of the
above generators are
\begin{equation}
\label{eqn:IC}
\begin{split}
I(E_\lambda)&=2(k_\lambda-1),\\
I(H_\lambda)&=2(k_\lambda-1)-1,\quad\quad k_\lambda>1.
\end{split}
\end{equation}

\begin{proposition}
\label{prop:step3}
Let $\Lambda$ be a closed admissible path of rotation number 1.  Then:
\begin{enumerate}
\item[(a)] $H_i^{(-2)}(\Lambda)=0$ for $i\neq 0$.  (See
\fullref{definition:subcomplex}.)
\item[(b)]
\[
\bigoplus_{j\neq-2}H_i^{(j)}(\Lambda) \simeq \left\{\begin{array}{cl} \Z,
& 0\le i \le 2(k_\Lambda-1),\\
0, & \mbox{otherwise.}
\end{array}
\right.
\]
This is generated
by the homology classes of $E_\lambda$ for $\lambda\le\Lambda$ and
$H_\lambda$ for $\lambda\le\Lambda$ nonconstant, with the relations
\eqref{eqn:ERelation} and \eqref{eqn:HRelation}.
\end{enumerate}
\end{proposition}

\begin{proof}
  We use induction on $k_\Lambda$.  If $k_\Lambda=1$ then
  $C_*(\Lambda)$ has only the single generator $E_\Lambda$ and the
  result is immediate.  Suppose $k_\Lambda>1$, let $c$ be a corner
  of $\Lambda$, and assume that the proposition holds for
  $\Lambda\setminus c$.  By \fullref{prop:closedExactSequence}
  there is a long exact sequence in homology
\begin{equation}
\label{eqn:step3ExactSequence}
\cdots\longrightarrow H_*(\Lambda\setminus c) \longrightarrow
H_*(\Lambda) \longrightarrow H_*(\Lambda^c)
\stackrel{\partial_c}{\longrightarrow} H_{*-1}(\Lambda\setminus c)
\longrightarrow \cdots.
\end{equation}
By \fullref{prop:step2}, $H_*(\Lambda^c)$ is freely generated
by $E_{\Lambda^c}$, $H_{\Lambda^c}$, and $Z_1(c,b)$ for a single $b$
with $b-c$ indivisible.  By regarding the open paths as closed paths,
these cycles lift to cycles in $C_*(\Lambda)$, namely $E_\Lambda$,
$H_\Lambda$, and a third cycle which we will also call $Z_1(c,b)$.
Hence the connecting homomorphism $\partial_c$ vanishes, and we have a
short exact sequence
\begin{equation}
\label{eqn:SES}
0 \longrightarrow H_*^{(j)}(\Lambda\setminus c) \longrightarrow
H_*^{(j)}(\Lambda)
\longrightarrow H_*^{(j)}(\Lambda^c) \longrightarrow 0.
\end{equation}

The index calculation \eqref{eqn:IC} implies that $E_{\Lambda}$ and
$H_{\Lambda}$ have $j\ge 0$, while $Z_1(c,b)$ has $j=-2$.  Also, the
definition of the index implies that $I(Z_1(c,b))=0$.  Thus part (a)
of the proposition follows immediately from the exact sequence
\eqref{eqn:SES} with $j=-2$.  For $j\neq -2$, the exact sequence
\eqref{eqn:SES} implies that the inclusion-induced map
\[
\Z\{E_\Lambda,H_\Lambda\}
\oplus
\bigoplus_{j\neq-2}H_*^{(j)}(\Lambda\setminus c)
\longrightarrow
\bigoplus_{j\neq
-2}H_*^{(j)}(\Lambda)
\]
is an isomorphism.  Part (b) of the proposition follows from this and
\fullref{lem:convenient}.
\end{proof}

At this point it is not hard to compute the rest of $H_*(\Lambda)$,
namely $H_0^{(-2)}(\Lambda)$, and also to take the direct limit over
$\Lambda$ to recover $H_*(2\pi;0)$ and $\wwtilde{H}_*(2\pi;0)$ as
described in \fullref{sec:combinatoricsIntro}.  (The notation $Z_1(c,b)$
here corresponds to $h(c,b)$ in \fullref{sec:combinatoricsIntro}.)  We
will do these calculations in greater generality in \fullref{sec:zero}
and \fullref{sec:HTildenProof} respectively.

\section{Flattening and applications}

The previous section did most of the calculation of
$\wwtilde{H}_*(2\pi n;0)$ when $n=1$, which entailed calculating
most of $H_*(\Lambda)$ where $\Lambda$ is a closed admissible path of
rotation number $1$.  In general $H_*(\Lambda)$ is much more
complicated when $\Lambda$ has rotation number $n>1$.  To simplify the
calculations for arbitrary $n$, in \fullref{sec:flattening} we introduce
a ``flattening'' technique, exemplified by
\fullref{prop:flattening} below, which reduces many
homological calculations to more manageable ones involving admissible
paths on the $x$--axis.  As a first application, in \fullref{sec:zero} we
compute $H_0^{(-2n)}(\Lambda)$ where $\Lambda$ is $n$--convex, which
means the following:

\begin{definition}
An {\em $n$--convex path\/} is a closed admissible path $\Lambda$ of
  rotation number $n$ which is the pullback, via the projection
  $\R/2\pi n\Z \to \R/2\pi \Z$, of a closed admissible path of
  rotation number $1$.  By \fullref{ex:simplestPath}, the latter
  corresponds to a convex polygonal region $P_\Lambda$ in $\R^2$ with
  corners in $\Z^2$, possibly a $2$--gon or a point.  Let $k_\Lambda$
  denote the number of lattice points in $P_\Lambda$.
\end{definition}

As a second application of the flattening technology, we will prove in
\fullref{sec:vanishing} that the homology vanishes for periodic paths
that are not closed:

\begin{theorem}
\label{thm:vanishing}
For any $n$, if $\Gamma\neq 0$, then
\[
H_*(2\pi n;\Gamma)=0.
\]
\end{theorem}

Note that this is equivalent to $\wwtilde{H}_*(2\pi n;\Gamma)=0$,
see \fullref{sec:variants}.  By the universal coefficient spectral
sequence (see \fullref{sec:UCSS}), this implies
\fullref{thm:main}(a).

\subsection{Flattening}
\label{sec:flattening}

Let $\Lambda$ be a closed admissible path of rotation number $n$, and
fix an angle $\theta\in(0,2\pi)$ with $\tan\theta$ irrational.  We now
define a subcomplex $C_*^\theta(\Lambda)$ of $C_*(\Lambda)$ which is
much smaller but has the same homology.  The first step is to define a
decomposition of the $\Z$--module $C_*(\Lambda)$ in terms of
``$\theta$--corner sequences''.

\begin{definition}
For $\lambda\le\Lambda$, the {\em $\theta$--corner sequence\/}
associated to $\lambda$ is the sequence $(p_0,p_1,\ldots,p_{2n}=p_0)$
of points in $\Z^2$ defined by
\[
p_i \eqdef \lambda(\theta + i\pi).
\]
\end{definition}

Let $S^\theta(\Lambda)$ denote the set of all $\theta$--corner
sequences that can arise for admissible paths $\lambda\le\Lambda$.
The set $S^\theta(\Lambda)$ can be characterized as follows.  Define
an ordering on $\Z^2$ by
\begin{equation}
\label{eqn:ordering}
p\le q \Longleftrightarrow 
\det\left(\begin{matrix}\begin{matrix}
\cos \theta \\
\sin \theta
\end{matrix}
& p-q
\end{matrix}
\right) \ge 0.
\end{equation}
If $p\le q$ and $i$ is odd, or $p\ge q$ and $i$ is even, let
$s^i_{p,q}$ denote the open admissible path parametrized by the
interval
\[
I_i\eqdef (\theta+i\pi, \theta+(i+1)\pi)
\]
with endpoints at $p$ and $q$ and with a single edge if $p\neq q$ and
no edges otherwise.  If $\lambda'$ and $\lambda$ are open admissible
paths parametrized by the same interval, but (unlike in
\fullref{def:left}) not necessarily having the same endpoints,
write
\[
\lambda' \le' \lambda
\]
if the inequality
\eqref{eqn:left} holds.

\begin{lemma}
\label{lem:cornerSequence}
$(p_0,p_1,\ldots,p_{2n}=p_0)\in S^\theta(\Lambda)$ if and only if for
each $i$, letting $p\eqdef p_i$ and $q\eqdef p_{i+1}$, we have
\begin{enumerate}
\item[(i)] $p \ge q$ if $i$ is even, and $p \le q$ if $i$ is odd.
\item[(ii)] $s^i_{p,q} \le' \Lambda|_{I_i}$.
\end{enumerate}
\end{lemma}

\begin{proof}
  $(\Rightarrow)$ If $(p_0,\ldots,p_{2n})$ is the $\theta$--corner
  sequence associated to $\lambda\le\Lambda$, then condition (i) holds
  since $\lambda$ is an admissible path, and condition (ii) holds
  because $s^i_{p_i,p_{i+1}}\le \lambda|_{I_i} \le' \Lambda|_{I_i}$.
  $(\Leftarrow)$ If $(p_0,\ldots,p_{2n})$ satisfies conditions (i) and
  (ii), then $\lambda=s^0_{p_0,p_1}\cdots s^{2n-1}_{p_{2n-1},p_{2n}}$
  is well defined by (i), satisfies $\lambda\le\Lambda$ by (ii), and
  has $(p_0,\ldots,p_{2n})$ as its $\theta$--corner sequence.
\end{proof}

For a given $i$, if $p$ and $q$ satisfy conditions (i) and (ii) above,
then there is a unique maximal open admissible path $\Lambda^i_{p,q}$
parametrized by $I_i$ with endpoints at $p$ and $q$ and with
$\Lambda^i_{p,q}\le' \Lambda|_{I_i}$.  The path $\Lambda^i_{p,q}$
traverses part of the boundary of the convex hull of the set of
lattice points enclosed by $s^i_{p,q}$, $\Lambda|_{I_i}$, and rays
from $p$ and from $q$ in the direction
$$\begin{pmatrix}\cos\theta\\\sin\theta\end{pmatrix}
\quad\text{if $i$ is even, and}\quad
\begin{pmatrix}-\cos\theta\\-\sin\theta\end{pmatrix}
\quad\text{if $i$ is odd.}$$

\begin{example}
The picture below shows an example where $n=1$ and $\theta$ is
slightly greater than zero.  Here the outer polygon is $\Lambda$, the
numbers indicate the ordering \eqref{eqn:ordering} of the $10$ lattice
points in $P_{\Lambda}$, and the inner path is $\Lambda^0_{8,3}$.
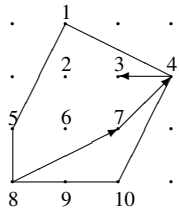
\begin{center}
\begin{picture}(77,80)(-8,-12)
\put(0,0){.}
\put(20,0){.}
\put(40,0){.}
\put(60,0){.}
\put(0,20){.}
\put(20,20){.}
\put(40,20){.}
\put(60,20){.}
\put(0,40){.}
\put(20,40){.}
\put(40,40){.}
\put(60,40){.}
\put(0,60){.}
\put(20,60){.}
\put(40,60){.}
\put(60,60){.}
\put(1.5,0.5){\line(1,0){40}}
\put(41.5,0.5){\line(1,2){20}}
\put(1.5,0.5){\line(0,1){20}}
\put(1.5,20.5){\line(1,2){20}}
\put(21.5,60.5){\line(2,-1){40}}
\put(19,64){${}_1$}
\put(19,44){${}_2$}
\put(39,44){${}_3$}
\put(59,44){${}_4$}
\put(-1,24){${}_5$}
\put(19,24){${}_6$}
\put(39,24){${}_7$}
\put(-1,-7){${}_8$}
\put(19,-7){${}_9$}
\put(39,-7){${}_{10}$}
\put(1.5,0.5){\vector(2,1){40}}
\put(41.5,20.5){\vector(1,1){20}}
\put(61.5,40.5){\vector(-1,0){20}}
\end{picture}
\end{center}
\end{example}

\begin{lemma}
As $\Z$--modules,
\begin{equation}
\label{eqn:flatteningDecomposition}
C_*(\Lambda) = \bigoplus_{(p_0,\ldots,p_{2n})\in S^\theta(\Lambda)}
\bigotimes_{i=0}^{2n-1} C_*\left(\Lambda^i_{p_i,p_{i+1}}\right).
\end{equation}
\end{lemma}

\begin{proof}
It is enough to show that cyclic concatenation induces a bijection
\[
\left\{\lambda\le\Lambda\right\}
=
\bigsqcup_{(p_0,\ldots,p_{2n})\in S^\theta(\Lambda)}
\prod_{i=0}^{2n-1} \left\{\lambda_i\le\Lambda^i_{p_i,p_{i+1}}\right\}.
\]
If $(p_0,\ldots,p_{2n})\in S^\theta(\Lambda)$ and $\lambda_i\le
\Lambda^i_{p_i,p_{i+1}}$ for $i=0,\ldots,2n-1$, then the cyclic
concatenation $\lambda=\lambda_0\lambda_1\cdots \lambda_{2n-1}$
satisfies $\lambda\le \Lambda$ and has $(p_0,\ldots,p_{2n})$ as its
$\theta$--corner sequence.  Conversely, any $\lambda\le\Lambda$ is
obtained this way, where $(p_0,\ldots,p_{2n})$ is $\lambda$'s
$\theta$--corner sequence and $\lambda_i=\lambda|_{I_i}$.
\end{proof}

\begin{lemma}
\label{lem:TPO}
For a given $i$, suppose that the pairs $(p,q)$ and $(p',q')$ satisfy
conditions (i) and (ii) above.  Suppose that under the ordering
\eqref{eqn:ordering}, the interval between $p$ and $q$ contains the
interval between $p'$ and $q'$.  Then
\begin{equation}
\label{eqn:TPO}
\Lambda^i_{p',q'} \le' \Lambda^i_{p,q}.
\end{equation}
\end{lemma}

\begin{proof}
  Since $\Lambda|_{I_i}$ is parametrized by an interval of length
  $\pi$, we can find a convex polygonal region $P$ such that
  $\Lambda|_{I_i}$ is part of the boundary of $P$.  Then the right
  side of \eqref{eqn:TPO} is part of the boundary of the convex hull
  of the set of lattice points in $P$ that are between $p$ and $q$ in
  the ordering \eqref{eqn:ordering}, inclusive.  The left side of
  \eqref{eqn:TPO} has an analogous description for a subset of these
  lattice points.  The relation \eqref{eqn:TPO} then follows as in
  \fullref{example:leftInclusion}.
\end{proof}

\begin{definition}
Let $E^i_{p,q}$ and $H^i_{p,q}$ respectively denote
$E_{\Lambda^i_{p,q}}, H_{\Lambda^i_{p,q}}\in C_*(\Lambda^i_{p,q})$,
see \fullref{sec:distinguishedCycles}.
In terms of the decomposition \eqref{eqn:flatteningDecomposition},
define
\[
C_*^\theta(\Lambda) \eqdef \bigoplus_{(p_0,\ldots,p_{2n}) \in
S^\theta(\Lambda)} \bigotimes_{i=0}^{2n-1}
\op{span}\left\{E^i_{p_i,p_{i+1}},H^i_{p_i,p_{i+1}}\right\} \subset
C_*(\Lambda).
\]
\end{definition}

\begin{lemma}
$C_*^\theta(\Lambda)$ is a subcomplex of $C_*(\Lambda)$.
\end{lemma}

\begin{proof}
  For $(p_0,\ldots,p_{2n})\in S^\theta(\Lambda)$, define the cyclic
  concatenation
\[
\Lambda(p_0,p_1,\ldots,p_{2n}) \eqdef \Lambda^0_{p_0,p_1} \cdots
  \Lambda^{2n-1}_{p_{2n-1},p_{2n}}.
\]
Recall from \fullref{sec:U} that $c_{\theta+i\pi}$ denotes the corner of
$\Lambda$ containing $\theta+i\pi$.  By \fullref{lem:concatenation},
it is enough to show that if $(p_0,\ldots,p_{2n})\in
S^\theta(\Lambda)$ and $c_{\theta+i\pi}$ is not a kink of
$\Lambda(p_0,p_1,\ldots,p_{2n})$, then
\[
\Lambda(p_0,\ldots,p_{2n})
\setminus
c_{\theta+i\pi} =
\Lambda(p_0,\ldots,p_{i-1},p_i',p_{i+1},\ldots,p_{2n})
\]
where
\[
p_i' \eqdef \left(\Lambda(p_0,\ldots,p_{2n})\setminus c_{\theta +
i\pi}\right) (\theta+i\pi).
\]
By locality of rounding it is enough to show that
\begin{equation}
\label{eqn:LOR}
\Lambda^{i-1}_{p_{i-1},p_i} \Lambda^i_{p_i,p_{i+1}} \setminus
c_{\theta+i\pi} = \Lambda^{i-1}_{p_{i-1},p_i'}
\Lambda^i_{p_i',p_{i+1}}.
\end{equation}
In \eqref{eqn:LOR} we have $\le$ by definition, and $\ge$ by
\fullref{prop:roundingMaximal}(b) and \fullref{lem:TPO},
hence $=$ by \fullref{prop:partialOrder}.
\end{proof}

\begin{example}
\label{example:xtheta}
If $\Lambda_0$ is on the $x$--axis, then $C_*^\theta(\Lambda_0) =
C_*(\Lambda_0)$.  Here $E^i_{p,q}$
and $H^i_{p,q}$ are single edges labeled `$e$' and `$h$' respectively
(when $p\neq q$).
\end{example}

\begin{proposition}
\label{prop:III}
For any closed admissible path $\Lambda$ of rotation number $n$, the
inclusion $C_*^\theta(\Lambda)\to C_*(\Lambda)$ induces an isomorphism
on homology.
\end{proposition}

\begin{proof}
By \fullref{prop:roundingSequence} and induction on length,
there are only finitely many admissible paths to the left of
$\Lambda$.  Hence only finitely many points in $\Z^2$ can appear in a
sequence in $S^\theta(\Lambda)$; denote these by $q_1,\ldots,q_k$, in
increasing order with respect to the ordering \eqref{eqn:ordering}.
If $\alpha$ is a generator of $C_*(\Lambda)$ with $\theta$--corner
sequence $(q_{j_0}, \ldots, q_{j_{2n}})\in S^\theta(\Lambda)$, define
the ``degree''
\[
\deg(\alpha) \eqdef \sum_{i=0}^{2n-1} (-1)^i j_i \ge 0.
\]
This defines an increasing filtration on $C_*(\Lambda)$.  Indeed the
differential $\delta$ on $C_*(\Lambda)$ splits as
\[
\delta = \delta_0 + \delta_1
\]
where $\delta_1$ is the contribution from rounding at the corners
$c_{\theta+i\pi}$; then $\delta_0$ preserves the degree while
$\delta_1$ decreases it (possibly by an arbitrarily large amount).
This filtration yields a spectral sequence $E^*_{*,*}$ converging to
$H_*(\Lambda)$.  The $E^1$ term is the homology of $\delta_0$.  By
equation \eqref{eqn:concatenationTensor}, the differential $\delta_0$
agrees with the tensor product differential on the right hand side of
\eqref{eqn:flatteningDecomposition}, with the standard sign if the
factors in the tensor product are arranged in the order
$i=2n-1,\ldots,0$.  So by \fullref{prop:step1}, the inclusion
$C_*^\theta(\Lambda) \to C_*(\Lambda)$ induces an isomorphism
\begin{equation}
\label{eqn:flattening1}
C_*^{\theta}(\Lambda) = E^1.
\end{equation}

We can filter the subcomplex $C_*^\theta(\Lambda)$ the same way to
obtain a spectral sequence $'E^*_{*,*}$ converging to the homology of
$C_*^\theta(\Lambda)$.  Now $\delta_0$ restricts to zero on
$C_*^\theta(\Lambda)$ by \fullref{lem:canonicalCycles}, so
\begin{equation}
\label{eqn:flattening2}
'E^1=C_*^{\theta}(\Lambda).
\end{equation}

Putting this together, the inclusion of filtered complexes
$C_*^\theta(\Lambda)\to C_*(\Lambda)$ induces a morphism of spectral
sequences, which by equations \eqref{eqn:flattening1} and
\eqref{eqn:flattening2} induces an isomorphism $'E^1\simeq E^1$.  Hence the
inclusion induces an isomorphism $'E^\infty\simeq E^\infty$, and
therefore an isomorphism on homology.
\end{proof}

We now specialize to the case where $\Lambda$ is $n$--convex.

\begin{lemma}
\label{lem:CICC}
Let $\Lambda, \Lambda'$ be $n$--convex paths with $k_\Lambda =
k_{\Lambda'}$.  Then there is a canonical isomorphism of $\Z$--graded
chain complexes
\[
C_*^\theta(\Lambda) = C_*^\theta(\Lambda').
\]
\end{lemma}

\begin{proof}
  We compute the chain complex $C_*^\theta(\Lambda)$ explicitly.
  Denote the lattice points enclosed by $\Lambda$ in increasing order
  with respect to the ordering \eqref{eqn:ordering} by $1,\ldots,k$.
  By \fullref{lem:cornerSequence}, $(j_0,\ldots,j_{2n}) \in
  S^\theta(\Lambda)$ if and only if $j_i \ge j_{i+1}$ for $i$ even and
  $j_i \le j_{i+1}$ for $i$ odd.  Moreover if $i$ is even and $j,j'' >
  j'$, then $c_{\theta + i\pi}$ is not a kink of $\Lambda^{i-1}_{j,j'}
  \Lambda^i_{j',j''}$, and
\[
\Lambda^{i-1}_{j,j'} \Lambda^i_{j',j''} \setminus c_{\theta+i\pi} =
\Lambda^{i-1}_{j,j'+1} \Lambda^i_{j'+1,j''}.
\]
Here $\ge$ holds and $c_{\theta+i\pi}$ is not a kink by
\fullref{lem:TPO} and Propositions \ref{prop:roundingMaximal}(b) and
\ref{prop:stuckAtKink}, while $\le$ holds by equation \eqref{eqn:LOR}
and \fullref{lem:TPO}.  Likewise, if $i$ is odd and $j,j'' < j'$
then
\[
\Lambda^{i-1}_{j,j'} \Lambda^i_{j',j''} \setminus c_{\theta+i\pi} =
\Lambda^{i-1}_{j,j'-1} \Lambda^i_{j'-1,j''}.
\]
So by \fullref{lem:concatenation}, the differential on
$C_*^\theta(\Lambda)$ operates on a length $2n$ cyclic string of $E$'s
and $H$'s according to the local (up to sign) rules
\begin{equation}
\label{eqn:d1}
\begin{split}
E^{i-1}_{j,j'}H^i_{j',j''},\;\, H^{i-1}_{j,j'}E^i_{j,j''} &\longmapsto
\pm E^{i-1}_{j,j'+1}E^i_{j'+1,j''},\\ H^{i-1}_{j,j'}H^i_{j',j''} &
\longmapsto \pm\left(E^{i-1}_{j,j'+1}H^i_{j'+1,j''} +
H^{i-1}_{j,j'+1}E^i_{j'+1,j''}\right),
\end{split}
\end{equation}
for $i$ even and $j,j'' > j'$, and similarly for $i$ odd.

The important point is that the above description of the chain complex
$C_*^\theta(\Lambda)$ depends only on $k$.  Thus we get a canonical
isomorphism of chain complexes $C_*^\theta(\Lambda) =
C_*^\theta(\Lambda')$.  This isomorphism respects the grading, as one
can see by using rounding and relabeling to inductively reduce to
generators involving constant corner sequences.
\end{proof}

\begin{definition}
\label{def:flattening}
  Let $\Lambda,\Lambda_0$ be $n$--convex paths with $k_\Lambda =
  k_{\Lambda_0}$ and with $P_{\Lambda_0}$ on the $x$--axis.  Define a
  chain map
\[
F_\theta\co  C_*(\Lambda_0) \longrightarrow C_*(\Lambda)
\]
to be the composition of canonical isomorphisms and inclusion
\[
C_*(\Lambda_0) = C_*^\theta(\Lambda_0) = C_*^\theta(\Lambda)
\longrightarrow C_*(\Lambda).
\]
\end{definition}

Explicitly, $F_\theta$ takes a generator of $C_*(\Lambda_0)$ and
replaces each `$e$' or `$h$' edge by a corresponding $E^i_{p,q}$ or
$H^i_{p,q}$ to obtain a cyclic concatenation of such in
$C_*^\theta(\Lambda)\subset C_*(\Lambda)$.

\begin{proposition}
\label{prop:flattening}
Let $\Lambda,\Lambda_0$ be $n$--convex paths with $k_\Lambda =
k_{\Lambda_0}$ and with $P_{\Lambda_0}$ on the $x$--axis.  Then:
\begin{enumerate}
\item[(a)]
The chain map $F_\theta\co C_*(\Lambda_0)\to C_*(\Lambda)$ induces an
isomorphism on homology
\[
(F_\theta)_*\co H_i^{(j)}(\Lambda_0)\stackrel{\simeq}{\longrightarrow}
H_i^{(j)}(\Lambda).
\]
\item[(b)] Let $\Lambda'$ and $\Lambda_0'$ be the $n$--convex paths
  obtained by rounding the distinguished corners $c_\theta$ of
  $P_\Lambda$ and $P_{\Lambda_0}$ respectively.  Then the diagram
\[
\begin{CD}
H_*^{(j)}(\Lambda_0') @>{(F_\theta)_*}>> H_*^{(j)}(\Lambda') \\
@VVV @VVV \\
H_*^{(j)}(\Lambda_0) @>{(F_\theta)_*}>> H_*^{(j)}(\Lambda)
\end{CD}
\]
commutes, where the vertical arrows are induced by inclusion.
\end{enumerate}
\end{proposition}

\begin{proof}
  Part (a) follows from \fullref{prop:III}; the upper index
  $j$ is preserved because $F_\theta$ preserves the number of `$h$'
  edges.  The diagram in part (b) commutes at the chain level, because
  if $q_k$ denotes the maximal lattice point in $P_\Lambda$ with
  respect to the ordering \eqref{eqn:ordering}, then
  $C_*^\theta(\Lambda')$ is the subcomplex of $C_*^\theta(\Lambda)$ in
  which $q_k$ does not appear in any of the $\theta$--corner sequences.
\end{proof}

\subsection{The special degree zero homology}
\label{sec:zero}

We now apply \fullref{prop:flattening} to compute
$H_0^{(-2n)}(\Lambda)$, where $\Lambda$ is an $n$--convex
path.  Generators for this homology are given explicitly as follows.

\begin{definition}
If $a,b\in P_\Lambda$ are distinct lattice points with $b-a$
indivisible, define $Z_n(a,b)\in C_*(\Lambda)$ to be the generator
that wraps $n$ times around the $2$--gon with vertices $a$ and $b$,
with all $2n$ edges labeled `$h$', ordered counterclockwise with an
edge from $a$ to $b$ coming first.  If $a,b\in P_\Lambda$ are distinct
lattice points with $m-1$ lattice points in the interior of the line
segment between them, define
\[
Z_n(a,b) \eqdef
\sum_{i=1}^{m}Z_n\left(a+\frac{i-1}{m}(b-a),a+\frac{i}{m}(b-a)\right).
\]
Define $Z_n(a,a)\eqdef 0$.  Observe that $\delta Z_n(a,b)=0$.
\end{definition}

\begin{definition}
\label{def:simpleTriangle}
A {\em simple triangle\/} is a triple $(a,b,c)$ of non-colinear points
in $\Z^2$, ordered counterclockwise, such that the triangle with
vertices $a,b,c$ encloses no other lattice points.
\end{definition}

\begin{lemma}
\label{lem:znRelations}
Let $\Lambda$ be an $n$--convex path.  Then for any lattice points
$a,b,c\in P_\Lambda$, the following relations hold in the homology
$H_*(\Lambda)$:
\begin{gather}
\label{eqn:znRelation1}
Z_n(a,b) + Z_n(b,a) = 0,\\
\label{eqn:znRelation2}
Z_n(a,b) + Z_n(b,c) + Z_n(c,a) = 0.
\end{gather}
\end{lemma}

\begin{proof}
  Equation \eqref{eqn:znRelation1} holds at the chain level
  because $Z_n(b,a)$ is obtained from $Z_n(a,b)$ by re-ordering the
  `$h$' edges, and a $2n$--cycle is an odd permutation.

Likewise, equation \eqref{eqn:znRelation2} holds at the chain level if
  $a,b,c$ are colinear.

Next we prove \eqref{eqn:znRelation2} when $(a,b,c)$ is a simple
  triangle.  Let $\lambda$ be the $n$--convex path that wraps $n$ times
  around this triangle.  Write the corners of $\lambda$ in
  counterclockwise order as $c_1,\ldots,c_{3n}$, starting with a
  corner that maps to $a$.  Define an ``admissible corner set'' to be
  a subset $I\subset\{1,\ldots,3n\}$ such that for distinct $i,j\in
  I$, the corners $c_i,c_j$ are not adjacent, ie $|i-j|\neq 1,
  3n-1$.  Let $\mc{C}$ denote the set of admissible corner sets.  If
  $I=\{i_1,\ldots,i_k\}\in\mc{C}$, let
\[
\lambda(I) \eqdef \lambda\setminus c_{i_1} \setminus
c_{i_2} \setminus \cdots \setminus c_{i_k}.
\]
Let $T(I)\in C_*(\Lambda)$ denote the generator with underlying path
$\lambda(I)$ and with all edges labeled `$h$'.  Order the `$h$' edges
counterclockwise, starting at $c_1$ if $1\notin I$ and starting at
$c_2$ if $1\in I$.  In this notation,
\begin{equation}
\label{eqn:znNotation}
\begin{aligned}
Z_n(a,b) &= T(\{3,6,\ldots,3n\}),\\
Z_n(b,c) &= T(\{1,4,\ldots,3n-2\}),\\
Z_n(c,a) &= -T(\{2,5,\ldots,3n-1\}).
\end{aligned}
\end{equation}
The differential of a generator $T(I)$ is given by
\begin{equation}
\label{eqn:DTI}
\delta T(I) = \sum_{I\cup\{i\}\in\mc{C}} (-1)^{\#\{j\notin I \mid
i<j\}}T(I\cup\{i\}).
\end{equation}
Now let $\mc{C}_0$ denote the set of admissible corner sets of the form
$I=\{i_1,\ldots,i_{n-1}\}$ with $i_1<i_2<\cdots < i_{n-1}$ and with
$i_1,i_3,\ldots$ odd and $i_2,i_4,\ldots$ even.  By
\eqref{eqn:znNotation} and \eqref{eqn:DTI},
\[
\delta\Bigl(\sum_{I\in\mc{C}_0}T(I)\Bigr)
 = Z_n(a,b) + Z_n(b,c) + Z_n(c,a).
\]
To see this, observe that a generator $T(I)$ with
$I=\{i_1<i_2<\cdots<i_n\}$ will appear exactly once in the left hand
side if the $i_j$'s alternate parity, twice with opposite signs if
there is exactly one $j$ such that $i_j$ and $i_{j+1}$ have the same
parity, and otherwise not at all.

For general noncolinear lattice points $a,b,c\in P_\Lambda$, we can
triangulate the triangle with vertices $a,b,c$ by simple triangles.
Adding the relations \eqref{eqn:znRelation2} for these simple
triangles, and using \eqref{eqn:znRelation1} to cancel interior edges,
proves \eqref{eqn:znRelation2} for $(a,b,c)$.
\end{proof}

The following is a useful way to understand the above relations.  If
$S$ is any set, let $\mc{I}(S)$ denote the set of finite formal sums
of elements of $S$ with integer coefficients such that the sum of the
coefficients is zero.  For $a,b\in S$, define $z(a,b)\eqdef a-b\in
\mc{I}(S)$.  We then have the following elementary fact, whose proof
is left to the reader:

\begin{lemma}
\label{lem:understandingRelations}
If $S$ is any set, then as a $\Z$--module, $\mc{I}(S)$ is generated by
$$\{z(a,b)\mid a,b\in S\},$$
with the relations $z(a,b)+z(b,a)=0$ and
$z(a,b)+z(b,c)+z(c,a)=0$.
\end{lemma}

\begin{proposition}
\label{prop:zeroHomology}
Let $\Lambda$ be an n-convex path.  Then there is an isomorphism
\[
\mc{I}(P_\Lambda\cap\Z^2)\simeq H_0^{(-2n)}(\Lambda)
\]
sending $z(a,b)\mapsto Z_n(a,b)$.
\end{proposition}

\begin{proof}
By Lemmas~\ref{lem:znRelations} and \ref{lem:understandingRelations},
the above map $\mc{I}(P_\Lambda\cap\Z^2)\to H_0^{(-2n)}$ is
well-defined.  To show that it is an isomorphism, let $k=k_\Lambda$
and let $\Lambda_0$ be the $n$--convex path on the $x$--axis with
$P_{\Lambda_0}=[1,k]\times\{0\}$. Choose any $\theta\in(0,\pi)$ with
$\tan\theta$ irrational. By \fullref{prop:flattening}(a),
there is an isomorphism on homology
\[
(F_\theta)_*\co  H_*^{(j)}(\Lambda_0) \stackrel{\simeq}{\longrightarrow}
H_*^{(j)}(\Lambda).
\]

Since every generator of $C_*(\Lambda_0)$ has nonnegative index
(cf equation \eqref{eqn:xIndex} below) and at most $2n$ edges
labeled `$h$', the generators of $C_*^{(-2n)}(\Lambda_0)$ are those
with $2n$ edges labeled `$h$' and with index $0$, which since all
edges are labeled `$h$' means that all edges have length one.  Thus
$C_*^{(-2n)}(\Lambda_0)=H_0^{(-2n)}(\Lambda_0)$ is freely generated by
the $k-1$ generators $Z_n((1,0),(2,0)), \ldots, Z_n((k-1,0),(k,0))$.

By the construction of the chain map $F_\theta$ in \fullref{sec:flattening},
\[
F_\theta(Z_n((i,0),(i+1,0))) = Z_n(q_i,q_{i+1})
\]
where $q_1,\ldots,q_k$ are the lattice points in $P_\Lambda$, ordered
by \eqref{eqn:ordering}.  Hence $H_0^{(-2n)}(\Lambda)$ is freely
generated by $Z_n(q_1,q_2), \ldots, Z_n(q_{k-1},q_k)$.  But it follows
from the definition of $\mc{I}$ that
$\mc{I}(P_\Lambda\cap\Z^2)$ is freely generated by
$z(q_1,q_2),\ldots,z(q_{k-1},q_k)$.
\end{proof}

\subsection{Vanishing of homology for $\Gamma\neq0$}
\label{sec:vanishing}

We now prove \fullref{thm:vanishing}.  Without loss of
generality,
$$\Gamma=\begin{pmatrix} k \\ 0 \end{pmatrix}$$
for some positive integer $k$.
This is justified by the following lemma.  An element $A\in SL_2\Z$
induces a diffeomorphism $S^1\to S^1$, where $S^1$ is the unit
circle in the $\R^2$ on which $A$ acts.  This diffeomorphism can be
lifted to a diffeomorphism ${f}\co \R\to\R$.  The set of pairs
$\left(A,{f}\right)$ forms a group $\widetilde{SL}_2\Z$
which is an infinite cyclic cover of $SL_2\Z$.

\begin{lemma}
\label{lem:symmetry}
A pair $\left(A,{f}\right) \in \widetilde{SL}_2\Z$ induces an
isomorphism
\[
\Phi_{\left(A,{f}\right)} \co  H_*(2\pi n;\Gamma)
\stackrel{\simeq}{\longrightarrow} H_*(2\pi n;A\Gamma).
\]
\end{lemma}

\begin{proof}
  If $\Lambda\co \R\setminus p^{-1}(T)\to\Z^2$ is an admissible path of
  rotation number $n$ and period $\Gamma$, define an admissible path
  $\Phi\Lambda$ of period $A\Gamma$ by
\[
\Phi\Lambda \eqdef A\circ \Lambda\circ {f}^{-1}.
\]
Then ${f}$ induces a bijection from the edges of $\Lambda$
to the edges of $\Phi\Lambda$, and pushing forward edge labels via
this bijection gives an isomorphism of $\Z$--modules
\[
\Phi_{\left(A,{f}\right)}\co
C_*(2\pi n;\Gamma) \longrightarrow C_*(2\pi n;A\Gamma).
\]
It follows from the definition of rounding that this is a chain map,
since the action of $SL_2\Z$ on $\Z^2$ preserves convex hulls and
(unlike the more general action of $GL_2\Z$) respects the signs in the
differential $\delta$.
\end{proof}

For $a,b\ge 1$, let $\Box(a,b;k)$ denote the closed admissible path of
rotation number $n$ and period $\Gamma$ whose restriction to $(0,2\pi
n]$ wraps $n$ times around the rectangle $[0,a-1]\times
[0,b-1]\subset\R^2$, except that the edge at $\theta=2\pi n$ has
length $a-1+k$.  For example, $\Box(1,1;k)$ has edges of length $k$
along the $x$--axis at angles $\theta=2\pi i n$ separated by kinks
parametrized by the intervals $(2\pi i n, 2\pi (i+1) n)$.  We now have
the following analogue of \fullref{prop:flattening}(a), which
replaces the path $\Box(a,b;k)$ by a path on the $x$--axis.

\begin{lemma}
\label{lem:GNEFlattening}
There is a chain map
\[
F\co  C_*(\Box(ab,1;kb)) \longrightarrow C_*(\Box(a,b;k))
\]
which induces an isomorphism on homology, preserves the relative
grading, and sends
\[
E_{\Box(1,1;kb)} \longmapsto E_{\Box(1,1;k)}.
\]
\end{lemma}

\begin{proof}
  If $\Lambda$ is any periodic admissible path of rotation number $n$
  and period $\Gamma$, and if $\theta\in(0,2\pi)$ is an angle with
  $\tan\theta$ irrational, then we can define a subcomplex
  $C_*^\theta(\Lambda)\subset C_*(\Lambda)$ as in
  \fullref{sec:flattening}.  The only difference is that now a
  $\theta$--corner sequence is an infinite sequence $\{p_i\mid
  i\in\Z\}$ of points in $\Z^2$ such that $p_{i+2n}=p_i+\Gamma$ for
  all $i$.  As in \fullref{prop:III}, the inclusion
  $C_*^\theta(\Lambda)\to C_*(\Lambda)$ induces an isomorphism on
  homology.
  
  Now let $\theta\eqdef \frac{\pi}{2} + \epsilon$ where $\epsilon>0$
  is small with respect to $b$.  Then for two consecutive points in a
  $\theta$--corner sequence for $\Box(a,b;k)$, the ordering
  \eqref{eqn:ordering} coincides with the lexicographic order on
  $\Z^2$.  There is then an isomorphism of chain complexes
\[
C_*^\theta(\Box(a,b;k)) \stackrel{\simeq}{\longrightarrow}
C_*^\theta(\Box(ab,1;kb)) = C_*(\Box(ab,1;kb)).
\]
This is defined via the bijection on $\theta$--corner sequences induced by the
map
\[
\begin{split}
\Z\times\{0,1,\ldots,b-1\} &\longrightarrow \Z\times \{0\},\\
(x,y) & \longmapsto (bx+y,0).
\end{split}
\]
This is a chain map and preserves the relative grading as in the proof
of \fullref{lem:CICC}, because for $\Lambda=\Box(a,b;k)$, locally
$\Lambda^{i-1}_{p_{i-1},p_i}$ and $\Lambda^i_{p_i,p_{i+1}}$ are the
same as they would be if $\Lambda$ were $n$--convex and $P_\Lambda$
were a rectangle.  Clearly this isomorphism sends $E_{\Box(1,1;k)}$
to $E_{\Box(1,1;kb)}$.
\end{proof}

The homology $H_*(\Box(a,b;k))$ may be complicated, but we will only
need to establish a lower bound on the index of nonvanishing homology
groups.  For $y_0\in\Z$, let $C_*(2\pi n;\Gamma;y\ge y_0)$ denote the
subcomplex of $C_*(2\pi n;\Gamma)$ spanned by generators whose
underlying admissible paths map to the half-plane $y\ge y_0$.

\begin{lemma}
\label{lem:indexBound}
If $\alpha\in H_*(2\pi n;\Gamma;y\ge 0)$ is nonzero and has pure
degree, then the relative index
\[
I\left(\alpha,E_{\Box(1,1;k)}\right) \ge -n.
\]
\end{lemma}

\begin{proof}
By horizontally translating a cycle representing $\alpha$ (which does
not affect the index by \eqref{eqn:indexAmbiguity}), we may assume
that $\alpha$ is contained in $H_*(\Box(a,b;k))$ for some $a,b$.  By
\fullref{lem:GNEFlattening}, $\alpha=F(\alpha_0)$ for some
$\alpha_0\in H_*(\Box(ab;1;kb))$, and
\[
I\left(\alpha,E_{\Box(1,1;k)}\right) =
I\left(\alpha_0,E_{\Box(1,1;kb)}\right).
\]
Let $\beta_0$ be a generator in a cycle representing $\alpha_0$.
Since $\beta_0$ and $E_{\Box(1,1;kb)}$ are on the $x$--axis, the
definition of the relative index \eqref{eqn:relativeIndex}
implies that
\[
I\left(\beta_0,E_{\Box(1,1;kb)}\right) = \ell(\beta_0) -
\#h(\beta_0) - kb.
\]
Since $\beta_0$ has period
$$\begin{pmatrix} kb \\ 0 \end{pmatrix},$$
we see that $\ell(\beta_0) - \#h(\beta_0) \ge kb-n$.
\end{proof}

\begin{proof}[Proof of \fullref{thm:vanishing}]  Suppose $\alpha\in
H_*(2\pi n;\Gamma)$ is nonzero and has pure degree.  Since
\[
C_*(2\pi n;\Gamma) = \lim_{y_0\to -\infty} C_*(2\pi n;\Gamma;y\ge
y_0),
\]
there exists $y_0$ such that
\[
\alpha\in H_*(2\pi n;\Gamma;y\ge y_0-m)
\]
for all $m\ge 0$. Let $\Psi\co C_*(2\pi n;\Gamma) \to C_*(2\pi n;\Gamma)$ be
the isomorphism of chain complexes that translates paths upward by one
unit, ie
\[
(\Psi\Lambda)(t) \eqdef \Lambda(t) + \begin{pmatrix} 0 \\ 1
\end{pmatrix}.
\]
By \fullref{lem:indexBound} and symmetry, for all $m\ge 0$,
\begin{equation}
\label{eqn:anotherIndexBound}
I\left(\alpha,\Psi^{y_0-m}E_{\Box(1,1;k)}\right) \ge -n.
\end{equation}
By the index ambiguity formula \eqref{eqn:indexAmbiguity},
\begin{align*}
I\left(\alpha,\Psi^{y_0-m}E_{\Box(1,1;k)}\right)
&{=} I\left(\alpha,\Psi^{y_0}E_{\Box(1,1;k)}\right){+}
I\left(\Psi^{y_0}E_{\Box(1,1;k)},\Psi^{y_0-m}E_{\Box(1,1;k)}\right)\\
&{=} I\left(\alpha,\Psi^{y_0}E_{\Box(1,1;k)}\right){-}2km.
\end{align*}
Combining this with \eqref{eqn:anotherIndexBound} gives a
contradiction when $m$ is sufficiently large.
\end{proof}

\section{Computation of $\wtilde{H}_*(2\pi n;0)$}
\label{sec:HTilden}

We now compute the homology for closed admissible paths of rotation
number $n$.  Recall the notations $\mc{I}(\Z^2)$ and $\Zmodule$ from
\fullref{sec:variantsIntro}.  Also recall the notation
$\wwtilde{H}_*^{(j)}(2\pi n;0)$ from
\fullref{definition:subcomplex}.  In this section we will
prove:

\begin{theorem}
\label{thm:HTilden}
\begin{enumerate}
\item[(a)]
As $\Z[\Z^2]$--modules,
\[
\wwtilde{H}_i^{(j)}(2\pi n;0) \simeq \left\{\begin{array}{cl}
\Zmodule, & i=2k,2k+1;
\;j=2k-2n+2; k\in\Z_{\ge 0},\\
\mc{I}(\Z^2), &
i=0,\; j=-2n,\\
0, & \mbox{otherwise}.
\end{array}\right.
\]
\item[(b)]
If $i\ge 2$, then the map
\[
U\co \wwtilde{H}_i(2\pi n;0) \longrightarrow
\wwtilde{H}_{i-2}(2\pi n;0)
\]
is an isomorphism between the $\Zmodule$ components.
\end{enumerate}
\end{theorem}

In particular, $\wwtilde{H}_*(2\pi n;0)$ is independent of $n$.
However, the homology generators look very different for different
$n$.  Among other things, increasing $n$ by $1$ increases the number
of `$h$' edges (in the chain complex generators that are summands in
cycles generating the homology) by $2$.  For example, if $\Lambda$ is
a closed admissible path of rotation number $n$, then the cycles
$E_\Lambda, H_\Lambda\in \wwtilde{C}_*(2\pi n;0)$ are homology
generators when $n=1$, but they are boundaries when $n>1$.

\subsection{Splicing}
\label{sec:splicing}

\begin{definition}
For $n\ge 1$, let $CX_*(n)$ denote the subcomplex of $C_*(2\pi n;0)$
spanned by generators with underlying admissible path on the $x$--axis
in $\R^2$.  Let
\[
CX_*^{(j)}(n) \eqdef CX_*(n)\cap C_*^{(j)}(2\pi n;0).
\]
\end{definition}

Thanks to the ``flattening'' technology of \fullref{sec:flattening},
to prove \fullref{thm:HTilden} it is mostly sufficient to compute the
homology $HX_*(n)$.  We will now compute $HX_*(n)$ by induction on
$n$, using a ``splicing'' isomorphism constructed below.

We will use the following notation for generators of $CX_*(n)$.
Similarly to \fullref{sec:flattening}, define a {\em corner sequence\/}
to be a sequence of integers $a_0,a_1,\ldots,a_{2n}=a_0$ with $a_0\ge
a_1\le a_2\ge\cdots$, and let $\mc{S}_n$ denote the set of all corner
sequences.  If $\alpha$ is a generator of $CX_*(n)$, then the
$x$--coordinates of the values of its underlying admissible path at
$\pi/2+i\pi$ determine a corner sequence $\{a_i\}$.  Thus as a
$\Z$--module,
\begin{equation}
\label{eqn:CX_*}
CX_*(n) = \bigoplus_{\{a_i\}\in\mc{S}_n} \bigotimes_{i=0}^{2n-1}
\Z\left\{e_{a_i}^{a_{i+1}},h_{a_i}^{a_{i+1}}\right\}.
\end{equation}
Here $e_a^b$ and $h_a^b$ denote edges from $(a,0)$ to $(b,0)$ labeled
`$e$' and `$h$' respectively.  If $a=b$, then $e_a^b$ denotes the lack
of an edge, and we interpret $h_a^b=0$.  (So if $a = b$, then $\Z\{e_a^b, h_a^b\}$ really means $\Z\{e_a^b\}$.)  This decomposition of
$CX_*(n)$ determines a natural convention for ordering the `$h$'
edges.  By the index formula \eqref{eqn:closedIndex}, the index of a
generator is the sum of the lengths of the edges, minus the number of
`$h$' edges.  That is, if $\alpha$ has corner sequence $\{a_i\}$, then
\begin{equation}
\label{eqn:xIndex}
I(\alpha) = \sum_{i=0}^{2n-1}\left| a_i - a_{i+1}\right| - \#h(\alpha).
\end{equation}
The differential on $CX_*(n)$ is given as in equation
\eqref{eqn:d1}.

\begin{definition}
\label{def:splicing}
Define the ``splicing map''
\[
S\co CX_*^{(j)}(n)\longrightarrow CX_*^{(j-2)}(n+1)
\]
as follows.  In terms of the decomposition \eqref{eqn:CX_*}, any
generator of $CX_*(n)$ has the form $w\tensor e_a^b\tensor e_b^c$,
$w\tensor e_a^b\tensor h_b^c$, $w\tensor h_a^b\tensor e_b^c$, or
$w\tensor h_a^b\tensor h_b^c$, where $w$ is a tensor product of $2n-2$
$e$'s or $h$'s, and $a\ge b \le c$.  Henceforth we will omit the
tensor product symbol when writing generators of $CX_*(n)$ this way.
Now $S$ sums over all ways of shortening the edges from $a$ to $b$ and
from $b$ to $c$ and splicing in two `$h$' edges between them so that
the grading is preserved.  Namely
\[
\begin{split}
S\left(we_a^be_b^c\right) & \eqdef \sum_{i=b}^{a}\sum_{j=b}^{c} w e_a^i
h_i^{i+j-b+1} h_{i+j-b+1}^j e_j^c,\\
S\left(wh_a^be_b^c\right) & \eqdef \sum_{i=b}^{a-1}\sum_{j=b}^{c} w h_a^i
h_i^{i+j-b+1} h_{i+j-b+1}^j e_j^c,
\end{split}
\]
and similarly for $S\left(we_a^bh_b^c\right)$ and $S\left(wh_a^bh_b^c\right)$.
\end{definition}

\begin{lemma}
$S$ is a degree $0$ chain map.
\end{lemma}

\begin{proof}
$S$ preserves the grading because
\[
(a-b) + (c-b)
=
(a-i) + ([i+j-b+1]-i) + ([i+j-b+1]-j) + (c-j) - 2.
\]
We now check that $\delta S = S \delta$.  Recall the notation
$c_\theta$ from \fullref{sec:U}.  Decompose the differential on
$CX_*{(n)}$ as $\delta = \delta_0+\delta_1$ where $\delta_1$ is the
contribution from rounding at the corner $c_{(2n-1/2)\pi}$.  Similarly
decompose the differential on $CX_*{(n+1)}$ as
$\delta=\delta_0+\delta_1$ where $\delta_1$ is the contribution from
rounding at the corners $c_{(2n-1/2)\pi}$, $c_{(2n+1/2)\pi}$, and
$c_{(2n+3/2)\pi}$.  Clearly $\delta_0S=S\delta_0$; the signs agree
because the number of `$h$' edges we are splicing in is always
even. So we just have to check that $\delta_1S=S\delta_1$.  We will
verify that $\delta_1S\alpha=S\delta_1\alpha$ when $\alpha =
we_a^bh_b^c$; the calculations for the other three types of generators
are very similar.  By the definitions of $S$ and $\delta_1$,
\[
\begin{split}
  \delta_1S\bigl(we_a^bh_b^c\bigr) & = \delta_1
  \sum_{i=b}^{a}\sum_{j=b}^{c-1} w
  e_a^i h_i^{i+j-b+1} h_{i+j-b+1}^j h_j^c\\ 
  &= w\sum_{j=b}^{c-1}\biggl(\sum_{i=b}^{a-1} e_a^{i+1}
    e_{i+1}^{i+j-b+1} h_{i+j-b+1}^j - \sum_{i=b+1}^a e_a^i e_i^{i+j-b}
    h_{i+j-b}^j\biggr)h_j^c\\
  &  + w\sum_{i=b}^a e_a^i \biggl(-\!\!\!\!\sum_{j=b+1}^{c-1}h_i^{i+j-b}
    e_{i+j-b}^j h_j^c + \sum_{j=b}^{c-2} h_i^{i+j-b+1}
    e_{i+j-b+1}^{j+1} h_{j+1}^c\biggr) 
\end{split}
\]
$\phantom{99}$\vspace{-15pt}
\[
\begin{split}
  &\phantom{=}  + w\sum_{i=b+1}^a\sum_{j=b}^{c-1} e_a^i h_i^{i+j-b+1}
  h_{i+j-b+1}^{j+1}
  e_{j+1}^c\hspace{1in}\\
  &=0+0+ w\sum_{i=b+1}^a\sum_{j=b+1}^c e_a^i h_i^{i+j-b}
  h_{i+j-b}^{j} e_{j}^c
  \\
  &= S\bigl(we_a^{b+1}e_{b+1}^c\bigr) = S\delta_1\bigl(we_a^bh_b^c\bigr).
\end{split}
\]
This completes the proof.
\end{proof}

We now come to the key argument which explains how the homologies for
different $n$ are related.

\begin{proposition}
\label{prop:xaxis}
For all $n\ge 1$, the chain map $S$ induces an isomorphism
\[
HX_*^{(j)}(n) \stackrel{\simeq}{\longrightarrow} HX_*^{(j-2)}(n+1).
\]
\end{proposition}

\begin{proof}
If $\alpha$ is a generator of $CX_*(n+1)$ corresponding to the
corner sequence $a_0,a_1,\ldots,a_{2n+2}=a_0$, define the
``degree''
\[
\deg(\alpha) \eqdef \sum_{i=0}^{2n-2}\left|a_i-a_{i+1}\right| +
\left|a_{2n+1}-a_{2n+2}\right|.
\]
In other words the degree is the sum of the lengths of the edges,
except for the two edges that could arise from splicing a generator of
$CX_*(n)$.  This ``degree'' defines an increasing filtration on
$CX_*(n+1)$ and hence gives rise to a spectral sequence $E^*_{*,*}$.
By equation \eqref{eqn:xIndex}, the filtration is bounded in terms of
the index by
\[
0 \le \deg(\alpha) \le I(\alpha) + 2n.
\]
It follows that the spectral sequence $E^*_{*,*}$ converges to
$HX_*(n+1)$.

For $i,j\in\Z$, let $C_*(i,j)$ denote the direct limit of
$C_*(\Lambda)$, where $\Lambda$ is an open admissible path on the
$x$--axis parametrized by the interval $(\pi/2+2\pi(n-1), \pi/2+2\pi
n)$ and with endpoints $(i,0)$ and $(j,0)$.  That is, generators of
$C_*(i,j)$ have underlying admissible paths that start at $(i,0)$, go
in the negative $x$ direction by some (possibly zero) amount, then
turn by $\pi$ and go in the positive $x$ direction by some (possibly
zero) amount to $(j,0)$.  In this notation, the $E^1$ term of the
spectral sequence is given by
\begin{multline}
\label{eqn:E1a}
E^1 = \bigoplus_{a_0,\ldots,a_{2n-1},a_{2n+1}} \bigotimes_{i=0}^{2n-2}
\Z\bigl\{ e_{a_i}^{a_{i+1}}, h_{a_i}^{a_{i+1}} \bigr\} \\[-10pt]
\tensor H_*(a_{2n-1},a_{2n+1}) \tensor
\Z\left\{e_{a_{2n+1}}^{a_{2n+2}},h_{a_{2n+1}}^{a_{2n+2}}\right\}
\end{multline}
where $a_i\ge a_{i+1}$ for $i\neq 2n$ even, $a_i\le a_{i+1}$ for $i\neq
2n-1$ odd, and $a_{2n+2}=a_0$.
By Propositions~\ref{prop:step1} and \ref{prop:step2},
\begin{equation}
\label{eqn:E1b}
H_*(i,j) = \left\{\begin{array}{cl} 0, & i\neq j,\\
\Z\left\{h_i^{i-1} h_{i-1}^i \right\}, & i=j.
\end{array}\right.
\end{equation}
By equations \eqref{eqn:CX_*}, \eqref{eqn:E1a} and \eqref{eqn:E1b},
there is an isomorphism of $\Z$--modules
\begin{equation}
\label{eqn:E1c}
CX_*(n) \simeq E^1.
\end{equation}
Now we can filter $CX_*(n)$ by the sum of the lengths of all the
edges.  This filtration gives rise to a spectral sequence $'E^*_{*,*}$
converging (at the third term) to $HX_*(n)$ with
\begin{equation}
\label{eqn:E1d}
'E^1 = CX_*(n).
\end{equation}
The splicing chain map $S\co CX_*(n)\to CX_*(n+1)$ respects the above
filtrations and therefore induces a morphism of spectral sequences
$'E^*_{*,*} \to E^*_{*,*}$.  On the first term, $S$ induces the
isomorphism $'E^1\stackrel{\simeq}{\longrightarrow} E^1$ given by
\eqref{eqn:E1c} and \eqref{eqn:E1d}, because by the definition of $S$,
the only term in which no old edge is shortened is the term in which
two `$h$' edges of length one are spliced in.  Therefore $S$ induces
an isomorphism $'E^\infty\stackrel{\simeq}{\longrightarrow} E^\infty$,
hence an isomorphism $HX_*(n) \stackrel{\simeq}{\longrightarrow}
HX_*(n+1)$.
\end{proof}

\begin{corollary}
\label{cor:HX}
For $(i,j)\neq (0,-2n)$,
\[
HX_i^{(j)}(n) 
\simeq
\left\{\begin{array}{cl} \Z, & i=2k,2k+1;
\;j=2k-2n+2; k\in\Z_{\ge 0},\\
0, & \mbox{otherwise}.
\end{array}\right.
\]
\end{corollary}

\begin{proof}
Applying \fullref{prop:step3} to $x$--axis polygons and taking
the direct limit proves the claim for $n=1$.  It follows by
\fullref{prop:xaxis} that the claim holds for all $n$.
\end{proof}

We will also need the following lemma in \fullref{sec:HTildenProof}:

\begin{lemma}
\label{lem:SU=US}
For appropriate choices of the angles used to define $U$,
\begin{equation}
\label{eqn:SU=US}
SU=US\co CX_i^{(j)}(n) \longrightarrow CX_{i-2}^{(j-4)}(n+1).
\end{equation}
\end{lemma}

\begin{proof}
   Define $U$ on $C_*(2\pi n;0)$ using $\theta=\pi/2+(2n-2)\pi$, and
  define $U$ on $C_*(2\pi(n+1);0)$ using $\theta=\pi/2+(2n-1)\pi$.  We
  will check that equation \eqref{eqn:SU=US} holds when applied to a
  generator of the form $we_a^bh_b^c$; the other cases are very
  similar. Note that $U$ acts on a generator of $CX_*(n)$ simply by
  shrinking the two edges adjacent to $c_\theta$ by one, preserving
  the edge labels.  Thus
\[
\begin{split}
US\bigl(we_a^bh_b^c\bigr) &= U\biggl(\sum_{i=b}^{a} \sum_{j=b}^{c-1}
we_a^i h_i^{i+j-b+1} h_{i+j-b+1}^j h_j^c\biggr) \\
&= \sum_{i=b+1}^{a} \sum_{j=b+1}^{c-1} we_a^i h_i^{i+j-b} h_{i+j-b}^j
h_j^c\\
&= S\bigl(we_a^{b+1}h_{b+1}^c\bigr) \\
&= SU\bigl(we_a^{b}h_{b}^c\bigr),
\end{split}
\]
which completes the proof.
\end{proof}

\subsection{Inclusion}
\label{sec:DZHI}

The proof of \fullref{thm:HTilden} will proceed by taking a direct
 limit, for which purpose we will need the following technical lemma
 about how $H_*(\Lambda)$ behaves under inclusion of polygonal regions.

\begin{lemma}
\label{lem:injectivity}
Let $\Lambda$ and $\Lambda'$ be $n$--convex paths with
$\Lambda'\le\Lambda$.  If $(i,j)\neq (0,-2n)$, and if $k_{\Lambda'}$ is
sufficiently large with respect to $i$ and $j$, then inclusion induces an
isomorphism
\[
H_i^{(j)}(\Lambda') \stackrel{\simeq}{\longrightarrow} H_i^{(j)}(\Lambda).
\]
\end{lemma}

To prove \fullref{lem:injectivity}, we will need to compute
$H_*(\Lambda)$ for certain special open admissible paths $\Lambda$ on
the $x$--axis.  Similarly to the notation in \fullref{sec:splicing}, we
can write
\[
\Lambda = \lambda_{a_0}^{a_1}\cdots \lambda_{a_{k-1}}^{a_k}
\]
where $\lambda_{a_i}^{a_{i+1}}$ denotes an edge from $(a_i,0)$ to
$(a_{i+1},0)$ occuring at angle $(i+i_0)\pi$ for some integer $i_0$.
Here we are not requiring that $a_0=a_k$.

\begin{definition}
  A {\em spiral\/} is a nonconstant open admissible path on the $x$--axis,
  $\Lambda=\lambda_{a_0}^{a_1}\cdots\lambda_{a_{k-1}}^{a_k}$, such that
  $k\ge 1$, and if $k>1$ then
\begin{equation}
\label{eqn:spiral}
\begin{split}
|a_0-a_1| & > |a_1-a_2|,\\
|a_i-a_{i+1}| &\ge |a_{i+1} - a_{i+2}|,
\quad i=1,\ldots,k-2.
\end{split}
\end{equation}
\end{definition}

\begin{lemma}
\label{lem:spiral}
If $\Lambda$ is a spiral, then
$H_*(\Lambda)=\Z\{E_{\Lambda},H_{\Lambda}\}$.
\end{lemma}

\begin{proof}
  If $\Lambda$ has only edge then the lemma is immediate.  Otherwise
  let $c$ be the corner of $\Lambda$ preceding the last edge.  Then
  $c$ splits $\Lambda=\Lambda_1\Lambda_2$ where $\Lambda_1$,
  $\Lambda_2$, and $\Lambda\setminus c$ are all spirals.  The lemma
  follows by induction using the long exact sequence
  \eqref{eqn:openExactSequence}, as in the proof of
  \fullref{prop:step1}.
\end{proof}


\begin{definition}
  A {\em semi-spiral\/} is an open admissible path on the $x$--axis,
  $\Lambda = \lambda_{a_0}^{a_1} \cdots \lambda_{a_{k-1}}^{a_k}
  \lambda_{a_k}^{a_{k+1}}$, with $k\ge 2$, satisfying the conditions
  \eqref{eqn:spiral} and
\[
|a_0 - a_1| \le |a_k - a_{k+1}|.
\]

If $\Lambda$ is a semi-spiral, define $V_\Lambda, W_\Lambda \in
C_*(\Lambda)$ by
\[
\begin{split}
  V_{\Lambda} &\eqdef E_{\lambda_{a_0}^{a_1}\cdots\lambda_{a_{k-2}}^{a_{k-1}}}
  e_{a_{k-1}}^{a_{k-1}} h_{a_{k-1}}^{a_{k+1}},\\
  W_{\Lambda} &\eqdef H_{\lambda_{a_0}^{a_1}\cdots\lambda_{a_{k-2}}^{a_{k-1}}}
  e_{a_{k-1}}^{a_{k-1}} h_{a_{k-1}}^{a_{k+1}}.
\end{split}
\]
Note that the generator $V_\Lambda$, and the generators in the sum
$W_\Lambda$, each have one kink.
\end{definition}

\begin{lemma}
\label{lem:EHXZ}
If $\Lambda$ is a semi-spiral,
then $H_*(\Lambda)=\Z\{E_{\Lambda},H_{\Lambda},V_\Lambda,W_\Lambda\}$.
\end{lemma}

\begin{proof}
  We use induction on the length of $\Lambda$.  If $\Lambda$
  has only two edges, so that $a_1=a_2=\cdots =a_k$, then the lemma
  follows by inspection.  Namely, $C_*(\Lambda)$ has only four generators $ee$,
  $eh$, $he$, $hh$ which correspond to $E_\Lambda$, $V_\Lambda$,
  $H_\Lambda-V_\Lambda$, and $W_\Lambda$ above.  The differential
  vanishes here because the corner between the two edges is a kink since
  $k\ge 2$.

If $\Lambda$ has more than two edges, let $c$ denote the corner
preceding the second-to-last edge.  Then $c$ splits
$\Lambda=\Lambda_1\Lambda_2$ where $\Lambda_2$ and $\Lambda\setminus
c$ are semi-spirals, while $\Lambda_1$ is a spiral.  We now use the
long exact sequence \eqref{eqn:openExactSequence}.  The connecting
homomorphism
\[
\partial_c\co  H_*(C_*(\Lambda_1)\tensor C_*(\Lambda_2)) \longrightarrow
H_*(\Lambda\setminus c)
\]
is given by the equations \eqref{eqn:connecting} together with
\[
\begin{split}
E_{\Lambda_1}\tensor V_{\Lambda_2} & \longmapsto 0,\\
H_{\Lambda_1}\tensor V_{\Lambda_2},\, -E_{\Lambda_1}\tensor
W_{\Lambda_2} & \longmapsto V_{\Lambda\setminus c},\\
H_{\Lambda_1}\tensor W_{\Lambda_2} &\longmapsto -W_{\Lambda\setminus c}.
\end{split}
\]
So $H_*(\Lambda)$ is freely generated by $E_{\Lambda_1}E_{\Lambda_2} =
E_\Lambda$, $E_{\Lambda_1}H_{\Lambda_2} + H_{\Lambda_1}E_{\Lambda_2} =
H_\Lambda$, $E_{\Lambda_1}V_{\Lambda_2} = V_\Lambda$, and
$H_{\Lambda_1}V_{\Lambda_2} + E_{\Lambda_1}W_{\Lambda_2} = W_\Lambda$.
\end{proof}

\begin{proof}[Proof of \fullref{lem:injectivity}]  By
\fullref{prop:roundingSequence}, $P_{\Lambda'}$ is obtained
from $P_\Lambda$ by a finite sequence of corner roundings.  We can
assume that $P_{\Lambda'}$ is obtained from $P_\Lambda$ by rounding a
single corner $c$ (the general case then follows by induction).
Choose $\theta\in\R/2\pi\Z$ such that $c$ is the distinguished corner
$c_\theta$ of $P_\Lambda$.  By \fullref{prop:flattening}, we
can assume that $\Lambda=\Lambda_0$ and $\Lambda'=\Lambda_0'$ are on
the $x$--axis.  Without loss of generality,
$P_{\Lambda_0}=[1,k]\times\{0\}$ and
$P_{\Lambda_0'}=[1,k-1]\times\{0\}$.

Let $c_1,\ldots,c_n$ denote the corners of $\Lambda_0$ at $(k,0)$, in
counterclockwise order.  For $p,q\in\{1,\ldots,n\}$ distinct, let
$\Lambda_{p,q}$ denote the open path given by the portion of
$\Lambda_0$ starting at $c_p$ and ending at $c_q$, with the
intermediate corners $c_{p+1},\ldots,c_{q-1}$ rounded.  There is then
a decomposition of $\Z$--modules
\begin{multline*}
C_*(\Lambda_0)=\\
C_*(\Lambda_0') \oplus \bigoplus_{m=1}^n
\bigoplus_{1\le p_1<\cdots<p_m\le n}
\!\!\!\!C_*(\Lambda_{p_1,p_2})\tensor
\cdots\tensor C_*(\Lambda_{p_{m-1},p_m})
\tensor C_*(\Lambda_{p_m,p_1}).
\end{multline*}
Furthermore, $m$ defines an increasing filtration on $C_*(\Lambda_0)$,
where we interpret the $C_*(\Lambda_0')$ summand as corresponding to
$m=0$.  (One can regard $m$ as a weighted count of the times that an
admissible path $\lambda\le\Lambda_0$ stops at the point $(k,0)$.)
Thus we obtain a spectral sequence $E^*_{*,*}$ converging to
$H_*(\Lambda_0)$.  The $E^1$ term is given by
\begin{multline}
\label{eqn:anotherE1Term}
E^1 = H_*(\Lambda_0') \oplus \bigoplus_{m=1}^n
\bigoplus_{1\le p_1<\cdots<p_m\le n}
H_*(\Lambda_{p_1,p_2})\tensor\cdots\\[-10pt]
\cdots\tensor H_*(\Lambda_{p_{m-1},p_m})
\tensor H_*(\Lambda_{p_m,p_1}).
\end{multline}
Here we have used the fact that $H_*(\Lambda_{p,q})$ has no torsion,
which is justified below.

The key now is to compute the indices of the generators of
$H_*(\Lambda_{p,q})$.  Let $n'=q-p$ if $p<q$, and let $n'=p+n-q$ if
$p\ge q$.  If $n'=1$, then $H_*(\Lambda_{p,q})$ is given by
\fullref{prop:step2}.  The indices of the generators are
\[
I(Z_1((k,0),(k-1,0)))=0,\quad I(E_{\Lambda_{p,q}})=2k-2, \quad
I(H_{\Lambda_{p,q}})=2k-3.
\]
If $n'>1$, then $\Lambda_{p,q}$ is a semi-spiral, so
$H_*(\Lambda_{p,q})$ is given by \fullref{lem:EHXZ}.  (Note that $k$
here is different from the $k$ in the definition of semi-spiral.)
The indices of the generators of $H_*(\Lambda_{p,q})$ are
\[
\begin{split}
I(E_{\Lambda_{p,q}}) &= 2n'(k-2)+2,\\
I(H_{\Lambda_{p,q}}) &= 2n'(k-2)+1,\\
I(V_{\Lambda_{p,q}}) &= 2(n'-1)(k-2)+1,\\
I(W_{\Lambda_{p,q}}) &= 2(n'-1)(k-2).
\end{split}
\]
In particular, the index of each generator, except for
$Z_1((k,0),(k-1,0))$, is a linearly increasing function of $k$.
Therefore the index of every $m>0$ homology class in the $E^1$ term
\eqref{eqn:anotherE1Term} is a linearly increasing function of $k$,
except for the product of $n$ $Z_1$'s in the $m=n$ piece, which lives
in $C_0^{(-2n)}(\Lambda_0)$.  It follows that if $k$ is sufficiently
large with respect to $(i,j)\neq (0,-2n)$, then the map
\[
H_i^{(j)}(\Lambda_0') \longrightarrow H_i^{(j)}(\Lambda_0)
\]
induced by the inclusion is an isomorphism.
\end{proof}

\subsection{The direct limit}
\label{sec:HTildenProof}

We can now complete the calculation of $\wwtilde{H}_*(2\pi n;0)$.

\begin{lemma}
\label{lem:xaxis}
If $(i,j)\neq (0,-2n)$, then the inclusion $CX_*^{(j)}(n) \to
C_*^{(j)}(2\pi n;0)$ induces an isomorphism
\[
HX_i^{(j)}(n) \stackrel{\simeq}{\longrightarrow} H_i^{(j)}(2\pi n;0).
\]
\end{lemma}

\begin{proof}
Let $\Lambda_0$ be an $n$--convex path on the $x$--axis, and consider
the commutative diagram
\begin{equation}
\label{eqn:CD}
\xymatrix{
{H_i^{(j)}(\Lambda_0)} \ar[r] \ar[dr] & {HX_i^{(j)}(n)} \ar[d] 
\\ & {H_i^{(j)}(2\pi n;0)} }
\end{equation}
where the arrows are induced by the inclusions of chain complexes.
The homology $H_i^{(j)}(2\pi n;0)$ is given by the direct limit
\begin{equation}
\label{eqn:subcomplexDirectLimit}
H_i^{(j)}(2\pi n;0) = \lim_{\longrightarrow} H_i^{(j)}(\Lambda)
\end{equation}
over $n$--convex paths $\Lambda$.  It follows by
\fullref{lem:injectivity} that if $k_{\Lambda_0}$ is sufficiently large
with respect to $(i,j)$, then the diagonal arrow in the diagram
\eqref{eqn:CD} is an isomorphism, and likewise the horizontal arrow is
an isomorphism.  Hence the vertical arrow is an isomorphism.
\end{proof}

\begin{proof}[Proof of \fullref{thm:HTilden}]  (a) To start, we obtain an
isomorphism of $\Z$--modules
\[
H_0^{(-2n)}(2\pi n;0) \simeq \mc{I}(\Z^2)
\]
by taking the direct limit
\eqref{eqn:subcomplexDirectLimit} and applying
\fullref{prop:zeroHomology}.  This is in fact an isomorphism
of $\Z[\Z^2]$--modules, because by the definition of $Z_n(a,b)$,
translation by $w\in\Z^2$ sends $Z_n(a,b)$ to $Z_n(a+w,b+w)$.

For the rest of the proof suppose that $(i,j)\neq (0,-2n)$.  By
\fullref{cor:HX} and \fullref{lem:xaxis}, we have
\[
H_i^{(j)}(2\pi n;0) \simeq \left\{\begin{array}{cl} \Z, & i=2k,2k+1\ge 0,
\;j=2k-2n+2,\\
0, & \mbox{otherwise}.
\end{array}\right.
\]
To complete the proof of part (a), we must show that translations act
by the identity on $H_i^{(j)}(2\pi n;0)$.  Taking the direct limit of
\fullref{prop:step3} applied to $x$--axis polygons shows that
$HX_i^{(j)}(1)$, if nonzero, is generated by $E_{\Lambda_0}$ or
$H_{\Lambda_0}$ where $\Lambda_0$ is on the $x$--axis. By
\fullref{lem:convenient}, the homology class of $E_{\Lambda_0}$ or
$H_{\Lambda_0}$ in $HX_i^{(j)}(1)$ depends only on the number of
lattice points enclosed by $\Lambda_0$.  Therefore translation in the
$x$ direction acts by the identity on $HX_i^{(j)}(1)$.  By
\fullref{def:splicing}, the splicing isomorphism $S$ commutes
with translation in the $x$ direction.  By
\fullref{prop:xaxis}, it follows that translation in the $x$
direction acts by the identity on $HX_i^{(j)}(n)$ for all $n$.  So by
\fullref{lem:xaxis}, translation in the $x$ direction acts by the
identity on $H_i^{(j)}(2\pi n;0)$.  By the symmetry of
\fullref{lem:symmetry}, all translations act by the identity on
$H_i^{(j)}(2\pi n;0)$.

(b) By \fullref{lem:xaxis} and part (a), it is enough to show that
if $i\ge 2$ and $j\neq 2-2n$ then $U$ induces an isomorphism
\begin{equation}
\label{eqn:UIsom}
U\co HX_i^{(j)}(n) \stackrel{\simeq}{\longrightarrow}
HX_{i-2}^{(j-2)}(n).
\end{equation}
When $n=1$, we know that \eqref{eqn:UIsom} is an isomorphism by the
 above description of the generators of $HX_i^{(j)}(1)$ and
 \fullref{prop:UProperties}(c).  It then follows from
 \fullref{prop:xaxis} and \fullref{lem:SU=US} that
 \eqref{eqn:UIsom} is an isomorphism for all $n$.
\end{proof}

\section{Calculation of $\wwbar{H}_*(2 \pi n; 0)$}
\label{sec:HBar}

This section is devoted to computing the homology of the complex
$\wwbar{C}_*(2\pi n;0)$ over $\Z$, defined in \eqref{eqn:CBar}, in
which we mod out by translation of polygons.  This will prove parts
(b) and (c) of \fullref{thm:main}.

\subsection{The universal coefficient spectral sequence}
\label{sec:UCSS}

The homology $\wwbar{H}_*(2\pi n;0)$ is partially but not
entirely determined by $\wwtilde{H}_*(2\pi n;0)$.  The precise
relation between $\wwtilde{H}_*$ and $\wwbar{H}_*$ is expressed
by the ``universal coefficient spectral sequence'', whose construction
we now recall.

In general let $R$ be a commutative ring, let $(C_*,\delta)$ be a
chain complex of projective $R$--modules, and let $A$ be an $R$--module.
The task at hand is to relate the homology of $(C_*\tensor_R
A,\delta\tensor 1)$ to the homology of $C_*$.  Let
\[
\cdots \stackrel{\partial}{\longrightarrow} P_2
\stackrel{\partial}{\longrightarrow} P_1
\stackrel{\partial}{\longrightarrow} P_0 \longrightarrow A
\longrightarrow 0
\]
be a projective resolution of $A$ in the category of $R$--modules.
Recall that if $M$ is another $R$--module then $\op{Tor}_*(M,A)$ is
defined to be the homology of the complex $(M\tensor_R P_*,
1\tensor\partial)$.  This satisfies $\op{Tor}_0(M,A)=M\tensor_R A$, and
if $M$ is projective then $\op{Tor}_i(M,A)=0$ for $i>0$.  Now form the
double complex
\[
C_{i,j} \eqdef C_j \tensor_R P_i.
\]
This has horizontal, vertical, and total differentials
\[
d_h\eqdef 1\tensor\partial,\quad\quad d_v\eqdef \delta\tensor 1, \quad
\quad d\eqdef d_h+(-1)^id_v.
\]
Filtering the double complex by $j$ gives a spectral sequence with
\[
E^1_{i,j} = \op{Tor}_i(C_j,A) = \left\{\begin{array}{cl} C_j\tensor_R
A, & i=0,\\ 0, & i>0,
\end{array}\right.
\]
so the homology of the double complex is $H_*(C_*\tensor_R A)$.  The
universal coefficient spectral sequence is obtained by filtering the
double complex by $i$ instead.  This spectral sequence satisfies
\[
\begin{split}
E^1_{i,j} &= H_j(C_*) \tensor_R P_i,\\
E^2_{i,j} &= \op{Tor}_i\left(H_j(C_*),A\right),
\end{split}
\]
and by the previous calculation converges to $H_*(C_*\tensor_R A)$.

We now specialize to the case $R=\Z[\Z^2]$,
$C_*=\wwtilde{C}_*(2\pi n;0)$, and $A=\Zmodule$.  In the rest of this
section, the tensor product is understood to be over $\Z[\Z^2]$.

 We begin by computing the relevant Tor's, ie the $E^2$ term of the
universal coefficient spectral sequence.  We fix the following projective
resolution of $\Zmodule$:
\begin{gather*}
0 \longrightarrow \Z[\Z^2]\{\gamma\} \stackrel{\partial}{\longrightarrow}
 \Z[\Z^2]\{\alpha,\beta\}
\stackrel{\partial}{\longrightarrow} \Z[\Z^2]\{\tau\} {\longrightarrow}
\Zmodule \longrightarrow
0,\\
\partial\gamma\eqdef (y-1)\alpha - (x-1)\beta, \quad\quad \partial\alpha\eqdef
(x-1)\tau, \quad\quad \partial\beta \eqdef (y-1)\tau.
\end{gather*}
After tensoring with $\Zmodule$, the differential becomes $0$, so 
\begin{equation}
\label{eqn:tor1}
\op{Tor}_i(\Zmodule,\Zmodule) \simeq \left\{\begin{array}{cl} \Zmodule, &
i=0,\\
\Zmodule\oplus\Zmodule, & i=1,\\
\Zmodule, & i=2,\\
0, & i>2.\end{array}\right.
\end{equation}
We compute $\op{Tor}_*(\mc{I}(\Z^2),\Zmodule)$ by the symmetry
$\op{Tor}_*(\mc{I}(\Z^2),\Zmodule)=\op{Tor}_*(\Zmodule,\mc{I}(\Z^2))$.
The $\Z[\Z^2]$--module $\mc{I}(\Z^2)$ has a presentation with two
generators $a=x-1$ and $b=y-1$ and the single relation
$(y-1)a-(x-1)b=0$.  Therefore
\begin{equation}
\label{eqn:tor2}
\op{Tor}_i(\mc{I}(\Z^2),\Zmodule)
=\op{Tor}_i(\Zmodule,\mc{I}(\Z^2))
\simeq
\left\{\begin{array}{cl} \Zmodule \oplus \Zmodule, &i=0,\\
\Zmodule, &i=1,\\
0, &i>0.
\end{array}\right.
\end{equation}
In terms of the projective resolution of $\Zmodule$, it turns out that
$\op{Tor}_0(\mc{I}(\Z^2),\Zmodule)$ is generated by $(x-1)\tensor\tau$
and $(y-1)\tensor\tau$, while $\op{Tor}_1(\mc{I}(\Z^2),\Zmodule)$ is
generated by $(y-1)\tensor\alpha-(x-1)\tensor\beta$.

Now the only possibly nonzero higher differential in the universal
coefficient spectral sequence is, for $j\ge 0$,
\begin{equation}
\label{eqn:UCSSDifferential}
\begin{split}
d_2\co E^2_{2,j} &\longrightarrow E^2_{0,j+1},\\
\Zmodule \simeq \op{Tor}_2\left(\wwtilde{H}_j(2\pi n;0),\Zmodule\right)
& \longrightarrow
\op{Tor}_0\left(\wwtilde{H}_{j+1}(2\pi n;0),\Zmodule\right)
\simeq \Zmodule.
\end{split}
\end{equation}
To compute this differential, the only explicit calculation we will
need is given by the following lemma, which will be proved in
\fullref{sec:explicit}.

\begin{lemma}
\label{lem:explicit}
If $p$ is a cycle generating $\wwtilde{H}_0^{(2-2n)}(2\pi
n;0)\simeq\Zmodule$, and if $s$ and $t$ are chains with $\delta
s = (x-1)p$ and $\delta t = (y-1)p$, then $(y-1)s - (x-1)t$
generates $\wwtilde{H}_1^{(2-2n)}(2\pi n;0)\simeq\Zmodule$.
\end{lemma}

\begin{lemma}
\label{lem:d2}
The differential \eqref{eqn:UCSSDifferential} is zero if $j$ is odd,
and an isomorphism if $j\ge 0$ is even.
\end{lemma}

\begin{proof}
We first derive a general formula for $d_2$.  Let
$p\in\wwtilde{C}_j(2\pi n;0)$ be a cycle generating the $\Zmodule$
component of $\wwtilde{H}_j(2\pi n;0)$.  Referring back to the
generator $\gamma$ of our projective resolution of $\Zmodule$, we want
to compute $d_2[p\tensor\gamma]$.  This is done by ``zig-zagging''.
First we calculate the horizontal differential
\[
d_h(p\tensor\gamma) = p\tensor ((y-1)\alpha-(x-1)\beta)\in C_{1,j}.
\]
Next we need to find an element
\[
-t\tensor\alpha + s\tensor \beta \in C_{1,j+1}
\]
with $-d_v(-t\tensor\alpha+s\tensor\beta)=d_h(p\tensor\gamma)$.  That
is, we need to find $s,t\in\wwtilde{C}_{j+1}(2\pi n;0)$ with
\begin{equation}
\label{eqn:pCondition}
\delta s =
(x-1)p, \quad\quad \delta t = (y-1)p
\end{equation}
in $\wwtilde{C}_j(2\pi n;0)$.  For $s$ and $t$ satisfying
\eqref{eqn:pCondition}, we then have
\begin{equation}
\label{eqn:completeComputation}
\begin{split}
d_2[p\tensor\gamma] &= [d_h(-t\tensor\alpha + s\tensor\beta)]\\
&= [(y-1)s-(x-1)t)\tensor\tau].
\end{split}
\end{equation}
If $j=0$, then \fullref{lem:explicit} and equation
\eqref{eqn:completeComputation} imply that
$d_2[p\tensor\gamma]$ is a generator of 
$\op{Tor}_0\left(\wwtilde{H}_{1}(2\pi n;0),\Zmodule\right) \simeq
\Zmodule$.

Suppose next that $j=2k$.  Since the chain map $U$ commutes with $x$
and $y$ according to \fullref{prop:UProperties}(b), it follows
by functoriality or by equation \eqref{eqn:completeComputation} that
the induced map on homology $U_*$ commutes with the differential
$d_2$.  By \fullref{thm:HTilden}(b), we are done in this case by
induction on $k$.

Finally, if $j=2k+1$, we show that the differential
\eqref{eqn:UCSSDifferential} is zero by considering the decomposition
into subcomplexes from \fullref{definition:subcomplex}.  First,
\fullref{thm:HTilden}(a) implies that
\[
\wwtilde{H}_j(2\pi n;0) = \wwtilde{H}_j^{(2k-2n+2)}(2\pi n;0),
\]
so we can take $p\in \wwtilde{C}_j^{(2k-2n+2)}(2\pi n;0)$ and
$s,t\in \wwtilde{C}_{j+1}^{(2k-2n+2)}(2\pi n;0)$.  But
\fullref{thm:HTilden}(a) also implies that
\[
\wwtilde{H}_{j+1}(2\pi n;0) = \wwtilde{H}_{j+1}^{(2k-2n+4)}(2\pi n;0).
\]
Hence the cycle $(y-1)s-(x-1)t$ is nullhomologous in
$\wwtilde{C}_*(2\pi n;0)$, so by equation
\eqref{eqn:completeComputation}, $d_2[p\tensor\gamma]=0$.
\end{proof}

\begin{proof}[Proof of \fullref{thm:main}(b)]  By equations
\eqref{eqn:tor1} and \eqref{eqn:tor2} and \fullref{lem:d2}, the
$E^2$ term of the universal coefficient spectral sequence looks like
this:
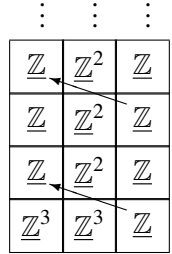
\begin{center}
\begin{picture}(60,95)(0,0)

\put(0,0){\framebox(20,20){$\underline{\Z}^3$}}
\put(20,0){\framebox(20,20){$\underline{\Z}^3$}}
\put(40,0){\framebox(20,20){$\underline{\Z}$}}

\put(0,20){\framebox(20,20){$\underline{\Z}$}}
\put(20,20){\framebox(20,20){$\underline{\Z}^2$}}
\put(40,20){\framebox(20,20){$\underline{\Z}$}}

\put(0,40){\framebox(20,20){$\underline{\Z}$}}
\put(20,40){\framebox(20,20){$\underline{\Z}^2$}}
\put(40,40){\framebox(20,20){$\underline{\Z}$}}

\put(0,60){\framebox(20,20){$\underline{\Z}$}}
\put(20,60){\framebox(20,20){$\underline{\Z}^2$}}
\put(40,60){\framebox(20,20){$\underline{\Z}$}}

\put(44,16){\vector(-3,1){29}}
\put(44,56){\vector(-3,1){29}}

\put(10,85){$\vdots$}
\put(30,85){$\vdots$}
\put(50,85){$\vdots$}

\end{picture}
\end{center}
Here the $d_2$ arrows drawn are isomorphisms, while the others are
zero.  It follows immediately that the spectral sequence converges to
$\underline{\Z}^3$ in each nonnegative degree.  This is
$\wwbar{H}_*(2\pi n;0)$ as a $\Z[\Z^2]$--module; as a $\Z$--module, it
is simply $\Z^3$ in each nonnegative degree.
\end{proof}

More explicitly, the proof of \fullref{thm:main}(b) shows the following.

\begin{proposition}
\label{prop:generatingCycles}
For $k\ge 0$, let $p_k$ be a cycle generating
$\wwtilde{H}_{2k}^{(2k-2n+2)}(2\pi n;0)$, and let $q_k$
be a cycle generating $\wwtilde{H}_{2k+1}^{(2k-2n+2)}(2\pi n;0)$.
Then $\wwbar{H}_*(2\pi n;0)$ is freely generated over $\Z$ by the
images in $\wwbar{C}_*(2\pi n;0)$ of the following chains in
$\wwtilde{C}_*(2\pi n;0)$:
\begin{itemize}
\item The index $2k$ cycle $p_k$, for each $k\ge 0$.
\item Two index $2k+1$ chains $s_k$ and $t_k$ with $\delta
  s_k=(x-1)p_k$ and $\delta t_k=(y-1)p_k$, for
  each $k\ge 0$.
\item The index $0$ cycles $u_0\eqdef Z_n((1,0),(0,0))$ and
  $v_0\eqdef Z_n((0,1),(0,0))$.
\item Two index $2k$ chains $u_k$ and $v_k$ with $\delta u_k =
  (x-1)q_{k-1}$ and $\delta v_k = (y-1)q_{k-1}$, for each
  $k\ge 1$.
\item An index $2k+1$ chain $w_k$ with $\delta w_k =
  (y-1)u_k-(x-1)v_k$ for each $k\ge 0$.
\end{itemize}
Moreover, the corresponding homology classes in
$\wwbar{H}_*(2\pi;0)$ do not depend on any of the choices of
$p_k,q_k,s_k,t_k,u_k,v_k,w_k$.
\qed
\end{proposition}

Note for example that the cycle $q_k$ maps to zero in
$\wwbar{H}_*(2\pi n;0)$, because we saw above that $\pm q_k$ is
homologous in $\wwtilde{C}_*(2\pi n;0)$ to $(y-1)s_k-(x-1)t_k$, and
$(y-1)$ and $(x-1)$ are in the kernel of the augmentation map
$\Z[\Z^2]\to\Z$.

\subsection{The action of $U$ on $\wwbar{H}_*(2\pi n;0)$}

We now prove \fullref{thm:main}(c).  Fix $\theta\in\R/2\pi n\Z$
with $\tan\theta$ irrational in the definition of $U$.  We
will use the following lemma which is proved in \fullref{sec:explicit}.

\begin{lemma}
\label{lem:explicitU}
There exist $q_0\in \wwtilde{C}_1^{(2-2n)}(2\pi n;0)$
generating $\wwtilde{H}_1^{(2-2n)}(2\pi n;0)$ and $u_1, v_1
\in
\wwtilde{C}_2^{(2-2n)}(2\pi n;0)$ such that
\begin{gather}
\label{eqn:URequirement1}
\delta u_1=(x-1)q_0, \quad \delta v_1 = (y-1)q_0,\\
\label{eqn:URequirement2}
Uu_1=Z_n((1,0),(0,0)),\quad Uv_1=Z_n((0,1),(0,0)).
\end{gather}
\end{lemma}

\begin{proof}[Proof of \fullref{thm:main}(c)]  It is enough to show that
for any $k\ge 0$, one can make the choices in
\fullref{prop:generatingCycles} for $k$ and $k+1$ such that
\begin{gather}
\label{eqn:U1}
Up_{k+1} = p_k, \quad Us_{k+1} = s_k, \quad
Ut_{k+1} = t_k,\\
\label{eqn:U2}
Uu_{k+1} = u_k, \quad Uv_{k+1}=v_k, \quad
Uw_{k+1} = w_k
\end{gather}
in $\wwtilde{C}_*(2\pi n;0)$.

First choose any $p_{k+1}, q_{k+1}, s_{k+1}, t_{k+1}, u_{k+1}, v_{k+1},
w_{k+1}$.  By \fullref{thm:HTilden}(b), we can choose
$p_k\eqdef Up_{k+1}$ and $q_k\eqdef Uq_{k+1}$.  Since $U$ is a
translation-invariant chain map, we can then choose $s_k\eqdef Us_{k+1}$ and
$t_k\eqdef Ut_{k+1}$, and if $k>0$ we can also choose $u_k\eqdef Uu_{k+1}$,
$v_k\eqdef Uv_{k+1}$, and $w_k\eqdef Uw_{k+1}$.

To complete the proof, it is enough to show that for suitable choices,
\begin{equation}
\label{eqn:suitableChoices}
Uu_{1}=u_0,\quad Uv_{1}=v_0,\quad Uw_1=w_0.
\end{equation}
We can obtain the first two conditions in \eqref{eqn:suitableChoices}
by \fullref{lem:explicitU}.  Then to obtain the third condition,
given any choice of $w_1$, we can choose $w_0\eqdef Uw_1$.
\end{proof}

\subsection{Some explicit homology generators and relations}
\label{sec:explicit}

This subsection is devoted to the proofs of Lemmas~\ref{lem:explicit}
and \ref{lem:explicitU} above.

\begin{definition}
\label{def:abtheta}
Let $a,b\in\Z^2$ be lattice points with $a-b$ indivisible.  Choose
$\theta\in\R$ such that the line from $a$ to $b$ has angle $\theta\mod
2\pi$.
\begin{itemize}
\item
Define
\[
p(a,\theta)\in\wwtilde{C}_0^{(2-2n)}(2\pi n;0)
\]
to be the generator which wraps $n-1$ times around the $2$--gon between
$a$ and $b$, with edges at angles $\theta, \theta+\pi,
\ldots\theta+(2n-3)\pi$, all labeled `$h$', in that order. 
Note that if $n=1$, then $p(a,\theta)$ is the constant path at $a$,
while if $n>1$, then $p(a,\theta)$ has a kink parametrized by the
interval $(\theta-3\pi,\theta)$ and mapping to $a$.  Also, $\delta
p(a,\theta) = 0$.
\item Define
\[
e(a,\theta)\in \wwtilde{C}_1^{(2-2n)}(2\pi n;0)
\]
  to be the generator obtained from $Z_n(a,b)$, see \fullref{sec:zero},
  by relabeling the edge from $b$ to $a$ at angle $\theta-\pi$ by
  `$e$' and ordering the $2n-1$ `$h$' edges counterclockwise.  Then
\begin{equation}
\label{eqn:deltaf}
\delta e(a,\theta) = p(b,\theta+\pi) - p(a,\theta).
\end{equation}
\item
Let
\[
q(a,b) \in \wwtilde{C}_1^{(2-2n)}(2\pi n;0)
\]
be the sum of all $2n$ generators that wrap $n$ times around the
$2$--gon between $a$ and $b$ with $2n-1$ edges labeled `$h$' and
ordered counterclockwise.  That is,
\[
q(a,b) \eqdef \sum_{i=0}^{n-1} \big( e(a,\theta+2i\pi) +
e(b,\theta+(2i+1)\pi) \big).
\]
Note that $q(a,b)=q(b,a)$.  By equation \eqref{eqn:deltaf}, $\delta
q(a,b)=0$.
\end{itemize}
\end{definition}

\begin{lemma}
\label{lem:simpleGenerators}
\begin{enumerate}
\item[(a)]
$\wwtilde{H}_0^{(2-2n)}(2\pi
n;0)\simeq\Zmodule$ is generated by $p(a,\theta)$.
\item[(b)]
$\wwtilde{H}_1^{(2-2n)}(2\pi
n;0)\simeq\Zmodule$ is generated by $q(a,b)$.
\end{enumerate}
\end{lemma}

\begin{proof}
  (a) This follows from the proof of \fullref{thm:HTilden}(a),
  since by the symmetry in \fullref{lem:symmetry} we may assume that
  $a$ and $b$ are on the $x$--axis with $\theta=2\pi$, and then the generator
  $p(a,\theta)$ is obtained by applying the splicing map $(n-1)$ times
  to a constant path which generates $\wwtilde{H}_0^{(0)}(2\pi;0)$.
  
  (b) Let $\Lambda'$ be the $n$--convex path that wraps $n$ times
  around the $2$--gon with corners $a$ and $b$.  Then
  $H_1^{(2-2n)}(\Lambda')\simeq\Z$ is generated by $q(a,b)$.  Indeed
  the summands in $q(a,b)$ are the only generators in
  $C_1^{(2-2n)}(\Lambda')$, since such a generator must have $2n-1$
  `$h$' edges, and their sum is the only cycle in
  $C_1^{(2-2n)}(\Lambda')$ by equation \eqref{eqn:deltaf}.  This
  cycle is not a boundary since $C_2^{(2-2n)}(\Lambda')=0$, as the
  only generator in $C_*(\Lambda')$ with $2n$ `$h$' edges has index
  zero.
      
  To finish the proof, we claim that the inclusion-induced map
  $$H_1^{(2-2n)}(\Lambda')\to H_1^{(2-2n)}(2\pi n;0)$$
  is an isomorphism.  This follows from the spectral sequence in the proof
  of \fullref{lem:injectivity}.  The reason is that in that spectral
  sequence, with $k$ and $m$ as in that proof, if $k\ge 3$ then there
  are no $m>0$ generators in the $E^1$ term in $C_1^{(2-2n)}(\Lambda)$
  or $C_2^{(2-2n)}(\Lambda)$.  Such a homology generator would have
  to be a sum of chain complex generators each having $2n-1$ or $2n$
  `$h$' edges, which means that it would be either a product of $n{-}1$
  $Z_1$'s and $1$ $H$ (which has index $2k-3\ge 3$), or a product of $n$
  $Z_1$'s (which has index $0$).
\end{proof}

\begin{lemmadef}
\label{lem-def:e}
  There exists a unique assignment, to each triple
  $(a,\theta,\theta')$ with $a\in\Z^2$ and $\theta,\theta'\in\R$ such that
  $\tan\theta,\tan\theta'\in\Q\cup\{\infty\}$ and $\theta\le\theta'$, of an
  equivalence class of chains
\[
f(a,\theta,\theta') \in
\frac{\wwtilde{C}_1^{(2-2n)}(2\pi n;0)}{\op{Im}(\delta)}
\]
such that:
\begin{enumerate}
\item[(i)]
For each $(a,\theta,\theta')$ as above,
\begin{equation}
\label{eqn:deltae}
\delta f(a,\theta,\theta') = p(a,\theta) -
p(a,\theta').
\end{equation}
\item[(ii)] Suppose that $0<\theta'-\theta\le\pi$, and that $b$ and
  $b'$ are defined from $(a,\theta)$ and $(a,\theta')$ as in
  \fullref{def:abtheta}.  Then $f(a,\theta,\theta')\in
  C_*(\Lambda)$, where $\Lambda$ wraps $n-1$ times around the triangle
  with vertices $a,b,b'$ (or $2$--gon with vertices $b$ and $b'$ when
  $\theta'-\theta=\pi$) and has a kink at $a$ parametrized by the
  interval $(\theta'-3\pi,\theta)$.
\item[(iii)]
If $\theta \le \theta' \le \theta''$, then
\begin{equation}
\label{eqn:eAdditive}
f(a,\theta,\theta'') = f(a,\theta,\theta') +
f(a,\theta',\theta'') \mod \op{Im}(\delta).
\end{equation}
\end{enumerate}
\end{lemmadef}

\begin{proof} The proof has four steps.
  
  \textbf{Step 1}\qua We first show that if $\theta'-\theta<\pi$ and if
  $(a,b,b')$ is a simple triangle (see
  \fullref{def:simpleTriangle}), then there exists a chain
  $f(a,\theta,\theta')$ satisfying (i) and (ii).  Reintroduce the
  notation from the proof of \fullref{lem:znRelations}, with
  $c\eqdef b'$.  Let $\mc{C}'$ denote the set of subsets
  $I\subset\{1,\ldots,3n\}$ such that for distinct $i,j\in I$, the
  corners $c_i$ and $c_j$ of $\lambda$ are not adjacent, except that
  we allow at most one adjacent pair involving $c_1$ or $c_{3n}$.  For
  $I=\{i_1,\ldots,i_k\}\in\mc{C}'$, let $T'(I)\in C_*(\Lambda)$ be the
  generator with underlying path $\lambda(I)$ and with all edges
  labeled `$h$', except that if $I$ does not contain an adjacent pair,
  then the edge that starts at $c_{3n}$ and/or ends at $c_1$ is
  labeled `$e$'.  Order the `$h$' edges counterclockwise, starting at
  $c_1$ if $1\notin I$, at $c_2$ if $1\in I$ and $2\notin I$, and at
  $c_3$ if $1,2\in I$.  Let $\theta_i$ denote the angle of the edge
  between corners $c_i$ and $c_{i+1}$.  We may assume that
  $\theta=\theta_1$ and $\theta'=\theta_0+\pi$.  Then in this
  notation,
\begin{equation}
\label{eqn:pNotation}
\begin{split}
p(a,\theta_1) & = T'(\{3,6,\ldots,3n-3,3n-1,3n\}),\\
p(a,\theta_0+\pi) & = T'(\{2,5,\ldots,3n-4,3n-1,3n\}).
\end{split}
\end{equation}
The differential of a generator $T'(I)$ is given by
\begin{equation}
\label{eqn:deltaT'}
\delta T'(I) = - \sum_{I\cup\{i\}\in\mc{C}'} (-1)^{\#\{j\notin I\mid
i<j\}} T'\left(I\cup\{i\}\right).
\end{equation}
Now let $\mc{C_0}'$ denote the set of
$I=\{i_1,\ldots,i_{m}\}\in\mc{C}'$ with $i_1<i_2<\cdots<i_{m}$
such that the $i_k$'s alternate parity with $i_1$ odd.  Define
\begin{equation}
\begin{aligned}
\label{eqn:eDefStep1}
f(a,\theta_1,\theta_0+\pi) &\eqdef
\sum_{I=\{i_1,\ldots,i_{n-2},3n-1,3n\}\in\mc{C}_0'}
T'(I)\\
&= T(\{3,6,\ldots,3n-9,3n-6,3n-1,3n\})\\
&\quad+ T(\{3,6,\ldots,3n-9,3n-4,3n-1,3n\})\\
&\quad+ \cdots +
T(\{5,8,\ldots,3n-7,3n-4,3n-1,3n\}).
\end{aligned}
\end{equation}
Then it follows from equations \eqref{eqn:pNotation} and
\eqref{eqn:deltaT'} that this satisfies condition (i), ie 
\[
\delta f(a,\theta_1,\theta_0+\pi) = p(a,\theta_1) -
p(a,\theta_0+\pi).
\]
Also, condition (ii) is satisfied since each term on the right side of
\eqref{eqn:eDefStep1} has the corners $c_{3n-1}$ and $c_{3n}$ rounded.

\textbf{Step 2}\qua We now show that if $0<\theta'-\theta<\pi$, then a chain
$f(a,\theta,\theta')$ satisfying (i) and (ii), if such exists, is
unique modulo $\op{Im}(\delta)$.  The difference between any two such
chains $f(a,\theta,\theta')$ is a cycle in $C_1^{(2-2n)}(\Lambda)$.
Thus it is enough to show that $H_1^{(2-2n)}(\Lambda)=0$.  Pick an
angle $\phi$ with irrational tangent between $\theta'-\pi$ and
$\theta$.  \fullref{prop:III} shows that
$H_1^{(2-2n)}(\Lambda) \simeq H_1^{(2-2n)}(C_*^\phi(\Lambda))$.  As in
\fullref{lem:CICC} and \fullref{example:xtheta},
$H_1^{(2-2n)}(C_*^\phi(\Lambda)) \simeq H_1^{(2-2n)}(\Lambda_0)$ where
$\Lambda_0$ wraps $n-1$ times around a $2$--gon and has a kink.  But
$C_1^{(2-2n)}(\Lambda_0)=0$ because a generator of
$C_1^{(2-2n)}(\Lambda_0)$ would have $2n-1$ edges labeled `$h$', but
generators of $C_*(\Lambda_0)$ have at most $2n-2$ edges.

\textbf{Step 3}\qua We now show that there exists an assignment
$f(a,\theta,\theta')$ satisfying (i), (ii), and (iii).  Let
$(a,\theta,\theta')$ be given.  If $\theta=\theta'$, define
$f(a,\theta,\theta') \eqdef 0$.  Otherwise choose
$\theta=\theta_0<\theta_1<\cdots<\theta_k=\theta'$ such that
$f(a,\theta_{i-1},\theta_{i})$ is defined by Step 1 for
$i=1,\ldots,k$.  Then define
\begin{equation}
\label{eqn:eDefinition}
f(a,\theta,\theta') \eqdef \sum_{i=1}^k
f(a,\theta_{i-1},\theta_i).
\end{equation}
As long as this is well-defined modulo $\op{Im}(\delta)$, it clearly
satisfies (i) and (ii) (by Step 1) and (iii) (by construction).

To show that \eqref{eqn:eDefinition} is well-defined modulo
$\op{Im}(\delta)$, let
$\theta=\theta_0'<\theta_1'<\cdots<\theta'_{k'}=\theta'$ be another
set of choices to define $f(a,\theta,\theta')$.  We need to
show that
\begin{equation}
\label{eqn:eSumDefined}
\sum_{i=1}^k f(a,\theta_{i-1},\theta_i)
=
\sum_{i=1}^{k'} f(a,\theta'_{i-1},\theta'_i)
 \mod\op{Im}(\delta).
\end{equation}
Without loss of generality, $\theta_1 < \theta_1'$.  By Step 2,
\begin{equation}
\label{eqn:byStep2}
f(a,\theta_0,\theta_1')
=
f(a,\theta_0,\theta_1) + 
f(a,\theta_1,\theta_1') \mod \op{Im}(\delta).
\end{equation}
Subtracting \eqref{eqn:byStep2} from \eqref{eqn:eSumDefined}, we see
that to prove equation \eqref{eqn:eSumDefined}, it is enough to show
that $f(a,\theta_1,\theta')$ is well defined modulo
$\op{Im}(\delta)$.  We are now done by induction on $k+k'$.

\textbf{Step 4}\qua An assignment $f(a,\theta,\theta')$ satisfying
(i), (ii), and (iii) is unique modulo $\op{Im}(\delta)$, because
condition (iii) forces it to satisfy equation \eqref{eqn:eDefinition},
and each summand on the right hand side of \eqref{eqn:eDefinition} is
unique modulo $\op{Im}(\delta)$ by conditions (i) and (ii) and Step 2.
\end{proof}

\begin{lemma}
\label{lem:sixTermRelation}
If $(a,b,c)$ is a simple triangle, then with $\theta_i$ defined as above,
\begin{equation}
\label{eqn:sixTermRelation}
\begin{split}
f(c,\theta_3,\theta_2+\pi) & + f(b,\theta_2,\theta_1+\pi) +
f(a,\theta_1,\theta_0+\pi) = \\ & = -
e(b,\theta_2) + e(a,\theta_0+\pi) -
e(a,\theta_1) \mod \op{Im}(\delta).
\end{split}
\end{equation}
\end{lemma}

\begin{proof}
In the notation of the previous proof, define
\begin{equation}
\label{eqn:T'}
T'_k \eqdef \sum_{I=\{i_1,\ldots,i_k\}\in\mc{C}_0'} T'(I).
\end{equation}
Equation \eqref{eqn:deltaT'} implies that
\begin{equation}
\label{eqn:deltaTn'}
\delta\bigl(-T'_{n-1}\bigr)
= T'_n - T'\bigl(\{2,5,\ldots,3n-1\}\bigr)
\end{equation}
because on the right hand side, terms $T'(I)$ in which the indices in
$I$ do not alternate parity will appear twice with opposite signs or
not at all, while terms $T'(I)$ in which the indices do alternate
parity will appear exactly once, and the only way to alternate parity
starting with an even index is $2,5,\ldots,3n-1$.  Now the right hand
side of equation \eqref{eqn:deltaTn'} equals
\begin{multline*}
\biggl(\sum_{I=\{1,2,i_3,\ldots,i_{n}\}\in\mc{C}_0'}{+}
\sum_{I=\{1,i_2,\ldots,i_{n-1},3n\}\in\mc{C}_0'} {+}
\sum_{I=\{i_1,\ldots,i_{n-2},3n-1,3n\}\in\mc{C}_0'}\biggr)
T'(I) + \\
+ T'(\{1,4,\ldots,3n-2\})
 - T'(\{2,5,\ldots,3n-1\}) + T'(\{3,6,\ldots,3n\}).
\end{multline*}
These six terms equal the six terms in the relation
\eqref{eqn:sixTermRelation}.
\end{proof}

\begin{lemma}
\label{lem:threeTermGenerator}
If $a$, $b$, and $\theta$ are as in \fullref{def:abtheta}, then
\[
f(a,\theta,\theta+2\pi) + e(a,\theta) + e(b,\theta+\pi) =
q(a,b) \in \wwtilde{H}_1^{(2-2n)}(2\pi n;0).
\]
\end{lemma}

\begin{proof}
Without loss of generality,
$$a=\begin{pmatrix}0\\0\end{pmatrix},\quad
b=\begin{pmatrix}1\\0\end{pmatrix},\quad \theta=0.$$
We can now do
the entire calculation on the $x$--axis, ie in the subcomplex
$CX_*(n)$, and use the notation of \fullref{sec:splicing}
to describe chains in this subcomplex.  In this notation,
\[
\begin{split}
e(a,0) &= \left(h_b^a h_a^b\right)^{n-1} e_b^a h_a^b,\\
e(b,\pi) &= \left(h_b^a h_a^b\right)^{n-1} h_b^a e_a^b.
\end{split}
\]
Let
$$c\eqdef\begin{pmatrix}-1\\0\end{pmatrix}.$$
We can take
\[
\begin{split}
f(a,0,\pi) &= \sum_{i=0}^{n-2}
\left(h_b^a h_a^b\right)^i h_b^c \left(h_c^a h_a^c\right)^{n-2-i}
h_c^a e_a^a h_a^b,\\
f(a,\pi,2\pi) &= \sum_{i=0}^{n-2}
h_a^c \left(h_c^a h_a^c\right)^i h_c^b \left(h_b^a
h_a^b\right)^{n-2-i} h_b^a e_a^a,
\end{split}
\]
because the right hand side of each equation satisfies conditions (i)
and (ii) in \fullref{lem-def:e}.  By the definition of
$\delta$ and straightforward manipulation of sums,
\begin{multline*}
\delta\biggl(\sum_{i=0}^{n-2} \sum_{j=0}^{n-2-i} \left(h_b^a
h_a^b\right)^i h_b^c
\left(h_c^a h_a^c\right)^j h_c^b \left(h_b^a h_a^b\right)^{n-1-i-j}\biggr) = \quad\quad \\
-f(a,0,\pi) - f(a,\pi,2\pi) - e(a,0) - e(b,\pi) + q(a,b),
\end{multline*}
which completes the proof.
\end{proof}

\begin{proof}[Proof of \fullref{lem:explicit}] Since this lemma is
computing a differential in the universal coefficient spectral
sequence, it is enough to verify the conclusion of the lemma for a
single choice of $p$, $s$, and $t$.  Introduce the lattice points
$a\eqdef(0,0)$, $b\eqdef(1,0)$, $c\eqdef(2,0)$, $d\eqdef(0,1)$,
$e\eqdef(1,1)$, and $f\eqdef(2,1)$.  By
\fullref{lem:simpleGenerators}(a), we can take
\[
p \eqdef p(a,0).
\]
By equations \eqref{eqn:deltaf} and \eqref{eqn:deltae}, we can take
\[
\begin{split}
s &\eqdef f(b,0,\pi) + e(a,0),\\
t & \eqdef -f(a,0,\pi/2) - e(d,-\pi/2) - f(d,-\pi/2,0).
\end{split}
\]
These chains are only defined mod $\op{Im}(\delta)$, which is fine
here since we just need to evaluate the homology class of
$(y-1)s-(x-1)t$.  By definition,
\begin{equation}
\label{eqn:long}
\begin{split}
(y-1)s -(x-1)t &= f(b,0,\pi/2) + e(e,-\pi/2) +
f(e,-\pi/2,0) \\
& \quad - f(a,0,\pi/2) - e(d,-\pi/2) - f(d,-\pi/2,0)
\\
& \quad + f(e,0,\pi) + e(d,0) - f(b,0,\pi) - e(a,0).
\end{split}
\end{equation}
By \fullref{lem:sixTermRelation},
\[
\begin{split}
f(b,3\pi/4,\pi) & + f(a,0,\pi/2) +
f(d,-\pi/2,-\pi/4) = \\ &= - e(a,0) + e(d,-\pi/4) - e(d,-\pi/2),\\
f(b,\pi/2,3\pi/4) & + f(d,-\pi/4,0) +
f(e,-\pi,-\pi/2) =\\&= - e(d,-\pi/4) + e(e,-\pi/2) - e(e,-\pi).
\end{split}
\]
Putting these two six-term relations into equation \eqref{eqn:long}
and repeatedly applying the relation \eqref{eqn:eAdditive} gives
\[
(y-1)s - (x-1)t = f(e,-\pi,\pi) + e(e,-\pi) + e(d,0).
\]
By Lemmas~\ref{lem:threeTermGenerator} and
\ref{lem:simpleGenerators}(b), this generates
$\wwtilde{H}_1^{(2-2n)}(2\pi n;0)\simeq\Zmodule$.
\end{proof}

\begin{lemma}
Let $(a,b,c)$ be a simple triangle, and let
$\theta_1,\ldots,\theta_{3n}$ be defined as previously.  Then there
exists a chain
\[
r(a,b,c) \in \wwtilde{C}_2^{(2-2n)}(2\pi n;0)
\]
such that if $U$ is defined using $\theta$, then
\begin{align}
\label{eqn:deltar}
\delta r(a,b,c) &= q(a,b) - q(a,c),\\
\label{eqn:Ur}
U r(a,b,c) &= \left\{\begin{array}{cl}
-Z_n(a,b), & \theta\in(\theta_{3i+1}-\pi,\theta_{3i}),\\
Z_n(b,c), & \theta\in(\theta_{3i},\theta_{3i+1}),\\
-Z_n(c,a), & \theta\in(\theta_{3i+1},\theta_{3i}+\pi),\\
0, & \mbox{otherwise.}
\end{array}
\right.
\end{align}
\end{lemma}

\begin{proof}
Suppose first that
\begin{equation}
\label{eqn:thetaAssumption}
\theta\notin (\theta_{3i},\theta_{3i+2}-\pi), \,
(\theta_{3i+1},\theta_{3i}+\pi), \quad i=0,\ldots,n-1.
\end{equation}
Let $\Lambda$ be the $n$--convex path that wraps $n$ times around the
triangle $(a,b,c)$.  There is an obvious action of the cyclic group
$\Z/3n$ on the chain complex $C_*(\Lambda)$, given by a chain map
\[
\eta\co C_*(\Lambda) \to C_*(\Lambda)
\]
which rotates everything
counterclockwise, replacing $\theta_i$ by $\theta_{i+1}$, etc.  Define
\begin{equation}
\label{eqn:defr}
r(a,b,c) \eqdef (\eta-1)\sum_{i=0}^{n-1}\eta^i T_{n-1}'
\end{equation}
where $T_{n-1}'$ is defined in equation \eqref{eqn:T'}.  We saw in the
proof of \fullref{lem:sixTermRelation} that
\begin{equation}
\label{eqn:STR}
\begin{split}
\delta\left(-T_{n-1}'\right) &= f(c,\theta_3,\theta_2+\pi) +
f(b,\theta_2,\theta_1+\pi) +
f(a,\theta_1,\theta_0+\pi) + 
\\ & \quad +
e(b,\theta_2) - e(a,\theta_0+\pi) +
e(a,\theta_1).
\end{split}
\end{equation}
In $\delta r(a,b,c)$, all the $f$'s cancel and
\[
\begin{split}
\delta r(a,b,c) &= \sum_{i=0}^{n-1}\Big(e(a,\theta_1+2\pi i) +
e(b,\theta_1+\pi+2\pi i) - \\
& \quad \quad\quad\quad - e(a,\theta_0+\pi+2\pi i) -
e(c,\theta_0+2\pi + 2\pi i)\Big)\\
&= q(a,b) - q(a,c).
\end{split}
\]
This proves \eqref{eqn:deltar}. To prove \eqref{eqn:Ur}, observe from
the definition of $U$ that
\[
U T'(I) = \left\{\begin{array}{cll}
    T\left(\{1\}\cup I\right), & 3n,1,2\notin I; &
    \theta\in(\theta_0,\theta_2-\pi),\\
    T(I\cup\{3n\}), & 3n-1,3n,1\notin I; &
    \theta\in(\theta_1-\pi,\theta_0),\\
    0, & \mbox{otherwise.} &
\end{array}
\right.
\]
In particular,
\[
U T_{n-1}' = \left\{\begin{array}{cl} T\left(\{3,6,\ldots,3n\}\right),
& \theta\in (\theta_1-\pi,\theta_0),\\
0, & \theta\in (\theta_2-\pi,\theta_{3n+1}-\pi).
\end{array}
\right.
\]
Note that $U T_{n-1}'$ is more complicated when $\theta\in
(\theta_0,\theta_2-\pi)$. But under our assumption
 \eqref{eqn:thetaAssumption} on $\theta$, we do not have to consider
 that case, and the above gives equation \eqref{eqn:Ur}.
 
 If the assumption \eqref{eqn:thetaAssumption} does not hold, then
 redefine $T_{n-1}'$ by summing over sequences that start with an even
 index instead of an odd one.  Then \eqref{eqn:STR} holds with a minus
 sign.  If we redefine $r(a,b,c)$ with a minus sign in equation
 \eqref{eqn:defr}, then the rest of the argument goes through.
\end{proof} 

\begin{proof}[Proof of \fullref{lem:explicitU}]
Let $a,b,c,d,e,f$ be the
lattice points in the proof of \fullref{lem:explicit}.
By \fullref{lem:simpleGenerators}(b), $\wwtilde{H}_1^{(2-2n)}(2\pi
n;0)$ is generated by
\[
q_0 \eqdef q(a,b).
\]
By equation \eqref{eqn:deltar}, the requirement
\eqref{eqn:URequirement1} is satisfied by
\[
\begin{split}
u_1 & \eqdef r(b,c,e) + r(b,e,a),\\
v_1 & \eqdef -r(a,b,d) - r(d,a,e).
\end{split}
\]
Choose $\theta\in(\pi/2,\pi)$.  Then by equation \eqref{eqn:Ur},
\begin{align*}
Ur(b,c,e) &= 0,\\
Ur(b,e,a) &= -Z_n(a,b)=Z_n(b,a),\\
Ur(a,b,d) &= 0,\\
Ur(d,a,e) &= -Z_n(d,a)=Z_n(a,d),
\end{align*}
so \eqref{eqn:URequirement2} holds as well.
\end{proof}

\section{Axioms for the chain complex}
\label{sec:uniqueness}

In this section we prove that the chain complex $(\wwtilde{C}_*(2\pi
n;\Gamma),\delta)$ defined in \fullref{sec:mgpc} is characterized by
certain axioms.  This will be used in \fullref{sec:ECH} to relate the
chain complex to the embedded contact homology of $T^3$.  It will
simplify some arguments below to consider all $\Gamma$ at once, so
introduce the notation
\begin{equation}
\label{eqn:gammaSum}
\wwtilde{C}_*(2\pi n) \eqdef
\bigoplus_{\Gamma\in\Z^2}\wwtilde{C}_*(2\pi n;\Gamma).
\end{equation}

\subsection{The axioms}
\label{sec:axioms}

Fix a positive integer $n$.  We now list a series of axioms for a
chain complex $(C_*,\partial)$ over $\Z[\Z^2]$.

\begin{enumerate}
\item[I] (Generators)\qua
$C_*=\wwtilde{C}_*(2\pi n)$ as a $\Z[\Z^2]$--module.
\item[II] (Index)\qua $\partial$ respects the decomposition
  \eqref{eqn:gammaSum} and has degree $-1$ with respect to the
  relative grading $I$ on $\wwtilde{C}_*(2\pi n;\Gamma)$ defined in
  \fullref{sec:mgpc}.
\end{enumerate}

To state the next axioms, let $\alpha$ and $\beta$ be generators of $C_*$ with
the same period $\Gamma$ and with underlying admissible paths
$\lambda$ and $\mu$.  By an ``edge'' of $\alpha$ or $\beta$, we mean
an edge of the corresponding admissible path $\lambda$ or $\mu$.  We
write $\beta\le\alpha$ if $\mu\le\lambda$, see
\fullref{sec:partialOrder}.

We say that two edges of $\alpha$ and $\beta$ ``agree'' (resp.\ 
``partially agree'') if they correspond to the same angle
$\theta\in\R/2\pi n\Z$, and if their adjacent corners map to the same
points in $\Z^2$ (resp. to points on the same line in $\R^2$) (for a
given lift of $\theta$ to $\R$ when $\Gamma\neq 0$).

Let $D(\alpha,\beta)$ denote the closure of the set of all
$t\in\R/2\pi n\Z$ such that $\lambda(t)$ and $\mu(t)$ are defined and
unequal (for a given lift of $t$ to $\R$ when $\Gamma\neq 0$).  Note
that a point $\theta\in\R/2\pi n\Z$ corresponding to an edge of
$\alpha$ (resp.\ $\beta$) is in $D(\alpha,\beta)$ if and only if this
edge does not agree with any edge of $\beta$ (resp.\ $\alpha$).

\begin{enumerate}
\item[III] (Nesting)\qua If $\langle\partial\alpha,\beta\rangle\neq 0$, then
$\beta\le\alpha$.
\item[IV] (Label Matching)\qua Suppose that
  $\langle\partial\alpha,\beta\rangle\neq 0$.  If two edges of
  $\alpha$ and $\beta$ agree, then the labels (`$e$' or `$h$') of the
  two edges are the same.  If two edges of $\alpha$ and $\beta$
  partially agree, and if the edge of $\beta$ is labeled `$h$', then
  the edge of $\alpha$ is also labeled `$h$'.
\item[V] (Connectedness)\qua If $\langle\partial\alpha,\beta\rangle\neq
0$ then the set $D(\alpha,\beta)$ is connected.
\item[VI] (No Double Rounding)\qua If
$\langle\partial\alpha,\beta\rangle\neq 0$ then $\alpha$ cannot have
three edges in $D(\alpha,\beta)$.
\end{enumerate}

The next axiom says essentially that the differential coefficient
$\langle\partial\alpha,\beta\rangle$ depends only on the local change
needed to get from $\alpha$ to $\beta$. To state it, suppose that
$\alpha$ and $\beta$ satisfy the Nesting and Label Matching conditions
above.  We construct generators $\alpha'$ and $\beta'$ as follows.  If
two edges of $\alpha$ and $\beta$ agree, remove them both.  If two
edges of $\alpha$ and $\beta$ partially agree (which by Nesting
implies that the edge of $\beta$ has smaller multiplicity than the
edge of $\alpha$), remove the edge of $\beta$ and shorten the edge of
$\alpha$ by the same amount, while preserving the number of `$h$'
labels.  (That is, if the edge of $\alpha$ is labeled `$h$' and the
edge of $\beta$ is labeled `$e$', then the edge of $\alpha'$ is
labeled `$h$'; otherwise the edge of $\alpha'$ is labeled `$e$'.)  In
particular the multiplicity functions of the underlying admissible
paths $\lambda$, $\mu$, $\lambda'$, $\mu'$ of $\alpha$, $\beta$,
$\alpha'$, $\beta'$ respectively satisfy
\[
m_\alpha - m_{\alpha'} = m_{\beta} - m_{\beta'}.
\]
So far we have only defined $\lambda'$ and $\mu'$ up to translation,
but there is a unique choice of $\lambda'$ and $\mu'$, up to
simultaneous translation of both, such that
\[
\lambda(t) - \mu(t) = \lambda'(t) - \mu'(t)
\]
for all $t$ not an edge.  We order the edges of $\alpha$ and $\beta$
such that the agreeing `$h$' edges are ordered first and in the same
order, and this determines an ordering of the `$h$' edges of $\alpha'$
and $\beta'$.  Note for future reference that
\begin{equation}
\label{eqn:indexLocal}
I(\alpha,\beta) = I(\alpha',\beta').
\end{equation}

\begin{enumerate}
\item[VII] (Locality)\qua
  Let $\alpha$ and $\beta$ satisfy the Nesting and Label Matching
  conditions, and let $\alpha'$ and $\beta'$ be obtained from $\alpha$
  and $\beta$ by removing matching edges as above.  Then
\[
\langle\partial\alpha,\beta\rangle =
\langle\partial\alpha',\beta'\rangle.
\]
\end{enumerate}

Suppose that $\beta$ is obtained from $\alpha$ by rounding a corner
and locally losing one `$h$', ie the differential coefficient
$\langle\delta\alpha,\beta\rangle=\pm 1$ as in \fullref{sec:mgpc}.  If
furthermore only one new edge is created by the rounding process, and
if this edge has multiplicity one (ie the corresponding segment in
$\Z^2$ contains no interior lattice points), then we say that $\beta$
is obtained from $\alpha$ by {\em simple rounding\/}.  If no edges at
all are created by the rounding process, ie if $\alpha$ turns by
angle $\pi$ at the rounded corner, then we say that $\beta$ is
obtained from $\alpha$ by {\em degenerate rounding\/}.

\begin{enumerate}
\item[VIII] (Simple Rounding)\qua If $\beta$ is obtained from $\alpha$ by
  simple rounding, then
$\langle\partial\alpha,\beta\rangle =
\langle\delta\alpha,\beta\rangle$.
\item[IX] (Degenerate Rounding)\qua If $\beta$ is obtained from $\alpha$
by degenerate rounding then $\langle\partial\alpha,\beta\rangle =
\langle\delta\alpha,\beta\rangle$.
\end{enumerate}

\subsection{Uniqueness of the chain complex}
\label{sec:uniquenessTheorem}

As usual let $\delta$ denote the differential on $\wwtilde{C}_*(2\pi
n)$ defined in \fullref{sec:mgpc}.  Of course, $(\wwtilde{C}_*(2\pi
n),\delta)$ satisfies the above axioms.

\begin{proposition}
\label{prop:uniqueness}
Let $\partial$ be a differential on $\wwtilde{C}_*(2\pi n)$
satisfying the axioms of \fullref{sec:axioms}.  Then $\partial=\delta$.
\end{proposition}

\begin{proof}
  Throughout this proof, $\alpha$ and $\beta$ will denote generators
  of $\wwtilde{C}_*(2\pi n)$, and $\lambda$ and $\mu$ will denote
  their underlying admissible paths.

\begin{lemma}
\label{lem:vanishing}
If $\langle\partial\alpha,\beta\rangle \neq 0$, then
$\langle\delta\alpha,\beta\rangle \neq 0$.
\end{lemma}

\begin{proof}
Suppose that $\langle\partial\alpha,\beta\rangle\neq 0$.  By the
Nesting axiom, $\beta\le\alpha$.  By
\fullref{prop:roundingSequence}, the polygon $\mu$ can be
obtained from $\lambda$ by a sequence of $k$ corner roundings for some
nonnegative integer $k$.  By \fullref{lem:pick},
\[
I(\alpha,\beta) = 2k -\#h(\alpha) + \#h(\beta).
\]
Let $D\eqdef D(\alpha,\beta)$ be defined as in \fullref{sec:axioms}.
Let $l$ denote the number of edges of $\alpha$ that are in $D$, ie
that do not agree with any edges in $\beta$.  Since $D$ is connected
by the Connectedness axiom, these $l$ edges of $\alpha$ are
consecutive.  Observe that
\[
l\le k+1.
\]
Otherwise, at least one of the corners between the edges of $\alpha$
in $D$ is not rounded in a sequence of $k$ roundings
from $\mu$ to $\lambda$.  If $D\neq \R/2\pi n\Z$, then
$D$ is separated by such a corner, contradicting the
Connectedness axiom.  If $D=\R/2\pi n\Z$, then two
corners in $\alpha$ are not rounded and these two corners separate
$D$.

By the Label Matching axiom, we can calculate $\#h(\alpha)-\#h(\beta)$
by considering only the edges of $\alpha$ and $\beta$ in
$D$, so
\[
\#h(\alpha)-\#h(\beta)\le l.
\]
Combining this with the previous inequality and equation gives
\[
I(\alpha,\beta) \ge k-1.
\]
By the Index axiom, $I(\alpha,\beta)=1$, so the only possibilities are
$k=0$, $k=1$, or $k=2$.

The case $k=0$ is impossible, as then the Label Matching axiom would
imply that $\alpha=\beta$ so that $I(\alpha,\beta)=0$.

If $k=2$ then equality must hold in the above inequalities so $l=3$.
But this is forbidden by the No Double Rounding axiom.

Therefore $k=1$, so $\mu$ is obtained from $\lambda$ by
rounding a corner.  By the index formula, $\#h(\alpha)-\#h(\beta)=1$.
Together with the Label Matching axiom, this implies that
$\langle\delta\alpha,\beta\rangle = \pm 1$.
\end{proof}

\begin{lemma}
\label{lem:backtrackingTrick}
If $\langle\delta\alpha,\beta\rangle\neq0$, then
\begin{equation}
\label{eqn:varepsilon}
\langle\partial\alpha,\beta\rangle =
\langle\delta\alpha,\beta\rangle.
\end{equation}
\end{lemma}

\begin{proof}
Suppose that $\langle\delta\alpha,\beta\rangle\neq 0$; we will show
that \eqref{eqn:varepsilon} holds.  The strategy is to use
$\partial^2=0$ to solve for the unknown differential coefficients.
The following is similar to an argument in our earlier paper
\cite[Section~3.8]{pfh3}, but because we are considering more general
polygonal paths here we can make some simplifications.

We know that $\beta$ is obtained from $\alpha$ by rounding a corner
$c$ and locally losing one `$h$'.  Let $\theta_1$ and $\theta_2$ be
the edges of $\alpha$ preceding and following $c$, respectively.  Let
\begin{equation}
\label{eqn:roundingDeterminant}
\Delta \eqdef
\det\begin{pmatrix}x_{\theta_1}&x_{\theta_2}\\
y_{\theta_1}&y_{\theta_2}\end{pmatrix} \in \Z.
\end{equation}
By the definition of rounding a corner, $\theta_2-\theta_1\in(0,\pi]$,
and in particular $\Delta\ge 0$.  We will now prove equation
\eqref{eqn:varepsilon} by induction on $\Delta$.

If $\Delta=0$, then \eqref{eqn:varepsilon} holds by the Degenerate
Rounding axiom.

If $\Delta=1$, then the triangle with vertices
$$\lambda(c)-\begin{pmatrix}x_{\theta_1}\\y_{\theta_1}\end{pmatrix},\quad
\lambda(c),\quad \lambda(c) +
\begin{pmatrix}x_{\theta_2}\\y_{\theta_2}\end{pmatrix}$$
is simple, so
\eqref{eqn:varepsilon} holds by the Simple Rounding axiom.

Now suppose that $\Delta>1$, and assume that the lemma holds for all
smaller values of $\Delta$.  Let $c'$ be the corner of $\alpha$
following the edge $\theta_2$.  Let $\theta_3$ be the edge of $\beta$
preceding $c'$, ie the last edge of $\beta$ created by rounding at
$c$.  By the Locality axiom, we may replace $\alpha$ and $\beta$ by
the generators $\alpha'$ and $\beta'$ in the statement of the Locality
axiom, and then add one new edge to each, to arrange the following:
\begin{itemize}
\item
The edges $\theta_1$ and $\theta_2$ of $\alpha$ have multiplicity $1$.
\item
$\alpha$ has only one edge other than $\theta_1$ and $\theta_2$, and
$\beta$ has only one edge other than the edges created by rounding the
corner $c$.  In both $\alpha$ and $\beta$, this additional edge
is at angle $\theta_3+\pi$ with multiplicity $1$.
\item
If the edge $\theta_3$ of $\beta$ is labeled `$h$', then the
edge $\theta_3+\pi$ of $\alpha$ and $\beta$ is labeled `$e$';
  otherwise the edge $\theta_3+\pi$ of $\alpha$ and $\beta$ is labeled
  `$h$'.
\end{itemize}

The last condition above ensures that it is possible to round the
corner $c'$ of $\beta$ and locally lose one `$h$' to obtain a
well-defined (up to sign) generator $\gamma$ with
$\langle\delta\beta,\gamma\rangle \neq 0$.  Since this rounding is
degenerate, we also know that
\[
\langle\partial\beta,\gamma\rangle =
\langle\delta\beta,\gamma\rangle \neq 0.
\]
An example of the admissible paths underlying $\alpha$ and $\gamma$ is
shown below.
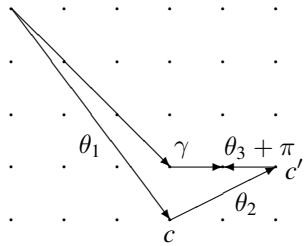
\begin{center}
\begin{picture}(117,90)(-8,-2)\small
\put(0,0){.}
\put(20,0){.}
\put(40,0){.}
\put(60,0){.}
\put(80,0){.}
\put(100,0){.}
\put(0,20){.}
\put(20,20){.}
\put(40,20){.}
\put(60,20){.}
\put(80,20){.}
\put(100,20){.}
\put(0,40){.}
\put(20,40){.}
\put(40,40){.}
\put(60,40){.}
\put(80,40){.}
\put(100,40){.}
\put(0,60){.}
\put(20,60){.}
\put(40,60){.}
\put(60,60){.}
\put(80,60){.}
\put(100,60){.}
\put(0,80){.}
\put(20,80){.}
\put(40,80){.}
\put(60,80){.}
\put(80,80){.}
\put(100,80){.}
\put(1.5,80.5){\vector(3,-4){60}}
\put(61.5,0.5){\vector(2,1){40}}
\put(101.5,20.5){\vector(-1,0){20}}
\put(1.5,80.5){\vector(1,-1){60}}
\put(61.5,20.5){\vector(1,0){20}}
\put(26,26){$\theta_1$}
\put(58,-8){$c$}
\put(85,4){$\theta_2$}
\put(104,16){$c'$}
\put(81,25){$\theta_3+\pi$}
\put(62,26){$\gamma$}
\end{picture}
\end{center}

In general the rounding of $\alpha$ at $c'$ is simple, because the
triangle with vertices
$$\lambda(c),\quad \lambda(c'),\quad
\lambda(c')-\begin{pmatrix}x_{\theta_3}\\y_{\theta_3}\end{pmatrix}$$
is simple, by the definition of rounding the corner of $\alpha$ at
$c$.  Thus there is a unique (up to sign) generator $\beta'$ obtained
from $\alpha$ by rounding the corner at $c'$ and locally losing one
`$h$'.  This generator $\beta'$ satisfies
$\langle\delta\alpha,\beta'\rangle\neq 0$, and because $\beta'$ is
obtained from $\alpha$ by simple rounding, we also know that
\[
\langle\partial\alpha,\beta'\rangle =
\langle\delta\alpha,\beta'\rangle \neq 0.
\]

Finally, the edge labels work out so that we can round $\beta'$ at $c$
and locally lose one `$h$' to obtain $\gamma$, and in particular
$\langle\delta\beta',\gamma\rangle \neq 0$.  Moreover the determinant
corresponding to this rounding as in \eqref{eqn:roundingDeterminant}
is less than $\Delta$, because the triangle with vertices
$$\lambda(c)-\begin{pmatrix} x_{\theta_1} \\ y_{\theta_1}
\end{pmatrix},\quad \lambda(c),\quad \lambda(c')-\begin{pmatrix}
x_{\theta_3} \\ y_{\theta_3} \end{pmatrix}$$
is a proper subset of the
triangle with vertices
$$\lambda(c) - \begin{pmatrix} x_{\theta_1} \\
y_{\theta_1} \end{pmatrix},\quad \lambda(c),\quad \lambda(c').$$
So by inductive hypothesis,
\[
\langle\partial\beta',\gamma\rangle =
\langle\delta\beta',\gamma\rangle \neq 0.
\]

By Lemmas~\ref{lem:deltaSquared}(a) and \ref{lem:vanishing}, there
does not exist a generator $\beta''$, other than $\pm\beta'$ and possibly
$\pm\beta$, with $ \langle\partial\alpha,\beta''\rangle,
\langle\partial\beta'',\gamma\rangle \neq 0$.  Using this
fact and then plugging in the previous three equations, we get
\[
\begin{split}
0 = \langle\partial^2\alpha,\gamma\rangle &=
\langle\partial\alpha, 
\beta\rangle\langle\partial\beta,
\gamma\rangle +
\langle\partial\alpha,\beta'\rangle
\langle\partial\beta',\gamma\rangle\\
&= 
\langle\partial\alpha, 
\beta\rangle\langle\delta\beta,
\gamma\rangle +
\langle\delta\alpha,\beta'\rangle
\langle\delta \beta',\gamma\rangle.
\end{split}
\]
By \fullref{lem:deltaSquared}(b),
\[
0 = \langle\delta^2\alpha,\gamma\rangle = \langle\delta\alpha, 
\beta\rangle\langle\delta\beta,
\gamma\rangle +
\langle\delta\alpha,\beta'\rangle
\langle\delta \beta',\gamma\rangle.
\]
Since all factors on the right hand side are nonzero, comparing this
equation with the previous one proves \eqref{eqn:varepsilon}.
\end{proof}

The above two lemmas prove \fullref{prop:uniqueness}.
\end{proof}

\section{$J$--holomorphic curves in $\R\times T^3$}
\label{sec:interlude}

Having completed the proofs of our algebraic theorems, we now gather
some fundamental facts about $J$--holomorphic curves in $\R\times T^3$, in
preparation for computing the embedded contact homology of $T^3$.
\fullref{sec:generalSetup} gives basic definitions.  In \fullref{sec:t3}
we establish a dictionary between some of the combinatorics of
\fullref{sec:preliminaries} and the geometry of $J$--holomorphic curves
in $\R\times T^3$.  In \fullref{sec:nesting} we prove a useful
restriction on the latter in terms of the partial order from
\fullref{sec:partialOrder}.  In \fullref{sec:taubes} and
\fullref{sec:zeroArea} we recall and prove some basic classification
results for $J$--holomorphic curves in $\R\times T^3$.

\subsection{$J$--holomorphic curves in symplectizations}
\label{sec:generalSetup}

Let $Y$ be a closed oriented 3--manifold with a contact form $\lambda$; see
\fullref{sec:motivation} for the basic contact terminology.

\begin{definition}
\label{def:orbitSet}
An {\em orbit set\/} is a finite set of pairs
$\alpha=\{(\alpha_i,m_i)\}$ where the $\alpha_i$'s are distinct embedded Reeb
orbits and the $m_i$'s are positive integers (``multiplicities'').
The homology class of $\alpha$ is defined by
\[
[\alpha] \eqdef \sum_im_i[\alpha_i] \in H_1(Y).
\]
\end{definition}

\begin{definition}
\label{def:affine}
If $\alpha=\{(\alpha_i,m_i)\}$ and $\beta=\{(\beta_j,n_j)\}$ are orbit
sets with $[\alpha]=[\beta]$, let $H_2(Y,\alpha,\beta)$ denote the set
of relative homology classes of $2$--chains $Z$ in $Y$ with
\[
\partial Z = \sum_im_i\alpha_i - \sum_jn_j\beta_j.
\]
Thus $H_2(Y,\alpha,\beta)$ is an affine space modelled on $H_2(Y)$.
\end{definition}

\begin{definition}
  An almost complex structure $J$ on $\R\times Y$ is {\em
    admissible\/} if $J$ is $\R$--invariant; $J(\xi)=\xi$ with
    $d\lambda(v,Jv)> 0$ for nonzero $v\in\xi$; and $J(\partial_s)$ is
    a positive multiple of $R$, where $s$ denotes the $\R$ coordinate.
\end{definition}

For our purposes, a {\em $J$--holomorphic curve\/} is a nonconstant map
$u\co C\to\R\times Y$, modulo reparametrization, where $C$ is a punctured
compact (possibly disconnected) Riemann surface with a complex
structure $j$, such that $du\circ j = J\circ du$.  When $u$ is an
embedding, we often identify $u$ with its image in $\R\times Y$.

For admissible $J$, if $\gamma\subset Y$ is an embedded Reeb orbit,
then $\R\times \gamma\subset\R\times Y$ is a $J$--holomorphic cylinder,
which we call a {\em trivial cylinder\/}.  Given a more general
$J$--holomorphic curve $u\co C\to\R\times Y$, a {\em positive end at
$\gamma$ of multiplicity $k$\/} is an end of $u$ asymptotic to
$\R\times\gamma^k$ as $s\to+\infty$, where $\gamma^k$ denotes the
$k$--fold connected covering of $\gamma$.  A {\em negative end\/} is
defined analogously with $s\to-\infty$.

\begin{definition}
\label{def:moduli}
If $\alpha=\{(\alpha_i,m_i)\}$ and $\beta=\{(\beta_j,n_j)\}$ are orbit
sets with $[\alpha]=[\beta]$, let $\mc{M}^J(\alpha,\beta)$ denote the
moduli space of $J$--holomorphic curves $u\co C\to\R\times Y$ as above
such that:
\begin{itemize}
\item
$u$ has positive ends at $\alpha_i$, whose multiplicities sum to
$m_i$.
\item
Similarly $u$ has negative ends at $\beta_j$ of total
multiplicity $n_j$.
\item
$u$ has no other ends.
\end{itemize}
\end{definition}

Note that $\R$ acts on $\mc{M}^J(\alpha,\beta)$ by translation in the
$\R$ direction in $\R\times Y$.  If $u\in\mc{M}^J(\alpha,\beta)$, then
the projection of $u$ from $\R\times Y$ to $Y$ has a well-defined
relative homology class
\[
[u] \in H_2(Y,\alpha,\beta).
\]

\begin{definition}
\label{def:MABZ}
If $Z\in H_2(Y,\alpha,\beta)$, let
\[
\mc{M}^J(\alpha,\beta;Z) \eqdef \left\{u\in \mc{M}^J(\alpha,\beta)
  \bigmid [u]=Z \right\}.
\]
\end{definition}

\begin{definition}
\label{def:action}
If $\alpha=\{(\alpha_i,m_i)\}$ is an orbit set, define the {\em
symplectic action\/}
\[
\mc{A}(\alpha) \eqdef \sum_i m_i \int_{\alpha_i}\lambda.
\]
\end{definition}

\begin{lemma}
\label{lem:action}
For an admissible almost complex structure $J$, if
$\mc{M}^J(\alpha,\beta)$ is non\-empty, then:
\begin{enumerate}
\item[(a)]
$\mc{A}(\alpha) \ge
\mc{A}(\beta)$.
\item[(b)]
If $\mc{A}(\alpha)=\mc{A}(\beta)$, then $\alpha=\beta$ and every element of
$\mc{M}^J(\alpha,\beta)$ maps to a union of trivial cylinders.
\end{enumerate}
\end{lemma}

\begin{proof}
  Suppose $u\in\mc{M}^J(\alpha,\beta)$.  Admissibility of $J$ implies
  that if $v$ is a tangent vector to a point in the domain $(C,j)$,
  then $u^*d\lambda(v,jv)\ge 0$.  Part (a) follows immediately from
  Stokes theorem.  Part (b) holds because $u^*d\lambda(v,jv)=0$ only
  if $du$ sends $v$ to the span of $\partial_s$ and $R$ in $T(\R\times
  Y)$.
\end{proof}

\subsection{Admissible paths and orbit sets in $T^3$}
\label{sec:t3}

Fix a positive integer $n$.  We now specialize to the example $Y=T^3$
with the contact form $\lambda_n$ defined by \eqref{eqn:T^3} and
\eqref{eqn:lambda_n}.  The Reeb orbits of $\lambda_n$ consist of
circles of Reeb orbits at each $\theta\in\Theta_n$, where
\[
\Theta_n \eqdef \{\theta\in\R/2\pi n\Z \mid \tan\theta\in\Q\cup\{\infty\}\}.
\]
Each Reeb orbit $\gamma$ in the circle at $\theta$ has homology class
\[
[\gamma] = (0,(x_\theta,y_\theta)) \in H_1(T^3).
\]
In this setting we define a {\em Morse--Bott orbit set\/} to be a
finite set of pairs $\alpha=\{(\alpha_i,m_i)\}$ where each $\alpha_i$
is a component of the space of embedded Reeb orbits and each $m_i$ is
a positive integer.  A Morse--Bott orbit set $\alpha$ with
$[\alpha]=\Gamma\in \Z^2 = H_1(T^2) \subset H_1(T^3)$ is equivalent to
a multiplicity function
\[
m\co \Theta_n\longrightarrow\Z_{\ge 0}
\]
which is finitely supported and which satisfies
\begin{equation}
\label{eqn:m}
\sum_{\theta\in \Theta_n} m(\theta) \begin{pmatrix} x_\theta \\
y_\theta \end{pmatrix} = \Gamma.
\end{equation}
So by the discussion in \fullref{sec:admissiblePaths}, there is a
canonical bijection
\begin{equation}
\label{eqn:OSAP}
\left\{\begin{array}{c}\mbox{Morse--Bott orbit}\\
\mbox{sets $\alpha$ with $[\alpha]=\Gamma$}\end{array}\right\} =
\left\{\begin{array}{c}\mbox{periodic admissible paths of}\\
\mbox{rotation number $n$ and period
$\Gamma$}\end{array}\right\}/\mbox{translation}.
\end{equation}

\begin{remark}
\label{remark:length=action}
Under this correspondence, the {\em length\/} of a periodic admissible
path, defined in \fullref{sec:partialOrder}, agrees with the symplectic
action of the corresponding Morse--Bott orbit set as in
\fullref{def:action}.
\end{remark}

There seems to be no natural way to resolve the translation ambiguity
in \eqref{eqn:OSAP} for a single path.  However, the relative
translation ambiguity of a pair of paths does have a geometric
interpretation, as we now explain.

\begin{definition}
  Let $m$ and $m'$ be finitely supported functions $\Theta_n\to\Z_{\ge
  0}$ satisfying \eqref{eqn:m}.  A {\em relative placement\/} of $m$
  and $m'$ is a locally constant map
\[
f\co (\R/2\pi n\Z)
\setminus (\op{supp}(m) \cup \op{supp}(m'))
\longrightarrow\Z^2
\]
satisfying the ``jumping condition''
\[
\frac{d f(t)}{dt} = \left(m(t)-m'(t)\right)
\sum_{\theta\in \Theta_n}\begin{pmatrix} x_\theta \\
y_\theta\end{pmatrix} \delta_\theta(t).
\]
Let $\mc{R}(m,m')$ denote the set of all such $f$;  this is an affine
space over $\Z^2$.
\end{definition}

The significance of this definition is that if $\Lambda$ and
$\Lambda'$ are periodic admissible paths of rotation number $n$ and
period $\Gamma$ with multiplicity functions $m$ and $m'$ respectively,
then
\[
f=\Lambda-\Lambda'
\]
is a relative placement of $m$ and $m'$.

On the geometric side, if $\alpha$ and $\alpha'$ are Morse--Bott
orbit sets with $[\alpha]=[\alpha']$, then $H_2(T^3,\alpha,\alpha')$
is a well-defined affine space over $H_2(T^3)/H_2(T^2)=\Z^2$.  We have
to mod out by $H_2(T^2)$ because each circle of Reeb orbits sweeps out
a surface $\{\theta\}\times T^2$ in $S^1\times T^2$.

\begin{lemmadef}
\label{lem:RHRP}
Let $\alpha$ and $\alpha'$ be Morse--Bott orbit sets corresponding to
multiplicity functions $m$ and $m'$.  Then there is a canonical
$\Z^2$--equivariant bijection between relative homology classes in $T^3$
and relative placements of periodic admissible paths,
\[
H_2(T^3,\alpha,\alpha') = \mc{R}(m,m').
\]
\end{lemmadef}

\begin{proof}
Let $T$ and $T'$ denote the supports of $m$ and $m'$.  For $Z\in
H_2(T^3,\alpha,\alpha')$ and $\theta_0\in(\R/2\pi n\Z)\setminus (T\cup
T')$, the intersection of $Z$ with the slice
$\{\theta=\theta_0\}\subset T^3$ has a well-defined homology class as
follows.  If $Z$ is represented by a smooth surface $\Sigma$
intersecting $\{\theta=\theta_0\}$ transversely, then the intersection
is a compact 1--manifold.  We orient the intersection so that if
$\{v,w\}$ is an oriented basis for the tangent space to $\Sigma$ at a
point and $v$ is a positively oriented tangent vector to the
intersection, then $w$ has positive $\partial_\theta$ component.  Now
define
\[
f(\theta_0)\eqdef\big[ Z\cap\{\theta=\theta_0\}\big]\in H_1(T^2)=\Z^2.
\]
Then $f\in\mc{R}(m,m')$, and this defines the required bijection.
\end{proof}

For an admissible almost complex structure $J$ as in
\fullref{sec:generalSetup} and $Z\in H_2(T^3,\alpha,\alpha')$, we can
define $\mc{M}^J(\alpha,\alpha',Z)$ by analogy with
\fullref{def:MABZ}.

\begin{definition}
If $\Lambda$ and $\Lambda'$ are periodic admissible paths of rotation
number $n$ and period $\Gamma$ corresponding to the Morse--Bott orbit
sets $\alpha$ and $\alpha'$, let
\[
\mc{M}^J(\Lambda,\Lambda') \eqdef
\mc{M}^J(\alpha,\alpha',\Lambda-\Lambda').
\]
\end{definition}

\subsection{Nesting of polygons and intersection positivity}
\label{sec:nesting}

There is a simple but important constraint on $J$--holomorphic curves
in $\R\times T^3$ in terms of the partial order $\le$ from
\fullref{sec:partialOrder}.

\begin{proposition}
\label{prop:nesting}
Let $J$ be an admissible almost complex structure on $\R\times T^3$
for the contact form $\lambda_n$.
Let $\Lambda$ and $\Lambda'$ be periodic admissible paths of rotation
number $n$ and period $\Gamma$.  Then
\[
\mc{M}^{J}(\Lambda,\Lambda')\neq\emptyset \Longrightarrow
\Lambda'\le\Lambda.
\]
\end{proposition}

\begin{proof}
Let $u\in \mc{M}^{J}(\Lambda,\Lambda')$.  We want to
show that for all $\theta\in\R$,
\begin{equation}
\label{eqn:inside}
\det\left(\begin{matrix}\begin{matrix}
\cos \theta \\
\sin \theta
\end{matrix}
&
\Lambda'(\theta)-\Lambda(\theta)
\end{matrix}
\right) \ge 0.
\end{equation}
Choose $\theta_0\in(\R/2\pi n\Z)\setminus (T\cup T')$ such that $u$ is
transverse to $\{\theta=\theta_0\}\subset\R\times T^3$.  Consider a
component of $u^{-1}\{\theta=\theta_0\}$, parametrized in an
orientation-preserving manner by
\[
\rho=(s,x,y)\co S^1\longrightarrow\{\theta=\theta_0\}.
\]
Let $\tau$ denote the $S^1$ coordinate.  By admissibility of $J$, the
transversality assumption, and our orientation convention for the
intersection, we have
\[
\sin(\theta_0)\frac{dx}{d\tau} - \cos(\theta_0)\frac{dy}{d\tau} > 0.
\]
(This can also be understood as positivity of intersections of $u$
with the leaves of the $J$--holomorphic foliation of
$\{\theta=\theta_0\}$ by $\R$ times the Reeb flow.)  Hence the
homology class $(\rho_x,\rho_y)\in\Z^2$ of this component satisfies
\[
\det
\begin{pmatrix}\rho_x & \cos \theta_0 \\ \rho_y & \sin
\theta_0\end{pmatrix}
> 0.
\]
Adding this up for all components of $u^{-1}\{\theta=\theta_0\}$
proves equation \eqref{eqn:inside} whenever $\theta\in\R$ is a lift of
$\theta_0\in(R/2\pi n\Z)\setminus (T\cup T')$.  By continuity,
equation \eqref{eqn:inside} holds for all $\theta\in\R$.
\end{proof}

\begin{remark}
\label{remark:emptySlice}
For a given $\theta_0 \notin T\cup T'$, the above argument shows that
if equality holds in \eqref{eqn:inside}, then $u$ does not intersect
the slice $\{\theta=\theta_0\}$.  \end{remark}

\subsection{Spheres with two or three punctures, and degenerate and simple
rounding}
\label{sec:taubes}

Of particular interest is the ``standard'' almost complex structure
$J_{\rm std}$ on $\R\times T^3$ defined by
\begin{equation}
\label{eqn:Jstd}
\begin{split}
  J_{\rm std}(\partial_s) &\eqdef \cos\theta\,\partial_x +
  \sin\theta\,\partial_y,\\
  J_{\rm std}(\partial_\theta) &\eqdef -\sin\theta\,\partial_x +
  \cos\theta\,\partial_y.
\end{split}
\end{equation}
It is easy to check that $J_{\rm std}$ is admissible.  Since
$J_{\rm std}$ does not depend on $s$, $x$, or $y$, the action of $\R\times
T^2$ on $\R\times T^3$ preserves $J_{\rm std}$ and hence induces an action
on $\mc{M}^{J_{\rm std}}(\Lambda,\Lambda')$.

\begin{definition}
  If $\Lambda,\Lambda'$ are periodic admissible paths with rotation
  number $n$, let $\mc{M}_0^{J}(\Lambda,\Lambda')$ denote the set of
  irreducible, genus zero curves $u\in\mc{M}^{J}(\Lambda,\Lambda')$.
\end{definition}

\begin{proposition}
\label{prop:taubes}
Let $J$ be any $T^2$--invariant admissible almost complex structure on
$\R\times T^3$ for the standard contact form $\lambda_n$.  Suppose
$\Lambda$ is a periodic admissible path of rotation number $n$, with
two edges.  Then:
\begin{enumerate}
\item[(a)] If $\Lambda'$ is obtained from $\Lambda$ by degenerate
  rounding (see \fullref{sec:axioms}), then $\R\times T^2$ acts
  transitively on $\mc{M}_0^{J}(\Lambda,\Lambda')$ with $S^1$
  stabilizer.
\item[(b)]
If $\Lambda'$ is obtained from $\Lambda$ by simple rounding, then
$\R\times T^2$ acts freely and transitively on
$\mc{M}_0^{J}(\Lambda,\Lambda')$.
\end{enumerate}
\end{proposition}

\begin{proof}
We will deduce the proposition from analogous results of Taubes
\cite{taubes02}, which hold for a similar contact form $\lambda_T$ on
$S^1\times S^2$ and a $T^2$--invariant admissible almost complex
structure $J_T$ on $\R\times S^1\times S^2$.  For this purpose we will
need to consider slightly more general contact forms on $T^3$.
Namely, consider
\[
\lambda = a_1(\theta)\,dx + a_2(\theta)\,dy
\]
where
\[
a = (a_1,a_2) \co  \R/2\pi n\Z \longrightarrow \R^2\setminus\{0\}
\]
has properties (i)--(iii) below.  (In the following, if $v=(v_1,v_2)$
and $w=(w_1,w_2)$ are vectors in $\R^2$, then $v\times w \eqdef
v_1w_2-v_2w_1$.)
\begin{enumerate}
\item[(i)]
The path $a$ has winding number $n$ around the origin in $\R^2$.
\item[(ii)]
$a\times a' >0$ for all $\theta$.
\item[(iii)]
$a'\times a''>0$ for all $\theta$.
\end{enumerate}
Condition (ii) ensures that $\lambda$ is a contact form.
The Reeb vector field is given by
\[
R = \frac{a_2'\,\partial_x - a_1'\,\partial_y}{a\times a'}.
\]
By (i) and (ii), $R$ has winding number $n$.  Condition (iii) implies
that $R$ turns to the left as $\theta$ increases.  Hence we can
reparametrize the $\theta$ coordinate (in exactly $n$ different ways)
so that
\begin{enumerate}
\item[(iv)] $R$ is a positive multiple of $\cos\theta\,\partial_x +
  \sin\theta\,\partial_y$ for all $\theta$.
\end{enumerate}
Of course the standard contact form $\lambda_n$ is recovered by taking
$a=(\cos\theta,\sin\theta)$, which satisfies conditions (i)--(iv)
above.  Also, any two contact forms satisfying (i)--(iv) above are
homotopic through such forms (by linear interpolation).  For such a
contact form, all of \fullref{sec:t3} and \fullref{sec:nesting} holds
verbatim (except for \fullref{remark:length=action}).  We will
prove the proposition for any such contact form.  We proceed in three
steps.

\textbf{Step 1}\qua Denote the edges of $\Lambda$ by $\theta_1$ and
$\theta_2$ with $\theta_2-\theta_1\in(0,\pi]$.  Then the subset
\[
[\theta_1,\theta_2]\times T^2\subset T^3
\]
can be identified with a subset of $S^1\times S^2$ between two
latitude lines, such that the pullback of Taubes's contact form
$\lambda_T$ extends to a contact form $\lambda_T'$ satisfying
(i)--(iv) above.  The pullback of Taubes's almost complex structure
$J_T$ extends to a $T^2$--invariant admissible $J_T'$ on $\R\times
T^3$.  By \fullref{remark:emptySlice}, any
$u\in\mc{M}^{J_T'}(\Lambda,\Lambda')$ maps to
$\R\times[\theta_1,\theta_2]\times T^2$, and an analogous argument
works for $J_T$--holomorphic curves in $\R\times S^1\times S^2$.  Hence
Taubes's results are applicable to $\mc{M}^{J_T'}(\Lambda,\Lambda')$.
In particular, \cite[Theorem~A.1(c)]{taubes02} proves (a), and
\cite[Theorem~A.2]{taubes02} proves (b), for $\lambda_T'$ and $J_T'$.

\textbf{Step 2}\qua Now consider another contact form $\lambda$ on $T^3$
satisfying (i)--(iv) above, and an admissible $T^2$--invariant
admissible $J$ on $\R\times T^3$.  We can deform $\lambda_T'$ to
$\lambda$ through contact forms satisfying (i)--(iv), and for this
family of contact forms we can find a family of $T^2$--invariant
admissible almost complex structures interpolating between $J_T'$ and
$J$.  The moduli spaces $\mc{M}_0(\Lambda,\Lambda')$ are smooth
manifolds of the expected dimension throughout the deformation, as in
\cite[Theorem~1.2]{taubes04}. (For more general automatic transversality
results see Wendl's doctoral thesis \cite[Section~4.5.5]{wendl}.)  By
Gromov compactness (see the paper \cite{behwz} by Bourgeois, Eliashberg,
Hofer, Wysocki and Zehnder) the moduli spaces
$\mc{M}_0(\Lambda,\Lambda')/\R$ are
compact throughout the deformation, because by
Propositions~\ref{prop:roundingMaximal}(b) and \ref{prop:nesting} and
\fullref{lem:action}(b), there are never any broken
pseudoholomorphic curves from $\Lambda$ to $\Lambda'$.  So the moduli
spaces $\mc{M}_0(\Lambda,\Lambda')/\R$ for $J_T'$ and for $J$ are
diffeomorphic.

\textbf{Step 3}\qua Consideration of the Reeb orbits that appear at the
ends of the $J$--holomorphic curves shows that in case (b), $\R\times
T^2$ acts freely on $\mc{M}_0^J(\Lambda,\Lambda')$.  This action must
then be transitive, or else $\mc{M}_0^J(\Lambda,\Lambda')$ would be
disconnected or not of the expected dimension, contradicting Step 2.
In case (a), $\R\times T^2/S^1$ acts freely on
$\mc{M}_0^J(\Lambda,\Lambda')$, where $S^1\subset T^2$ is generated by
the vector $(x_{\theta_1},y_{\theta_1})$.  So $\R\times T^2$ must act
transtively with $S^1$ stabilizer, or else again there would be a
contradiction of Step 2.
\end{proof}

The above proposition can also be deduced from work of Parker
\cite{parker}, which classifies genus zero pseudoholomorphic curves in
$\R\times T^3$ for a degeneration of $J_{\op{\rm std}}$, in terms of
certain labeled graphs in $\R^2\setminus\{0\}$.

\subsection{The zero area constraint}
\label{sec:zeroArea}

We now show that the sets of Reeb orbits that can appear at the ends
of a $J_{\rm std}$--holomorphic curve in $\R\times T^3$ satisfy a
codimension one constraint.  This will be used in \fullref{sec:echt3} to
establish the No Double Rounding axiom for the embedded contact
homology of $T^3$.  To state the constraint, for $\theta\in\Theta_n$
we explicitly identify the circle of Reeb orbits in $\{\theta\}\times
T^2$ with $S^1=\R/\Z$ via a map $\varphi$ defined as follows.  If the
Reeb orbit $\gamma$ contains a point $(\theta,x,y)\in T^3$, then we
define
\begin{equation}
\label{eqn:varphi}
\varphi(\gamma) \eqdef x_\theta y - y_\theta x + \frac{x_\theta
  y_\theta}{2} \in \R/\Z.
\end{equation}
If $u\in\mc{M}^J(\alpha,\beta)$, let $E_+(u)$
and $E_-(u)$ denote the set of positive and negative ends of $u$,
respectively.  For $e\in E_\pm(u)$, let
$\gamma(e)$ denote the corresponding embedded Reeb orbit and $m(e)$
the multiplicity of the end as defined in \fullref{sec:generalSetup}.

\begin{proposition}
\label{prop:conservedQuantity}
Let $u\in\mc{M}^{J_{\rm std}}(\Lambda,\Lambda')$.  Then
\begin{equation}
\label{eqn:area}
\sum_{e\in E_+(u)}m(e)\varphi(\gamma(e)) - \sum_{e\in
  E_-(u)}m(e)\varphi(\gamma(e)) = 0 \in \R/\Z.
\end{equation}
\end{proposition}

\begin{proof}
(Conpare our earlier paper \cite[Lemma~A.2]{pfh3}) It follows from \eqref{eqn:Jstd} that
the 2--form $-ds\,d\theta + dx\,dy$ on $\R\times T^3$ annihilates any
pair of tangent vectors of the form $(v,J_{\rm std}v)$.  Therefore
\[
\int_C u^*(dx\,dy) = \int_C u^*(ds\,d\theta).
\]
Now $\int_C u^*(ds\,d\theta)=0$ by Stokes' theorem, because the 1--form
$s\,d\theta$ vanishes along the Reeb orbits.  Therefore $\int_C
u^*(dx\,dy)=0$, ie the projection of $u$ to the $(x,y)$--torus has area
zero.  It follows from the identification \eqref{eqn:varphi} that this
area is congruent modulo $\Z$ to the left side of equation
\eqref{eqn:area}.
\end{proof}

\section{Embedded contact homology}
\label{sec:ECH}

We now (in \fullref{sec:defineECH}) outline the definition of the
embedded contact homology of a contact 3-manifold.  The idea is to
count $J$--holomorphic curves with $I=1$, where $I$ is a certain upper
bound on the index introduced in \fullref{sec:II}.  In \fullref{sec:echt3}
we explain the correspondence between the embedded contact homology of
$T^3$ and our combinatorial chain complexes.  This will prove
\fullref{thm:echt3}.

\subsection{The index inequality}
\label{sec:II}

As in \fullref{sec:generalSetup}, let $Y$ be a closed oriented
3--manifold, let $\lambda$ be a contact 1--form on $Y$, and let $J$ be
an admissible almost complex structure on $\R\times Y$.

If $\gamma$ is a Reeb orbit passing through a point $y\in Y$, then the
linearization of the Reeb flow $R$ on the contact planes along
$\gamma$ determines a linearized return map $P_\gamma:\xi_y\to \xi_y$.
This is a symplectic linear map whose eigenvalues do not depend on
$y$.  The Reeb orbit $\gamma$ is {\em nondegenerate\/} if $P_\gamma$
does not have $1$ as an eigenvalue.  Assume now that all Reeb orbits,
including multiply covered ones, are nondegenerate.

A Reeb orbit $\gamma$ is called {\em elliptic\/} or positive (resp.\
negative) {\em hyperbolic\/} when the eigenvalues of $P_\gamma$ are on
the unit circle or the positive (resp.\ negative) real line
respectively.  If $\tau$ is a trivialization of $\xi$ over $\gamma$,
one can then define the {\em Conley--Zehnder index\/}
$\mu_\tau(\gamma)\in\Z$.  In our three-dimensional situation this is
given explicitly as follows.  For a positive integer $k$, let
$\gamma^k$ denote the $k^{th}$ iterate of $\gamma$.  If $\gamma$ is
elliptic, then there is an irrational number $\phi\in\R$ such that
$P_\gamma$ is conjugate in $\op{SL}(2,\R)$ to a rotation by angle
$2\pi\phi$, and
\begin{equation}
\label{eqn:CZEll}
\mu_\tau(\gamma^k) = 2\lfloor k\phi\rfloor +1.
\end{equation}
Here $2\pi \phi$ is the total rotation angle with respect to $\tau$ of
the linearized flow around the orbit.  If $\gamma$ is positive (resp.\
negative) hyperbolic, then there is an even (resp.\ odd) integer $r$
such that the linearized flow around the orbit rotates the eigenspaces
of $P_\gamma$ by angle $\pi r$ with respect to $\tau$, and
\begin{equation}
\label{eqn:CZHyp}
\mu_\tau(\gamma^k) = kr.
\end{equation}

Let $\alpha=\{(\alpha_i,m_i)\}$ and $\beta=\{(\beta_j,n_j)\}$ be orbit
sets as in \fullref{sec:generalSetup}.  Suppose that $[\alpha]=[\beta]$
and let $Z\in H_2(Y,\alpha,\beta)$.

\begin{definition}[Compare Eliashberg--Givental--Hofer
\cite{eliashberg-givental-hofer00}]
\label{def:ind}
If
$u\in\mc{M}^J(\alpha,\beta;Z)$, define the {\em SFT index\/}
\begin{multline}
\label{eqn:ind}
\op{ind}(u) \eqdef \\ -\chi(C) + 2c_1(u^*\xi,\tau) + \sum_{e\in
E_+(u)}\mu_\tau\bigl(\gamma(e)^{m(e)}\bigr) - \sum_{e\in
E_-(u)}\mu_\tau\bigl(\gamma(e)^{m(e)}\bigr).
\end{multline}
Here $\tau$ is a trivialization of the 2--plane bundle $\xi$ over the
$\alpha_i$'s and $\beta_j$'s, and $c_1$ denotes the relative first
Chern class with respect to $\tau$, see
our earlier paper \cite[Section~2]{pfh2}.
\end{definition}

The following proposition is the $3$--dimensional case of a formula from
\cite{eliashberg-givental-hofer00} which is proved in the paper by Dragnev
\cite{dragnev}, using an index calculation by Schwarz \cite{schwarz}.

\begin{proposition}
\label{prop:dimension}
If $J$ is generic, and if $u\in\mc{M}^J(\alpha,\beta)$ has no multiply
covered components, then $\mc{M}^J(\alpha,\beta)$ is a manifold near
$u$ of dimension $\op{ind}(u)$. \qed
\end{proposition}

\begin{definition}[Hutchings \cite{pfh2}]
Define the {\em ECH index\/}
\begin{equation}
\label{eqn:ECHIndex}
I(\alpha,\beta,Z) \eqdef c_1(\xi|_Z,\tau) + Q_\tau(Z) + 
\sum_i\sum_{k=1}^{m_i}\mu_\tau\left(\alpha_i^k\right) -
\sum_j\sum_{k=1}^{n_j}\mu_\tau\left(\beta_j^k\right).
\end{equation}
Here $Q_\tau$ denotes the ``relative intersection pairing'', which is
defined in
\cite[Section~2]{pfh2}.  If $u\in\mc{M}^J(\alpha,\beta,Z)$, write
$I(u)\eqdef I(\alpha,\beta,Z)$.
\end{definition}

The following basic properties of $I$ are proved in \cite{pfh2}.
First, $I(\alpha,\beta,Z)$ does not depend on the choice of $\tau$.
Second, $I$ is additive in the sense that
\begin{equation}
\label{eqn:IAdditive}
I(\alpha,\beta,Z) + I(\beta,\gamma,W) = I(\alpha,\gamma,Z+W).
\end{equation}
Third, $I$ depends on $Z$ via the ``index ambiguity formula''
\begin{equation}
\label{eqn:IAmbiguity}
I(\alpha,\beta,Z) - I(\alpha,\beta,W) = \left\langle c_1(\xi)+
2\op{PD}(\Gamma), Z-W\right\rangle.
\end{equation}
Fourth, the index mod 2 is given by
\[
I(\alpha,\beta,Z) \equiv \#h(\alpha) - \#h(\beta) \mod 2
\]
where $\#h(\alpha)$ denotes the number of positive hyperbolic Reeb
orbits in $\alpha$.

The key, nontrivial property of $I$ is the inequality \eqref{eqn:II}
below which bounds the SFT index in terms of the ECH index.

\begin{proposition}
Suppose that $u\in\mc{M}^J(\alpha,\beta)$ does not multiply cover any
component of its image and that the image of $u$ contains no trivial
cylinders.  Then
\begin{equation}
\label{eqn:II}
\op{ind}(u) \le I(u) - 2\delta(u).
\end{equation}
Moreover, if $T$ is a union of (possibly multiply covered) trivial
cylinders, then
\begin{equation}
\label{eqn:TC}
I(u) \le I(u\cup T) - 2\#(u\cap T).
\end{equation}
\end{proposition}
Here $\delta(u)$ is a count of the singularities of $u$ with positive
integer weights; in particular $\delta(u)=0$ iff $u$ is an embedding.
Also, `$\#$' denotes the algebraic intersection number in $\R\times
Y$.  By intersection positivity (see McDuff \cite{mcduff}), $\#(u\cap T)\ge 0$,
with equality iff $u\cap T=\emptyset$.

\begin{proof}
Equation \eqref{eqn:II} follows from \cite[Equation~(18) and
Proposition~6.1]{pfh2}, and equation \eqref{eqn:TC} holds as in
\cite[Proposition~7.1]{pfh2}.  Note that these results in \cite{pfh2} are proved in a
slightly different setting, where $Y$ is a mapping torus and a ``local
linearity'' assumption is made.  The asymptotic analysis needed to
carry over these results to the present setting is done in
Siefring's doctoral thesis \cite{siefring}.
\end{proof}

The above proposition leads to strong restrictions on curves of low
ECH index:

\begin{corollary}
\label{cor:lowIndex}
Suppose $J$ is generic and $u\in\mc{M}^J(\alpha,\beta)$.  Then:
\begin{enumerate}
\item[(a)]
$I(u)\ge 0$.
\item[(b)]
If $I(u)=0$, then the image of $u$ is a union of trivial cylinders.
\item[(c)]
If $I(u)=1$, then $u$ contains one embedded component $u_1$ with
$\op{ind}(u_1)=I(u_1)=1$.  All other components of $u$ map to trivial
cylinders that do not intersect $u_1$.
\end{enumerate}
\end{corollary}

\begin{proof}
The image of $u$ consists of a union of $k$ irreducible
non-multiply-covered $J$--holomorphic curves $u_i$, covered by $u$ with
multiplicity $d_i$.  Let $u'$ be the union of $d_i$ different
translates of $u_i$ in the $\R$ direction, for each $i$ such that
$u_i$ is not a trivial cylinder.  Let $T$ be the union of the
components of $u$ that map to trivial cylinders.  Note that if $u_i$
is a trivial cylinder then $\op{ind}(u_i)=0$.  So by equations
\eqref{eqn:II} and \eqref{eqn:TC},
\begin{equation}
\label{eqn:MC}
\sum_{i=1}^k d_i\op{ind}(u_i) = \op{ind}(u') \le I(u') - 2\delta(u')
\le I(u) - 2\#(u'\cap T) - 2\delta(u').
\end{equation}
Since $J$ is generic, each nontrivial $u_i$ has $\op{ind}(u_i)>0$ by
\fullref{prop:dimension}, since $\R$ acts nontrivially on the
moduli space containing $u_i$.  Also, $\delta(u')=0$ only if $u'$ is
embedded, which implies that all of the nontrivial $u_i$'s are
embedded.  We can now read off the conclusions (a), (b), and (c) from
the inequality \eqref{eqn:MC}.
\end{proof}

\subsection{The definition of embedded contact homology}
\label{sec:defineECH}

Continue to assume that all Reeb orbits are nondegenerate.

\subsubsection{The chain complex}

\begin{definition}
\label{def:AOS}
An orbit set $\{(\alpha_i,m_i)\}$ is {\em admissible\/} if $m_i=1$
whenever $\alpha_i$ is hyperbolic.
\end{definition}

\begin{definition}
\label{def:untwisted}
If $\Gamma\in H_1(Y)$, then $C_*(Y,\lambda;\Gamma)$ is the free
$\Z$--module generated by admissible orbit sets $\alpha$ such that
$[\alpha]=\Gamma$, and an ordering of the positive hyperbolic Reeb
orbits in $\alpha$ is chosen.  We declare that changing this ordering
multiplies the generator by the sign of the reordering permutation.
\end{definition}

Let $N$ denote the divisibility of the image of $c_1(\xi) +
2\op{PD}(\Gamma)$ in $\Hom(H_2(Y),\Z)$.  It follows from
\eqref{eqn:IAmbiguity} and \eqref{eqn:IAdditive} that
$I(\alpha,\beta,Z)\mod N$ does not depend on $Z$ and defines a
relative $\Z/N$ grading on ${C}_*(Y,\lambda;\Gamma)$.

There is also a twisted chain complex defined for any subgroup
$G\subset H_2(Y)$.  Fix a ``reference cycle'', consisting of an
oriented $1$--dimensional submanifold $\rho\subset Y$ such that
\begin{equation}
\label{eqn:referenceCycle}
[\rho] = \Gamma \in H_1(Y).
\end{equation}

\begin{definition}
Let $\wwtilde{C}_*(Y,\lambda;\Gamma,G)$ be the free $\Z$--module
generated by pairs $(\alpha,[W])$ where $\alpha$ is generator of
$C_*(Y,\lambda;\Gamma)$ and $[W]\in H_2(Y,\rho,\alpha)/G$.
\end{definition}
The $H_2(Y)$ action on $H_2(Y,\rho,\alpha)$ makes
$\wwtilde{C}_*(Y,\lambda;\Gamma,G)$ into a free module over the
group ring $\Z[H_2(Y)/G]$, with one generator for each admissible
orbit set in the homology class $\Gamma$.  If $G=H_2(Y)$, then
$\wwtilde{C}_*$ reduces to the ``untwisted'' complex in
\fullref{def:untwisted}.  The ``fully twisted'' version has
$G=\{0\}$.

\subsubsection{The differential}

To define the differential, we first briefly review how to orient the
relevant moduli spaces of $J$--holomorphic curves following
Bourgeois--Mohnke \cite{bourgeois-mohnke00}.  For each Reeb orbit $\gamma$, there is a
determinant line $\mc{O}_\gamma$ associated to $\overline{\partial}$
operators on the plane with asymptotics determined by the linearized
Reeb flow along $\gamma$.  (When $\gamma$ is an even multiple cover of
a negative hyperbolic orbit, $\mc{O}_\gamma$ is only defined if one
also chooses a marked point on the image of $\gamma$.)  For each
$\gamma$ we choose an orientation of $\mc{O}_\gamma$.  When $\gamma$
is elliptic, there is a canonical ``complex'' orientation which we
choose, cf Floer--Hofer \cite[Theorem~2]{floer-hofer}.  By
\cite{bourgeois-mohnke00}, the above choices determine a sign for any
transversely cut out $\op{ind}=1$ curve provided that the ends at
positive hyperbolic orbits are ordered, and there are no ends at even
covers of negative hyperbolic orbits.

In the following, assume that admissible orbit sets have orderings of
the positive hyperbolic orbits chosen.  If $J$ is generic, and if
$\alpha$ and $\beta$ are homologous admissible orbit sets with
$I(\alpha,\beta,Z)=1$, define a count
\begin{equation}
\label{eqn:count}
\#\frac{\mc{M}^J(\alpha,\beta,Z)}{\R}\in\Z
\end{equation}
as follows.  Declare two curves $u,u'\in\mc{M}^J(\alpha,\beta,Z)/\R$
to be equivalent if their embedded components from
\fullref{cor:lowIndex}(c) are the same up to translation, and if
their other components cover each embedded trivial cylinder
$\R\times\gamma$ with the same total multiplicity.  In other words,
$u$ and $u'$ define the same current (modulo translation) in $\R\times
Y$.  The compactness argument of our earlier paper \cite[Section~9.4]{pfh2} shows that there
are only finitely many equivalence classes.  For each equivalence
class, if we discard the multiply covered trivial cylinders, then the
resulting embedded curve has a sign by the previous paragraph.  The
count \eqref{eqn:count} is now the sum over the equivalence classes of
the corresponding signs.

\begin{definition}
Define the differential
\[
\partial\co \wwtilde{C}_*(Y,\lambda;\Gamma,G)
\longrightarrow \wwtilde{C}_{*-1}(Y,\lambda;\Gamma,G)
\]
as follows.  If $\alpha$ is an admissible orbit set with
$[\alpha]=\Gamma$, then
\[
\partial(\alpha,[W]) \eqdef
\sum_{I(\alpha,\beta,Z)=1}\#\frac{\mc{M}^J(\alpha,\beta,Z)}{\R} \cdot
(\beta,[W+Z]).
\]
Here the sum is over admissible orbit sets $\beta$ and relative
homology classes $Z\in H_2(Y,\alpha,\beta)$.
\end{definition}

For technical reasons, we also need to consider, for a positive real
number $L$, the subcomplex $\wwtilde{C}_*^{< L}(Y,\lambda;\Gamma,G)$
generated by orbit sets $\alpha$ with symplectic action
$\mc{A}(\alpha) < L$.  By \fullref{lem:action}(a), the differential
$\partial$ sends $\wwtilde{C}_*^{<L}$ to itself.

A proof of the following is in
preparation.

\begin{conjecture}
\label{conj:ECHDefined}
\begin{enumerate}
\item[(a)] $\partial^2=0$.  
\item[(b)]
The homology of $\wwtilde{C}_*(Y,\lambda;\Gamma,G)$, which we denote by
$\widetilde{ECH}_*(Y,\lambda;\Gamma,G)$, does not depend on $J$.
\item[(c)]
The homology of $\wwtilde{C}_*^{<L}(Y,\lambda;\Gamma,G)$, which we denote by
$\widetilde{ECH}_*^{<L}(Y,\lambda;\Gamma,G)$, is invariant under
deformations of $\lambda$ during which all orbits of action $<L$ are
nondegenerate and no orbit has its action increase or decrease past
$L$.
\end{enumerate}
\end{conjecture}

\subsubsection{Morse--Bott version}

Suppose now that $\lambda$ has not only nondegenerate Reeb orbits but
also $S^1$--families of Reeb orbits which are nondegenerate in the
Morse--Bott sense.  In principle one could define ECH in this situation
along the lines of Bourgeois \cite{bourgeois04}, without perturbing $\lambda$.
However, the following definition is simpler to state.  Each
$S^1$--family of Reeb orbits, by a small perturbation of
$\lambda$, can be replaced by two embedded Reeb orbits, one elliptic
and one positive hyperbolic.  Since there are typically infinitely
many $S^1$--families of Reeb orbits, we cannot expect to perturb them
all this way simultaneously.  However, we can do this for all circles
of Reeb orbits of symplectic action $< L$.  So
\fullref{conj:ECHDefined} implies that
$\widetilde{ECH}_*^{<L}(Y,\lambda;\Gamma,G)$ is well-defined, and we
then define ECH as the direct limit
\[
\widetilde{ECH}_*(Y,\lambda;\Gamma,G) \eqdef \lim_{L\to\infty}
\widetilde{ECH}_*^{<L}(Y,\lambda;\Gamma,G).
\]
In both the nondegenerate and Morse--Bott cases, we denote the
``untwisted'' ECH by
\[
ECH_*(Y,\lambda;\Gamma)
\eqdef \widetilde{ECH}_*(Y,\lambda;\Gamma,H_2(Y)).
\]

\subsection{The example of $T^3$}
\label{sec:echt3}

We now explain why the untwisted embedded contact homology of $T^3$,
for the standard contact form $\lambda_n$ defined in \eqref{eqn:T^3}
and \eqref{eqn:lambda_n}, is computed by the combinatorial chain
complex $\wwbar{C}_*(2\pi n;\Gamma)$.  Also, the combinatorial
chain complex $\wwtilde{C}_*(2\pi n;\Gamma)$ computes a partially
twisted version of the embedded contact homology of $T^3$.  Below,
$H_*(T^2)$ denotes the homology of an $x,y$ torus in $T^3$.

\begin{theorem}
\label{thm:correspondence}
Assume \fullref{conj:ECHDefined}, so that ECH is
well-defined.  Then for $\Gamma\in H_1(T^2)=\Z^2$, the embedded
contact homology of $T^3$ is related to the combinatorial chain
complexes by
\begin{align}
\label{eqn:claimTilde}
\widetilde{ECH}_*(T^3,\lambda_n;\Gamma,H_2(T^2)) &\simeq
\wwtilde{H}_*(2\pi n;\Gamma),\\
\label{eqn:claimBar}
{ECH}_*(T^3,\lambda_n;\Gamma)
& \simeq
\wwbar{H}_*(2\pi n;\Gamma).
\end{align}
\end{theorem}

\begin{proof}
This is similar to the computation of the periodic Floer homology of a
Dehn twist on a cylinder in our earlier paper \cite{pfh3}, because the mapping torus
flow for a negative Dehn twist on a cylinder is isomorphic to the Reeb
flow on a subset of $T^3$ where $\theta$ ranges over an interval of
length less than $\pi$.  Thus we will carry over some lemmas from
\cite{pfh3}.  In making the translation, note that because the results
in \cite{pfh3} are stated for positive Dehn twists, positive ends of
$J$--holomorphic curves here correspond to negative (or ``incoming''')
ends there, and vice-versa.  We now prove the theorem in five steps.

\textbf{Step 1}\qua We begin by defining an isomorphism of relatively
graded $\Z[\Z^2]$--modules
\begin{equation}
\label{eqn:GCL}
\wwtilde{C}_*^{< L}(T^3,\lambda_n;\Gamma,H_2(T^2)) \simeq
\wwtilde{C}_*^{< L}(2\pi n;\Gamma).
\end{equation}
Here the right hand side denotes the subcomplex of
$C_*(2\pi n;\Gamma)$ generated by admissible paths of length $<L$, as
defined in \eqref{eqn:defineLength}.

Recall from \fullref{sec:t3} that for every $\theta\in\Theta_n$ there
is an $S^1$ family of embedded Reeb orbits, such that each Reeb orbit
$\gamma$ in the family has homology class
\[
[\gamma] = (0,(x_\theta,y_\theta)) \in H_1(T^3) = H_1(S^1) \oplus
H_1(T^2).
\]
After perturbation of $\lambda_n$, this family becomes an elliptic
orbit $e_\theta$ and a positive hyperbolic orbit $h_\theta$, of
approximately the same symplectic action.

To define \eqref{eqn:GCL}, we first define an isomorphism
of $\Z$--modules
\begin{equation}
\label{eqn:GCLBar}
{C}_*^{< L}(T^3,\lambda_n;\Gamma) =
\wwbar{C}_*^{< L}(2\pi n;\Gamma).
\end{equation}
Given a generator $\alpha$ of
${C}_*^{<L}(T^3,\lambda_n;\Gamma)$, define a multiplicity
function
\[
m\co  \R/2\pi n\Z \longrightarrow \Z_{\ge 0}
\]
by setting $m(\theta)$ equal to the total multiplicity of $e_\theta$
and $h_\theta$ in $\alpha$.  By \eqref{eqn:OSAP}, this defines an
admissible path $\Lambda$ of rotation number $n$ and period $\Gamma$,
up to translation.  By \fullref{remark:length=action}, the length
of $\Lambda$ is less than $L$.  If $\theta$ is an edge of $\Lambda$,
label it `$h$' if $h_\theta$ appears in $\alpha$, and label it `$e$'
otherwise.  (By \fullref{def:AOS}, $h_\theta$ cannot have
multiplicity greater than $1$ in $\alpha$.)  For agreement with the
SFT sign conventions
of Eliashberg--Givental--Hofer~\cite{eliashberg-givental-hofer00} and
Bourgeois--Mohnke \cite{bourgeois-mohnke00}, we order the
`$h$' edges by the reverse of the ordering of the $h_\theta$ orbits in
$\alpha$.  This completes the definition of the isomorphism
\eqref{eqn:GCLBar}.

We next lift \eqref{eqn:GCLBar} to an isomorphism of
$\Z[\Z^2]$--modules \eqref{eqn:GCL}. (The possible lifts according to
the prescription below will form an affine space over $\Z^2$.) To
specify a lift, first choose a reference admissible path $\Lambda_0$
of rotation number $n$ and period $\Gamma$.  Then choose the reference
cycle $\rho$ as in \eqref{eqn:referenceCycle} to be a union of Reeb
orbits corresponding to $\Lambda_0$.  By
\fullref{lem:RHRP}, this determines an isomorphism
\eqref{eqn:GCL}.

The left side of \eqref{eqn:GCL} has a well-defined relative
$\Z$--grading, by \eqref{eqn:IAmbiguity}.  We claim that this agrees
with the relative grading on the right hand side of \eqref{eqn:GCL}
defined in equation \eqref{eqn:relativeIndex}.  To see this, let
$\alpha$ and $\beta$ be generators of $\wwtilde{C}_*^{< L}(2\pi
n;\Gamma)$.  Denote the corresponding orbit sets by
$\{(\alpha_i,m_i)\}$ and $\{(\beta_j,n_j)\}$.  The correspondence of
\fullref{lem:RHRP} then defines a relative homology class
\[
Z \in H_2(T^3,\{(\alpha_i,m_i)\},\{(\beta_j,n_j)\})/H_2(T^2).
\]
We need to show that with this $Z$, the right hand sides of
\eqref{eqn:ECHIndex} and \eqref{eqn:relativeIndex} agree.  Observe
that the contact $2$--plane field $\xi$ has a nonvanishing section
$\partial_\theta$ over $T^3$, and this gives rise to a global
trivialization $\tau$ of $\xi$.  We then have
\begin{equation}
\label{eqn:ctau}
c_1(\xi|_Z,\tau) = 0.
\end{equation}
Also, if $\alpha$ or $\beta$ contains $e_\theta$ or $h_\theta$ with
multiplicity $k$, then for a sufficiently small perturbation of the
Morse--Bott contact form, $0<\phi<1/k$ in equation \eqref{eqn:CZEll}
and $r=0$ in equation \eqref{eqn:CZHyp}, so
\begin{equation}
\label{eqn:cztau}
\mu_\tau(e_\theta^k) = 1, \quad \quad \mu_\tau(h_\theta^k)=0.
\end{equation}
Therefore
\begin{equation}
\label{eqn:mutau}
\sum_i\sum_{k=1}^{m_i}\mu_\tau\left(\alpha_i^k\right) -
\sum_j\sum_{k=1}^{n_j}\mu_\tau\left(\beta_j^k\right)
=
(\ell(\alpha)-
\#h(\alpha))-(\ell(\beta)- \#h(\beta)).
\end{equation}
By equations \eqref{eqn:ctau} and \eqref{eqn:mutau}, to complete the
proof that the relative indices agree, we must show that
\begin{equation}
\label{eqn:Qtau}
Q_\tau(Z) = 2\int_{P}x\,dy.
\end{equation}
This follows as in our earlier paper \cite[Lemma~3.7]{pfh2} when
the admissible path underlying $\beta$ is obtained from that of
$\alpha$ by nondegenerate rounding.  The case of degenerate rounding
follows by an easy generalization of this.  By induction using
\fullref{prop:roundingSequence}, equation \eqref{eqn:Qtau} holds for
any two generators $\alpha$ and $\beta$.

\textbf{Step 2}\qua Choose a small perturbation of $\lambda_n$, a generic almost
complex structure $J$, and orientations of the determinant lines
$\mc{O}_{h_\theta}$ needed to define the ECH differential $\partial$ on the
left hand side of \eqref{eqn:GCL}. We assume that the perturbed
contact form agrees with $\lambda_n$ away from an
$\varepsilon$--neighborhood of the circles of Reeb orbits of
$\lambda_n$ with action $< L$, where $\varepsilon$ is small with
respect to $L$.  We claim that the differential $\partial$, regarded
as a differential on the right hand side of \eqref{eqn:GCL}, satisfies
the Nesting, Connectedness, Label Matching, and Locality axioms of
\fullref{sec:axioms} where applicable, ie whenever $\alpha$ and
$\beta$ have length $< L$.

If $\alpha$ and $\beta$ are generators of the right hand side of
\eqref{eqn:GCL}, let $\mc{M}^J(\alpha,\beta)$ denote the
$J$--holomorphic curves counted by the differential coefficient
$\langle\partial\alpha,\beta\rangle$.  To prove the Nesting axiom,
suppose there exists $C\in\mc{M}^J(\alpha,\beta)$ where $\alpha$ and
$\beta$ have length $< L$ and $\beta\not\le\alpha$.  Then there exists
$\theta_0\in\R/2\pi n\Z$ such that
\[
\det\left(\begin{matrix}\begin{matrix}
\cos \theta_0 \\
\sin \theta_0
\end{matrix}
&
\beta(\theta_0)-\alpha(\theta_0)
\end{matrix}
\right) < 0.
\]
By continuity we can choose $\theta_0$ such that $\tan\theta_0$ is not a
rational number of denominator $\le L$.  We can assume that $\varepsilon$
above is sufficiently small that the perturbed contact form agrees
with $\lambda_n$ when $\theta=\theta_0$.  We then get a contradiction
as in \fullref{prop:nesting}.

The Connectedness axiom holds because if $D(\alpha,\beta)$ is
disconnected, then as in \fullref{remark:emptySlice}, if
$\varepsilon$ is sufficiently small, then a $J$--holomorphic curve in
$\mc{M}^J(\alpha,\beta)$ has at least two non-trivial components.  By
\fullref{cor:lowIndex}, such a curve cannot exist unless
$I(\alpha,\beta)\ge 2$, whence $\langle\partial\alpha,\beta\rangle =
0$.

Before continuing, we need some restrictions on the topological complexity
of the $J$--holomorphic curves counted by $\partial$. Suppose that
$u\in\mc{M}^J(\alpha,\beta)$ and $I(\alpha,\beta)=1$. By
\fullref{cor:lowIndex}, $u$ has one component $u_1$ which does
not map to a trivial cylinder, with $\op{ind}(u_1)=1$.  Let $g(u_1)$
denote the genus of the domain of $u_1$, let $e_+(u_1)$ denote the
number of positive ends of $u_1$ at elliptic Reeb orbits, and let
$h(u_1)$ denote the number of positive or negative ends of $u_1$ at
hyperbolic Reeb orbits.  Since $\op{ind}(u_1)=1$, it follows from
equations \eqref{eqn:ind}, \eqref{eqn:ctau}, and \eqref{eqn:cztau}
that
\begin{equation}
\label{eqn:threeDollar}
2g(u_1) + 2e_+(u_1) + h(u_1) = 3.
\end{equation}
We claim now that for each $J$--holomorphic curve counted by
$\partial$, the nontrivial component $u_1$ has genus zero.  By
equation \eqref{eqn:threeDollar} the only other possibility is that
$g(u_1)=1$; $u_1$ has one positive end, which is hyperbolic; and all
negative ends of $u_1$ are elliptic.  By Nesting, $u_1$ has only one
negative end, which corresponds to the same edge as its positive end;
then $u_1\in\mc{M}^J(\alpha,\beta)$ with $I(\alpha,\beta)=-1$,
contradicting \fullref{cor:lowIndex}(a).

We now prove the Label Matching axiom.  To prove the first sentence of
the axiom, if $\langle\partial\alpha,\beta\rangle\neq 0$ and if two
edges of $\alpha$ and $\beta$ at angle $\theta_0$ agree but have
different labels, then as in \fullref{remark:emptySlice}, the
contributing $J$--holomorphic curves include nontrivial components
living in an $\varepsilon$--neighborhood of the slice
$\{\theta=\theta_0\}$.  By equation \eqref{eqn:threeDollar}, the only
such nontrivial curves that can arise are cylinders with a positive
end at $e_{\theta_0}$ and a negative end at $h_{\theta_0}$.  By
Morse--Bott theory, cf Bourgeois \cite{bourgeois04}, these count with opposite
signs as in the Morse homology of $S^1$.  The second sentence of the
Label Matching axiom holds because if $\alpha$ and $\beta$ fail this
condition, then the nontrivial component of any
$u\in\mc{M}(\alpha,\beta)$ would have a negative hyperbolic end, a
positive elliptic end, and at least one other positive end, violating
equation \eqref{eqn:threeDollar}. 

To prove the Locality axiom, let $\alpha'$ and $\beta'$ be defined as
in the statement of the axiom.  By \eqref{eqn:indexLocal} and the
analogue of \cite[Lemma~3.9]{pfh3}, taking the union with trivial
cylinders defines a map $\mc{M}(\alpha',\beta')/\R \to
\mc{M}(\alpha,\beta)/\R$, which is a bijection on the equivalence
classes of curves that the differential counts.  Our ordering
convention in Step 1 ensures that this bijection is
orientation-preserving, so $\langle\partial\alpha,\beta\rangle =
\langle\partial\alpha',\beta'\rangle$.

\textbf{Step 3}\qua The proof of \fullref{lem:vanishing} then shows that
we can write $\partial=\partial_0+\partial_1$, where
$\langle\partial_0\alpha,\beta\rangle\neq 0$ only if
$\langle\delta\alpha,\beta\rangle\neq 0$, and
$\langle\partial_1\alpha,\beta\rangle\neq 0$ only if $\beta$ is
obtained from $\alpha$ by ``double rounding'', ie rounding two
adjacent corners and losing three `$h$'s.  Also $\partial^2=0$ implies
that $\partial_0^2=0$, because $\wwtilde{C}_*^{< L}(2\pi n; \Gamma)$
is filtered by $I-\#h$, and $\partial_0$ is the differential on the
associated graded complex.

We henceforth orient all of the $\mc{O}_{h_\theta}$'s as follows.  As
mentioned above there are two $J$--holomorphic cylinders from
$e_\theta$ to $h_\theta$ which count with opposite signs.  The
projections of these cylinders to $T^2$ have areas of opposite sign.
(The areas differ by $1$.)  We choose the orientation of
$\mc{O}_{h_\theta}$ so that the cylinder whose projection to $T^2$ has
positive area counts with positive sign.

With the above orientation choices, as in \cite[Lemma~3.15(b)]{pfh3},
$\partial_0$ does not depend on the small perturbation of the contact
form, $J$, or $L$.  In conclusion, $\partial_0$ is a well-defined
differential on all of $\wwtilde{C}_*(2\pi n)$ satisfying the
Nesting, Connectedness, Label Matching, Locality, and No Double
Rounding axioms.

\textbf{Step 4}\qua We claim now that, possibly after changing some signs
in the isomorphism \eqref{eqn:GCL}, the differential $\partial_0$ also
satisfies the Degenerate Rounding and Simple Rounding axioms.  To
prove either of these axioms, by Locality we may assume that $\alpha$
has only two edges.  By the invariance of $\partial_0$, we may assume
that $J$ is close to the almost complex structure $J_{\rm std}$ defined in
equation \eqref{eqn:Jstd}.  Up to signs, the Degenerate Rounding and
Simple Rounding axioms now follow from \fullref{prop:taubes}
by using Morse--Bott theory as in \cite[Section~3.8]{pfh3}.

To understand the signs, recall from \fullref{lem:symmetry} that
$\widetilde{SL}_2\Z$ acts on $\wwbar{C}_*(2\pi n)$.  In fact, an
element $(A,f)\in\widetilde{SL}_2\Z$ gives rise to a diffeomorphism of
$T^3$ sending
\[
(\theta,(x,y))  \longmapsto (f(\theta),A(x,y))
\]
and preserving the Reeb direction, and this induces the action on the
generators of $\wwbar{C}_*(2\pi n)$.  As in \cite[Lemma~3.16]{pfh3},
with the orientation choices of Step 3, the cofficients of $\partial_0$ are
$\widetilde{SL}_2\Z$--invariant.  Thus there are only two degenerate
rounding coefficients (depending on whether the `$h$' edge comes
before or after the rounded corner) and three simple rounding
coefficients, each of which is $+1$ or $-1$.

We claim that the two degenerate rounding coefficients have opposite
signs, and the three simple rounding coefficients are related
schematically by
\begin{equation}
\label{eqn:schematically}
\langle \partial_0(eh),e\rangle =
-\langle\partial_0(he),e\rangle = \langle\partial_0(hh),h\rangle.
\end{equation}
A shortcut to checking these signs is to consider the fully twisted
chain complex with its differential $\widetilde{\partial}$, cf
\fullref{sec:fullyTwisted}.  Our previous discussion of the cylinders
from $e_\theta$ to $h_\theta$ implies that if $\Lambda$ has one edge
then
\[
\widetilde{\partial}E_\Lambda = t^k(1-t)H_\Lambda,
\]
where $t$ is a group ring generator corresponding to $H_2(T^2)$ and
$k$ is some integer depending on $\Lambda$.  (One can arrange that
$k=0$, but this is not necessary here.)  Each degenerate rounding
coefficient is now plus or minus a power of $t$.  Then applying
$\widetilde{\partial}^2=0$ to a generator $E_\Lambda$ where $\Lambda$
has two edges and a corner of angle $\pi$ implies that the two
degenerate rounding coefficients have opposite sign.  Next, applying
$\partial^2=0$ to simple triangles with two edges labeled `$h$'
establishes the relations \eqref{eqn:schematically} between the three
simple rounding coefficients.

Thus the Degenerate Rounding and Simple Rounding axioms hold up to a
global sign in each.  To make both of these signs positive, consider the
automorphism $\phi_h$ of $\wwtilde{C}_*(2\pi n)$ that sends
$\alpha\mapsto (-1)^{\#h(\alpha)}\alpha$, and similarly let $\phi_e$
be the automorphism of $\wwtilde{C}_*(2\pi n)$ that multiplies a
generator $\alpha$ by $(-1)$ to the number of elliptic orbits in
$\alpha$.  Then composing the isomorphism \eqref{eqn:GCL} with
$\phi_e$ will change the degenerate rounding sign but not the simple
rounding sign, while composing the isomorphism \eqref{eqn:GCL} with
$\phi_h$ will change both signs.

\textbf{Step 5}\qua We now complete the proof of the theorem.  By the
  previous steps, the chain complex $(\wwtilde{C}_*(2\pi
  n),\partial_0)$ satisfies all the axioms of \fullref{sec:axioms}.  By
  \fullref{prop:uniqueness}, $\partial_0=\delta$.  It follows
  from \fullref{prop:conservedQuantity} as in \cite[Lemma~A.1(a)]{pfh3}
  that for any $L$, the perturbation of $\lambda_n$ and
  $J_{\rm std}$ can be chosen so that $\partial_1=0$.  Hence
\[
\wwtilde{H}_*^{< L}(T^3,\lambda_n;\Gamma,H_2(T^2)) =
\wwtilde{H}_*^{< L}(2\pi n;\Gamma).
\]
Taking the direct limit as $L\to\infty$ proves the isomorphism
\eqref{eqn:claimTilde}.  The isomorphism \eqref{eqn:claimBar} follows
because all of the chain maps in the proof of \eqref{eqn:claimTilde}
are $\Z^2$--equivariant.
\end{proof}

\section{Concluding remarks}
\label{sec:concludingRemarks}

\subsection{Additional structure on ECH}
\label{sec:additionalStructure}

We now briefly describe some additional structures on ECH in general
and their combinatorial manifestations in the example of $T^3$,
assuming \fullref{conj:ECHDefined}.

\subsubsection{The fully twisted ECH of $T^3$}
\label{sec:fullyTwisted}

Similarly to \fullref{thm:correspondence}, the fully twisted
embedded contact homology of $T^3$ is described combinatorially by
\begin{equation}
\label{eqn:FTC}
\widetilde{ECH}_*(T^3,\lambda_n;\Gamma,0) \simeq
H_*\left(\wwtilde{C}_*(2\pi
n;\Gamma)\tensor\Z[t,t^{-1}],\widetilde{\delta}\right)
\end{equation}
where $\widetilde{\delta}$ is defined below.  First define a map
\[
\delta'\co \wwtilde{C}_*(2\pi n;\Gamma) \longrightarrow
\wwtilde{C}_{*-1}(2\pi n;\Gamma)
\]
as follows.  If $\alpha$ is a generator of $\wwtilde{C}_*(2\pi
n;\Gamma)$, define $\delta'(\alpha)$ to be the sum of all ways of
relabeling an `$e$' edge of $\alpha$ by `$h$' and making it last in
the ordering.  For example, $\delta'(E_\Lambda)=H_\Lambda$.  Note that
$\delta'$ is essentially a special case of the operator
$K_{\theta_1,\theta_2}$ defined in \fullref{sec:U}, with
$\theta_2=\theta_1+2\pi n$.  So as in
\fullref{prop:UChainHomotopy},
\begin{equation}
\label{eqn:deltadelta'}
\delta'\delta + \delta\delta'
= 0.
\end{equation}
We now define
\begin{equation}
\label{eqn:deltatilde}
\widetilde{\delta} \eqdef \delta + (1-t)\delta'.
\end{equation}
It is easy to see that $\left(\delta'\right)^2=0$.  Together with
$\delta^2=0$ and equation \eqref{eqn:deltadelta'}, this implies that
$\widetilde{\delta}^2=0$.

In the correspondence \eqref{eqn:FTC}, $t$ is an extra group ring
generator corresponding to a generator of $H_2(T^2)$.  The
$(1-t)\delta'$ term in the differential arises from the twisted Morse
complex of the circles of Reeb orbits of $\lambda_n$.

Similarly to \fullref{thm:HTilden},
\begin{equation}
\label{eqn:fullyTwisted}
H_*\left(\wwtilde{C}_*(2\pi
n;\Gamma)\tensor\Z[t,t^{-1}],\widetilde{\delta}\right)
\simeq
\left\{\begin{array}{cl} \mc{I}(\Z^3), & \Gamma=0,\; *=0,\\
\Zmodule, & \Gamma=0,\; *=1,3,\ldots,\\
0, & \mbox{otherwise.}
\end{array}\right.
\end{equation}
Here $\mc{I}(\Z^3)$ denotes the augmentation ideal in $\Z[\Z^3]$, and
$\Zmodule$ denotes the $\Z[\Z^3]$--module with one generator on which
$\Z^3$ acts by the identity.

\subsubsection{The contact element}
\label{sec:contactInvariant}

In general there is a canonical homology class
\begin{equation}
\label{eqn:CHC}
c(\lambda) \in\widetilde {ECH}_0(Y,\lambda;0,G).
\end{equation}
This is the homology class of the chain complex generator
$\alpha=\emptyset$, namely the empty set of Reeb orbits.  (This is
well-defined in the twisted cases if we choose the reference cycle
$\rho=0$ in \eqref{eqn:referenceCycle}.)  Note that $\partial\alpha=0$,
because by convexity any $J$--holomorphic curve in $\R\times Y$ has at
least one positive end.  We conjecture that the homology class
\eqref{eqn:CHC} depends only on the contact structure $\xi$.

In the untwisted ECH of $T^3$, the class $c(\lambda_n)$ corresponds to
the homology class of a $0$--gon in the combinatorial homology
$\wwbar{H}_0(2\pi n;0)$.  This homology class is a generator if
$n=1$, and $0$ if $n>1$.  On the other hand, in the fully twisted ECH
of $T^3$, the isomorphism \eqref{eqn:fullyTwisted} can be chosen so
that
\[
c(\lambda_n) = (1-t)^n \in \mc{I}(\Z^3).
\]

\subsubsection{The action of $H_1$}
\label{sec:h1action}

An element
\[
\zeta \in H_1(Y)/\op{Tors} = \Hom\left(H^1(Y;\Z),\Z\right) =
\Hom\left(H_2(Y),\Z\right)
\]
induces a degree $-1$ map
\begin{equation}
\label{eqn:degree-1Map}
\partial_\zeta\co  \widetilde{ECH}_*(Y,\lambda;\Gamma,G) \longrightarrow
\widetilde{ECH}_{*-1}(Y,\lambda;\Gamma,G).
\end{equation}
The map $\partial_\zeta$ is defined by an algebraic operation on the
fully twisted chain complex, by analogy with a construction due to
Ozsv\'ath and Szab\'o \cite{ozsvath-szabo0101}.  We will also give an
equivalent geometric definition in \fullref{sec:ECHU}.

In general, suppose we are given a free chain complex
$(\wwtilde{C}_*,\partial)$ over a group ring $\Z[H]$ and a
homomorphism $\zeta\co H\to\Z$.  Then $\zeta$ induces a $\Z$--linear map
$\widetilde{\zeta}\co \Z[H]\to\Z[H]$ sending $\sum_ha_hh\mapsto
\sum_h\zeta(h)a_hh$ and satisfying
\begin{equation}
\label{eqn:zetaProduct}
\widetilde{\zeta}(xy)=\widetilde{\zeta}(x)y + x\widetilde{\zeta}(y).
\end{equation}
Choose a basis $\{x_i\mid i\in I\}$ for $\wwtilde{C}_*$ over
$\Z[H]$, and define a $\Z[H]$--linear map
$\partial_\zeta\co  \wwtilde{C}_* \longrightarrow
\wwtilde{C}_{*-1}$
by
\[
\partial_\zeta(x_i) \eqdef
\sum_{j\in I}
\widetilde{\zeta}\big(\langle\partial x_i,x_j\rangle\big)x_j.
\]
Then $\partial^2=0$ and equation \eqref{eqn:zetaProduct} imply that
$\partial\circ\partial_\zeta + \partial_\zeta\circ\partial=0$.
Equation \eqref{eqn:zetaProduct} further implies that the map that
$\partial_\zeta$ induces on homology is natural and hence does not
depend on the choice of basis.  Also $\partial_{\zeta_1+\zeta_2} =
\partial_{\zeta_1} + \partial_{\zeta_2}$, and
$2\partial_\zeta\circ\partial_\zeta=0$ on homology.

Specializing this to the fully twisted ECH, if
$\zeta\in \Hom(H_2(Y),\Z)$, then we obtain a chain map
\[
\partial_\zeta\co  \wwtilde{C}_*(Y,\lambda;\Gamma,0) \longrightarrow
\wwtilde{C}_{*-1}(Y,\lambda;\Gamma,0).
\]
Modding out by $G$ and passing to homology gives the map
\eqref{eqn:degree-1Map}.

We now consider the example of $T^3$ and compute the map
\[
\partial_\zeta\co  ECH_*(T^3,\lambda_n;0) \longrightarrow
ECH_{*-1}(T^3,\lambda_n;0)
\]
in terms of the generators in \fullref{prop:generatingCycles}.
Recall that we have been using a basis $\{t,x,y\}$ for $H_2(T^3)$; we
denote the dual basis of $H_1(T^3)$ by the same letters $\{t,x,y\}$.
It then follows from \eqref{eqn:deltatilde} that $\partial_t$ is
induced by $-\delta'$.  Observe that $\delta'$ commutes with the
splicing chain map $S$ defined in \fullref{sec:splicing}.  It follows by
induction on $n$ that in \fullref{prop:generatingCycles} one
can take $q_{k-1}\eqdef\delta'p_k$ for $k>0$. By equation
\eqref{eqn:deltadelta'}, one can then take $u_k\eqdef\delta's_k$ and
$v_k\eqdef\delta't_k$, whence
\begin{equation}
\label{eqn:partialt}
\partial_t(s_k) = -u_k, \quad \partial_t(t_k) = -v_k
\end{equation}
for $k>0$.  It is not hard to obtain \eqref{eqn:partialt} for $k=0$ as
well.  From the bigrading and $(\delta')^2=0$, we find that
$\partial_t$ of all other generators is zero.  We also read off from
\fullref{prop:generatingCycles} that
\begin{gather*}
\partial_x(s_k) = p_k, \quad \partial_x(w_k) = -v_k,\\
\partial_y(t_k) = p_k, \quad \partial_y(w_k) = u_k,
\end{gather*}
and $\partial_x$ and $\partial_y$ of all other generators is zero.

\subsubsection{The homology operation $U$}
\label{sec:ECHU}

We now describe a degree $-2$ operation on the embedded contact homology
\[
U\co  \widetilde{ECH}_*(Y,\lambda;\Gamma,G) \longrightarrow
\widetilde{ECH}_{*-2}(Y,\lambda;\Gamma,G).
\]
Fix a point $z\in Y$ which is not on any Reeb orbit.  Let
$\mc{M}^J(\alpha,\beta,Z)^z$ denote the set of curves
$u\in \mc{M}^J(\alpha,\beta,Z)$ with a marked point mapping to
$(0,z)\in\R\times Y$.  For a suitable orientation on
$\mc{M}^J(\alpha,\beta,Z)^z$, define
\[
\begin{split}
U_z\co  \wwtilde{C}_*^{< L}(Y,\lambda;\Gamma,G) &\longrightarrow
\wwtilde{C}_{*-2}^{< L}(Y,\lambda;\Gamma,G),\\
(\alpha,[W]) & \longmapsto
\sum_{I(\alpha,\beta,Z)=2}\#\mc{M}^J(\alpha,\beta,Z)^z\cdot(\beta,[W+Z]).
\end{split}
\]
We expect to prove similarly to \fullref{conj:ECHDefined} that
$U_z$ is a chain map, and a generic path $P$ from $z$ to $z'$ induces
a chain homotopy $K_P$ between $U_{z}$ and $U_{z'}$.  The chain
homotopy counts $J$--holomorphic curves with $I=1$ that contain a
marked point mapping to $\{0\}\times P\subset\R\times Y$.  Then $U_z$
induces a well defined map $U$ on ECH.

We remark that if $P$ is a loop, then $K_P$ is equivalent to the map
$\partial_{[P]}$ defined in \fullref{sec:h1action}.

For $Y=T^3$, if we take $z=(\theta,x,y)$, then the geometric chain map
$U_z$ defined above is related to the combinatorial chain map
$U_\theta$ defined in
\fullref{sec:U} as follows.  Similarly to
\fullref{thm:correspondence}, under the isomorphism
\eqref{eqn:GCL} we have
\[
U_z = U_\theta + U',
\]
where $\langle U'\alpha,\beta\rangle\neq 0$ only if $\beta$ is
obtained from $\alpha$ by rounding two consecutive corners and losing
two `$h$'s, or rounding three consecutive corners and losing four
`$h$'s.  Without knowing anything more about the ``error term'' $U'$,
we can show that $U_z$ and $U_\theta$ induce the same map on
$ECH_*(T^3,\lambda_n;0)$.  That is, the generators in
\fullref{prop:generatingCycles} can be chosen so that
equations \eqref{eqn:U1} and \eqref{eqn:U2} hold in homology with
$U=U_\theta$ replaced by $U_z$.  We obtain the equations
\eqref{eqn:U1}, ie $U'p_{k+1}=U's_{k+1}=U't_{k+1}=0$, just by
counting the number of `$h$'s in the generators.  We then obtain the
equations \eqref{eqn:U2} by noting that $U_z$ commutes with the map
$\partial_t$ defined in \fullref{sec:h1action}, and using equation
\eqref{eqn:partialt}.

\subsection{Some other 3--manifolds}
\label{sec:otherManifolds}

\subsubsection{$S^1\times S^2$}

The methods of this paper can be modified to compute the ECH of
$S^1\times S^2$ with the contact form $\lambda_T$ studied by Taubes in
\cite{taubes02}.  Apparently
\[
\widetilde{ECH}_*(S^1\times S^2,\lambda_T;\Gamma,0) \simeq
\left\{\begin{array}{cl}\Zmodule, & \Gamma=[S^1]\times[pt],\;
*=i_0,i_0+2,\ldots,\\ 0, & \mbox{otherwise}.
\end{array}\right.
\]
Here $i_0$ is a certain odd value of the grading.  A generator in
degree $i_0$ is given (after perturbation from the Morse--Bott setting)
by a hyperbolic orbit in $S^1$ cross the equator of $S^2$.  This
calculation is relevant to Taubes's program \cite{taubes98,taubes01},
provides more evidence for \fullref{conj:big}, and shows that
ECH need not vanish for an overtwisted contact form.

\subsubsection{Torus bundles}

Let $Y$ be the $T^2$--bundle over $S^1$ with monodromy $A^{-1}\in
  SL_2\Z$.  Choose a lift $\left(A,{f}\right)\in\widetilde{SL_2\Z}$ of
  $A$, as in \fullref{sec:vanishing}, such that $f(\theta)>\theta$ for
  all $\theta$, and $f$ has rotation number in
  $(2\pi(n-1),2\pi n]$ with $n>0$.  Also choose a lift
  $\left((A^{-1})^T,{g}\right)$ of the inverse transpose $(A^{-1})^T$
  corresponding to the same $n$.  Then the diffeomorphism
\[
\begin{split}
  \R\times T^2 & \longrightarrow \R\times T^2,\\
  \left(\theta,(x,y)\right) &
  \longmapsto \left({g}(\theta),A(x,y)\right)
\end{split}
\]
preserves the contact structure given by the kernel of the standard
contact form \eqref{eqn:lambda_n}, and thus defines a contact
structure on the quotient, which is diffeomorphic to $Y$.  This
contact structure is the kernel of a $\theta$--dependent rescaling of
the contact form \eqref{eqn:lambda_n}.  The rescaling can be chosen so
that the Reeb vector field rotates to the left as $\theta$ increases,
cf \fullref{sec:taubes}.  Similarly to
\fullref{thm:correspondence}, the ECH of $Y$ for such a contact
form and for
\[
\Gamma\in\Z^2/\op{Im}(1-A)\subset H_1(Y)
\]
is computed by a ``twisted'' variant of the combinatorial
complex $\wwbar{C}_*(2\pi n;\Gamma)$.  In this chain complex, which we
denote by $\wwbar{C}_*(A,n;\Gamma)$, the periodicity condition in
\fullref{def:PAP} is replaced by the conditions
\begin{gather*}
\frac{d\Lambda}{d\theta} \circ {f} = A \circ
\frac{d\Lambda}{d\theta},\\
[\Lambda(f(\theta))-\Lambda(\theta)] = \Gamma.
\end{gather*}
It is an interesting problem to compute the homology of this complex.

\bibliographystyle{gtart}
\bibliography{link}

\begin{thebibliography}{}
\providecommand\bibmarginpar{\leavevmode\marginpar}
\def\urlstyle#1{{\tt #1}}

\bibitem{ach}
\textbf{C Abbas}, \textbf{K Cieliebak}, \textbf{H Hofer}, \emph{The {W}einstein
  conjecture for planar contact structures in dimension three}, Comment. Math.
  Helv. 80 (2005) 771--793 \xox{MR}{2182700}

\bibitem{bourgeois04}
\textbf{F Bourgeois}, \emph{A {M}orse--{B}ott approach to contact homology},
  from: ``Symplectic and contact topology: interactions and perspectives
  (Toronto, ON/Montreal, QC, 2001)'', Fields Inst. Commun. 35, Amer. Math.
  Soc., Providence, RI (2003)  55--77 \xox{MR}{1969267}

\bibitem{bourgeois-colin}
\textbf{F Bourgeois}, \textbf{V Colin},
  \href{http://dx.doi.org/10.2140/gt.2005.9.299} {\emph{Homologie de contact
  des vari\'et\'es toro\"\i dales}}, Geom. Topol. 9 (2005) 299--313
  \xox{MR}{2116317}

\bibitem{behwz}
\textbf{F Bourgeois}, \textbf{Y Eliashberg}, \textbf{H Hofer}, \textbf{K
  Wysocki}, \textbf{E Zehnder}, \href{http://dx.doi.org/10.2140/gt.2003.7.799}
  {\emph{Compactness results in symplectic field theory}}, Geom. Topol. 7
  (2003) 799--888 \xox{MR}{2026549}

\bibitem{bourgeois-mohnke00}
\textbf{F Bourgeois}, \textbf{K Mohnke},
  \href{http://dx.doi.org/10.1007/s00209-004-0656--x} {\emph{Coherent
  orientations in symplectic field theory}}, Math. Z. 248 (2004) 123--146
  \xox{MR}{2092725}

\bibitem{dragnev}
\textbf{D\,L Dragnev}, \href{http://dx.doi.org/10.1002/cpa.20018}
  {\emph{Fredholm theory and transversality for noncompact pseudoholomorphic
  maps in symplectizations}}, Comm. Pure Appl. Math. 57 (2004) 726--763
  \xox{MR}{2038115}

\bibitem{eliashberg-givental-hofer00}
\textbf{Y Eliashberg}, \textbf{A Givental}, \textbf{H Hofer},
  \emph{Introduction to symplectic field theory}, Geom. Funct. Anal.  (2000)
  560--673 \xox{MR}{1826267}

\bibitem{floer-hofer}
\textbf{A Floer}, \textbf{H Hofer}, \emph{Coherent orientations for periodic
  orbit problems in symplectic geometry}, Math. Z. 212 (1993) 13--38
  \xox{MR}{1200162}

\bibitem{pfh2}
\textbf{M Hutchings}, \href{http://dx.doi.org/10.1007/s100970100041} {\emph{An
  index inequality for embedded pseudoholomorphic curves in symplectizations}},
  J. Eur. Math. Soc. $($JEMS$)$ 4 (2002) 313--361 \xox{MR}{1941088}

\bibitem{pfh3}
\textbf{M Hutchings}, \textbf{M Sullivan},
  \href{http://dx.doi.org/10.2140/agt.2005.5.301} {\emph{The periodic {F}loer
  homology of a {D}ehn twist}}, Algebr. Geom. Topol. 5 (2005) 301--354
  \xox{MR}{2135555}

\bibitem{pfh1}
\textbf{M Hutchings}, \textbf{M Thaddeus}, \emph{Periodic Floer homology},
  $($in preparation$)$

\bibitem{kronheimer-mrowka06}
\textbf{P Kronheimer}, \textbf{T Mrowka}, \emph{Floer homology for
  Seiberg--Witten monopoles}, (in preparation)

\bibitem{kronheimer-mrowka97}
\textbf{P\,B Kronheimer}, \textbf{T\,S Mrowka},
  \href{http://dx.doi.org/10.1007/s002220050183} {\emph{Monopoles and contact
  structures}}, Invent. Math. 130 (1997) 209--255 \xox{MR}{1474156}

\bibitem{kmos}
\textbf{P Kronheimer}, \textbf{T Mrowka}, \textbf{P Ozsv\'ath}, \textbf{Z
  Szab\'o}, \emph{Monopoles and lens space surgeries}, Annals of Math. (to
  appear) \xox{arXiv}{math.GT/0310164}

\bibitem{lipshitz}
\textbf{R Lipshitz}, \emph{A cylindrical reformulation of Heegaard Floer
  homology}  \xox{arXiv}{math.SG/0502404}

\bibitem{mcduff}
\textbf{D McDuff}, \emph{Singularities and positivity of intersections of
  $J$--holomorphic curves}, from: ``Holomorphic curves in symplectic
  geometry'', Progr. Math. 117, Birkh\"auser, Basel (1994)  191--215
  \xox{MR}{1274930}

\bibitem{ozsvath-szabo0110}
\textbf{P Ozsv{\'a}th}, \textbf{Z Szab{\'o}},
  \href{http://dx.doi.org/10.1016/S0001-8708(02)00030-0} {\emph{Absolutely
  graded {F}loer homologies and intersection forms for four-manifolds with
  boundary}}, Adv. Math. 173 (2003) 179--261 \xox{MR}{1957829}

\bibitem{ozsvath-szabo0101}
\textbf{P Ozsv{\'a}th}, \textbf{Z Szab{\'o}},
  \href{http://projecteuclid.org/getRecord?id=euclid.annm/1105737568}
  {\emph{Holomorphic disks and topological invariants for closed
  three-manifolds}}, Ann. of Math. $(2)$ 159 (2004) 1027--1158
  \xox{MR}{2113019}

\bibitem{ozsvath-szabo0210}
\textbf{P Ozsv{\'a}th}, \textbf{Z Szab{\'o}},
  \href{http://projecteuclid.org/getRecord?id=euclid.dmj/1121448863}
  {\emph{Heegaard {F}loer homology and contact structures}}, Duke Math. J. 129
  (2005) 39--61 \xox{MR}{2153455}

\bibitem{parker}
\textbf{B Parker}, \emph{Holomorphic curves in Lagrangian torus fibrations},
  PhD thesis, Stanford University (2005)

\bibitem{schwarz}
\textbf{M Schwarz}, \emph{Cohomology operations from $S^1$ cobordisms in Floer
  homology}, PhD thesis, ETH Z\"urich (1995)

\bibitem{siefring}
\textbf{R Siefring}, \emph{Intersection theory of finite energy surfaces}, PhD
  thesis, New York University (2005)

\bibitem{taubes04}
\textbf{C\,H Taubes}, \emph{Pseudoholomorphic punctured spheres in
  $\mathbb{R}\times(S^1\times S^2)$: properties and existence}, $($preprint$)$

\bibitem{taubes98}
\textbf{C\,H Taubes}, \emph{The geometry of the {S}eiberg--{W}itten
  invariants}, Doc. Math.  (1998) 493--504 \xox{MR}{1648099}

\bibitem{taubes01}
\textbf{C\,H Taubes}, \emph{Seiberg--{W}itten invariants, self-dual harmonic
  2--forms and the {H}ofer--{W}ysocki--{Z}ehnder formalism}, from: ``Surveys in
  differential geometry'', Surv. Differ. Geom., VII, Int. Press, Somerville, MA
  (2000)  625--672 \xox{MR}{1919438}

\bibitem{taubes00}
\textbf{C\,H Taubes}, \emph{Seiberg {W}itten and {G}romov invariants for
  symplectic 4--manifolds}, First International Press Lecture Series 2,
  International Press, Somerville, MA (2000) \xox{MR}{1798809}

\bibitem{taubes02}
\textbf{C\,H Taubes}, \href{http://dx.doi.org/10.2140/gt.2002.6.657} {\emph{A
  compendium of pseudoholomorphic beasts in {$\Bbb R\times (S^1\times S^2)$}}},
  Geom. Topol. 6 (2002) 657--814 \xox{MR}{1943381}

\bibitem{turaev97}
\textbf{V Turaev}, \emph{Torsion invariants of $\mathrm{Spin}^c$--structures on
  3--manifolds}, Math. Res. Lett. 4 (1997) 679--695 \xox{MR}{1484699}

\bibitem{wendl}
\textbf{C Wendl}, \emph{Finite energy foliations and surgery on transverse
  links}, PhD thesis, New York University (2005)

\end{thebibliography}

\end{document}